\documentclass[11pt,twoside]{article}
\usepackage[T1]{fontenc}
\usepackage[english]{babel}
\usepackage{amsmath,amsfonts,amssymb,amsthm,fullpage}
\usepackage{dsfont}
\usepackage{fancyhdr,graphicx,color}
\usepackage[pdftex]{hyperref}
\usepackage{epsfig}
\usepackage{yhmath}

\newtheorem{thm}{Theorem}[section]
\newtheorem{lem}[thm]{Lemma}

\newtheorem{prop}[thm]{Proposition}
\newtheorem{rem}[thm]{Remark}
\newenvironment{dem}{\textbf{Proof}~:\\}{\flushright$\blacksquare$\\}
\numberwithin{equation}{section} 

\DeclareMathOperator{\card}{card}
\DeclareMathOperator{\carde}{card_{\textrm{e}}}
\DeclareMathOperator{\cardv}{card_{\textrm{v}}}
\DeclareMathOperator{\cyl}{cyl}
\DeclareMathOperator{\diam}{diam}
\DeclareMathOperator{\intt}{int}
\DeclareMathOperator{\ext}{ext}
\DeclareMathOperator{\slab}{slab}
\DeclareMathOperator{\hyp}{hyp}
\DeclareMathOperator{\flow}{flow}

\newcommand{\eps}{\varepsilon}
\newcommand{\II}{\mathds{1}}

\def\PP{\mathbb{P}}
\def\RR{\mathbb{R}}
\def\EE{\mathbb{E}}
\def\NN{\mathbb{N}}
\def\ZZ{\mathbb{Z}}
\def\QQ{\mathbb{Q}}
\def\SS{\mathbb{S}}
\def\HH{\mathbb{H}}
\def\H{\mathcal{H}}
\def\E{\mathcal{E}}
\def\C{\mathcal{C}}

\def\v{\vec{v}}
\def\w{\vec{w}}
\def\ind{{\mathds{1}}_}
\def\G{\mathcal{G}}
\def\GG{\mathfrak{G}}
\def\pe{\partial_{\textrm{e}}}
\def\pv{\partial_{\textrm{v}}}
\def\wphi{\widetilde{\phi}}
\def\wnu{\widetilde{\nu}}

\title{\huge Existence and continuity of the flow constant in first passage percolation\footnote{Research was partially supported by the ANR project PPPP (ANR-16-CE40-0016).}}
\author{Rapha\"el Rossignol\footnote{Univ. Grenoble Alpes, CNRS, Institut Fourier, F-38000 Grenoble, France} and Marie Th\'eret\footnote{LPSM UMR 8001, Universit\'e Paris Diderot, Sorbonne Paris Cit\'e, CNRS, F-75013 Paris, France. }}
\date{}

\begin{document}
\maketitle

\thispagestyle{empty}

%
%
%
%
\noindent
{\bf Abstract:}
We consider the model of i.i.d. first passage percolation on $\ZZ^d$, where we associate with the edges of the graph a family of i.i.d. random variables with common distribution $G$ on $ [0,+\infty]$ (including $+\infty$). Whereas the time constant is associated to the study of $1$-dimensional paths with minimal weight, namely geodesics, the flow constant is associated to the study of $(d-1)$-dimensional surfaces with minimal weight. In this article, we investigate the existence of the flow constant under the only hypothesis that $G(\{+\infty\} ) < p_c(d)$ (in particular without any moment assumption), the convergence of some natural maximal flows towards this constant, and the continuity of this constant with regard to the distribution $G$.\\

\noindent
{\it AMS 2010 subject classifications:} primary 60K35, secondary 82B20.

\noindent
{\it Keywords:} First passage percolation, maximal flow, minimal cutset, continuity.\\


\section{Introduction}

First passage percolation was introduced by Hammersley and Welsh \cite{HammersleyWelsh} in 1965. They defined in this model a random pseudo-metric that has been intensively studied since then. We will say a few words about it in Section \ref{s:T}, but this random metric is not the subject of this paper. The study of maximal flows in first passage percolation on $\ZZ^d$ has been initiated by Grimmett and Kesten \cite{GrimmettKesten84} in 1984 for dimension $2$ and by Kesten \cite{Kesten:flows} in 1987 for higher dimensions. This interpretation of the model of first passage percolation has been a lot less studied than the one in terms of random distance. One of the reasons is the added difficulty to deal with this interpretation, in which the study of the random paths that are the geodesics is replaced by the study of random cuts or hypersurfaces, objects which should be thought of as $(d-1)$-dimensional.

Consider a large piece of material represented by a large box in $\ZZ^d$ and equip the edges with random i.i.d. capacities representing the maximum amount of flow that each edge can bear. Typically, one is interested in the maximal flow that can cross the box from one side to the other. This question was adressed notably in \cite{Kesten:flows} and \cite{RossignolTheret08b} where one can find laws of large numbers and large deviation principles when the dimensions of the box grow to infinity. We refer to section~\ref{sec:background} for a more precise picture of the background, but let us stress for the moment that in those works, moment assumptions were made on the capacities. It is however interesting for modelling purposes to remove this assumption, allowing even infinite capacities which would represent microscopic defects where capacities are of a different order of size than elsewhere. The first achievement of the present work, Theorem~\ref{t:CV}, is to prove a law of large numbers for maximal flows without any moment assumption, allowing infinite capacities under the assumtion that the probability that an edge has infinite capacity is less than the critical parameter of percolation in $\ZZ^d$.

Once such a result is obtained, one may wonder in which way the limit obtained in this law of large numbers, the so-called flow constant, depends on the capacity distribution  put on the edges. The second achievement of this article, Theorem~\ref{thmcont}, is to show the continuity of the flow constant. One application of this continuity result could be the study of maximal flows in an inhomogeneous environment when capacities are not identically distributed but their distribution
depends smoothly (at the macroscopic scale) on the location of the edges.

The rest of the paper is organized as follows. In section~\ref{sec:background}, we give the necessary definitions and background, state our main results and explain in detail the strategy of the proof. The law of large numbers is proved in section~\ref{s:CV} and the continuity result is shown in section~\ref{s:cont}. Between those two sections lies in section~\ref{s:ssadd} a technical intermezzo 
devised to express the flow constant as the limit of a subbadditive object. The reason why we need it will be decribed at length in section~\ref{s:T}.


\section{Definitions, background and main results}
\label{sec:background}

\subsection[Maximal flows]{Definition of the maximal flows}

We use many notations introduced in \cite{Kesten:flows} and \cite{RossignolTheret08b}. Given a probability measure $G$ on $[0,+\infty]$, we equip the graph $(\ZZ^d, \EE^d)$ with an i.i.d. family $(t_G(e), e\in \EE^d)$ of random variables of common distribution $G$. Here $\EE^d$ is the set of all the edges between nearest neighbors in $\ZZ^d$ for the Euclidean distance. The variable $t_G(e)$ is interpreted as the maximal amount of water that can cross the edge $e$ per second. Consider a finite subgraph $\Omega=(V_\Omega,E_\Omega)$ of $(\ZZ^d,\EE^d)$ (or a bounded subset of $\RR^d$ that we intersect with $(\ZZ^d, \EE^d)$ to obtain a finite graph), which represents the piece of rock through which the water flows, and let $\GG^1$ and $\GG^2$ be two disjoint subsets of vertices in $\Omega$: $\GG^1$ (resp. $\GG^2$) represents the sources (resp. the sinks) through which the water can enter in (resp. escapes from) $\Omega$. A possible stream inside $\Omega$ between $\GG^1$ and $\GG^2$ is a function $\vec f : \EE^d \mapsto \RR^d$ such that for all $e\in \EE^d$,
\begin{itemize}
\item $\| \vec f(e) \|_2$ is the amount of water that flows through $e$ per second,
\item $\vec f(e) / \| \vec f(e)\|_2$ is the direction in which the water flows through $e$.
\end{itemize}
For instance, if the endpoints of $e$ are the vertices $x$ and $y$, which are at Euclidean distance $1$, then $\vec f(e) / \| \vec f(e)\|_2$ can be either the unit vector $\vec{xy}$ or the unit vector $\vec{yx}$. A stream $\vec f$ inside $\Omega$ between $\GG^1$ and $\GG^2$ is $G$-admissible if and only if it satisfies the following constraints:
\begin{itemize}
\item {\em the node law:} for every vertex $x$ in $\Omega \smallsetminus (\GG^1 \cap \GG^2)$, we have
$$\sum_{y\in \ZZ^d \,:\, e=\langle x,y\rangle \in \EE^d \cap \Omega } \| \vec f (e)\|_2 \left( \mathds{1}_{\vec f (e)/\| \vec f (e)\|_2 = \vec{xy}} -  \mathds{1}_{\vec f (e)/\| \vec f (e)\|_2 = \vec{yx}} \right) \,=\, 0 \,,$$
{\em i.e.}, there is no loss of fluid inside $\Omega$;
\item {\em the capacity constraint:} for every edge $e$ in $\Omega$, we have
$$ 0 \,\leq \, \| \vec f(e) \|_2 \,\leq \, t_G(e) \,,$$
{\em i.e.}, the amount of water that flows through $e$ per second cannot exceed its capacity $t_G(e)$.
\end{itemize}
Since the capacities are random, the set of $G$-admissible streams inside $\Omega$ between $\GG^1$ and $\GG^2$ is also random. With each such $G$-admissible stream $\vec f$, we associate its flow defined by
$$ \flow (\vec f) \,=\, \sum_{x \in \GG^1} \,\sum_{y \in \Omega \smallsetminus \GG^1 \,:\, e=\langle x,y\rangle \in \EE^d } \| \vec f (e) \|_2 \left( \mathds{1}_{\vec f (e) /\| \vec f (e)\|_2= \vec{xy}} -  \mathds{1}_{\vec f (e)/\| \vec f (e)\|_2 = \vec{yx}} \right)  \,.$$
This is the amount of water that enters in $\Omega$ through $\GG^1$ per second (we count it negatively if the water escapes from $\Omega$). By the node law, equivalently, $\flow (\vec f)$ is equal to the amount of water that escapes from $\Omega$ through $\GG^2$ per second: 
$$ \flow (\vec f) \,=\, \sum_{x \in \GG^2} \,\sum_{y \in \Omega \smallsetminus \GG^2 \,:\, e=\langle x,y\rangle \in \EE^d } \| \vec f (e) \|_2 \left( \mathds{1}_{\vec f (e)/\| \vec f (e)\|_2 = \vec{yx}} -  \mathds{1}_{\vec f (e)/\| \vec f (e)\|_2 = \vec{xy}} \right)  \,.$$
The maximal flow from $\GG^1$ to $\GG^2$ in $\Omega$ for the capacities $(t_G(e), e\in \EE^d)$, denoted by $\phi_G (\GG^1 \rightarrow \GG^2 \textrm{ in }\Omega)$, is the supremum of the flows of all admissible streams through $\Omega$:
$$ \phi_G (\GG^1 \rightarrow \GG^2 \textrm{ in }\Omega ) \,=\, \sup \{  \flow (\vec f) \,:\, \vec f \textrm{ is a $G$-admissible stream inside $\Omega$ between $\GG^1$ and $\GG^2$} \} \,.$$

It is not so easy to deal with admissible streams, but there is an alternative description of maximal flows we can work with. We define a path from $\GG^1$ to $\GG^2$ in $\Omega$ as a finite sequence $(v_0, e_1,  v_1, \dots ,e_n, v_n)$ of vertices $(v_i)_{0\leq i \leq n}$ and edges $(e_i)_{1\leq i \leq n}$ such that $v_0\in\GG^1$, $v_n\in\GG^2$ and $e_i = \langle v_{i-1},v_{i}\rangle\in E_\Omega$ for any $1\leq i \leq n$. We say that a set of edges $E \subset E_\Omega$ cuts $\GG^1$ from $\GG^2$ in $\Omega$ (or is a cutset, for short) if there is no path from $\GG^1$ to $\GG^2$ in $(V_\Omega,E_\Omega\setminus E)$. We associate with any set of edges $E$ its capacity $T_G(E)$ defined by $T_G(E) = \sum_{e\in E} t_G(e)$. The max-flow min-cut theorem (see \cite{Bollobas}), a result of graph theory, states that
$$  \phi_G (\GG^1 \rightarrow \GG^2 \textrm{ in }\Omega) \,=\, \min \{ T_G(E) \,:\, E \textrm{ cuts $\GG^1$ from $\GG^2$ in $\Omega$} \}\,.  $$
The idea of this theorem is quite intuitive: the maximal flow is limited by edges that are jammed, {\em i.e.}, that are crossed by an amount of water per second which is equal to their capacities. These jammed edges form a cutset, otherwise there would be a path of edges from $\GG^1$ to $\GG^2$ through which a higher amount of water could circulate. Finally, some of the jammed edges may not limit the flow since other edges, before or after them on the trajectory of water, already limit the flow, thus the maximal flow is given by the minimal capacity of a cutset. 

Kesten \cite{Kesten:flows} presented this interpretation of first passage percolation as a higher dimensional version of classical first passage percolation. To understand this point of view, let us associate with each edge $e$ a small plaquette $e^*$, {\em i.e.}, a $(d-1)$-dimensional hypersquare whose sides have length $1$, are parallel to the edges of the graph, such that $e^*$ is normal to $e$ and cuts $e$ in its middle. We associate with the plaquette $e^*$ the capacity $t_G (e)$ of the edge $e$ to which it corresponds. With a set of edges $E$ we associate the set of the corresponding plaquettes $E^*=\{ e^* \,:\, e \in E \}$. Roughly speaking, if $E$ cuts $\GG^1$ from $\GG^2$ in $\Omega$ then $E^*$ is a "surface" of plaquettes that disconnects $\GG^1$ from $\GG^2$ in $\Omega$ -- we do not try to give a rigorous definition of the term surface here.  In dimension $2$, the plaquette $e^*$ associated to the edge $e$ is in fact the dual edge of $e$ in the dual graph of $\ZZ^2$. A "surface" of plaquettes is thus very similar to a path in the dual graph of $\ZZ^2$ in dimension $2$. The study of maximal flows in first passage percolation is equivalent, through the max-flow min-cut theorem, to the study of the minimal capacities of cutsets. When we compare this to the classical interpretation of first passage percolation, the study of geodesics (which are paths) is replaced by the study of minimal cutsets (which are rather hypersurfaces). In this sense, the study of maximal flow is a higher dimensional version of classical first passage percolation.

We now define two specific maximal flows through cylinders that are of particular interest. Let $A$ be a non-degenerate hyperrectangle, {\em i.e.}, a rectangle of dimension $d-1$ in $\RR^d$. Let $\v$ be one of the two unit vectors normal to $A$. For a positive real $h$, denote by $\cyl (A,h)$ the cylinder of basis $A$ and height $2h$ defined by
\begin{equation}
\label{e:defcyl}
 \cyl (A,h) \,=\, \{ x + t \v \,:\, x \in A \,, \, t \in [-h, h] \} \,.
 \end{equation}
Let $B_1(A,h)$ (resp. $B_2(A,h)$) be (a discrete version of) the top (resp. the bottom) of this cylinder, more precisely defined by
\begin{align*}
B_1(A,h) & \,=\, \{ x \in \ZZ^d \cap \cyl (A,h) \,:\, \exists y  \notin \cyl(A,h) \,,\, \langle x ,y \rangle \in \EE^d \textrm{ and $\langle x ,y \rangle$ intersects } A+h \v \}\,,\\
B_2(A,h) & \,=\, \{ x \in \ZZ^d \cap \cyl (A,h) \,:\, \exists y  \notin \cyl(A,h) \,,\, \langle x,y \rangle \in \EE^d \textrm{ and $\langle x ,y \rangle$ intersects } A-h \v \}\,.
\end{align*}
We denote by $\phi_G(A,h)$ the maximal flow from the top to the bottom of the cylinder $\cyl(A,h)$ in the direction $\v$, defined by
$$ \phi_G(A,h)\,=\, \phi_G(B_1(A,h) \rightarrow B_2(A,h) \textrm{ in } \cyl(A,h) ) \,. $$
We denote by $\H^{d-1}$ the Hausdorff measure in dimension $d-1$: for $A = \prod_{i=1}^{d-1} [k_i, l_i] \times \{c\}$ with $k_i  < l_i, c \in \RR$, we have $\H^{d-1} (A) = \prod_{i=1}^{d-1} (l_i - k_i)$. We can expect that $\phi_G (A,h)$ grows asymptotically linearly in $\H^{d-1} (A)$ when the dimensions of the cylinder go to infinity, since $\H^{d-1} (A)$ is the surface of the area through which the water can enter in the cylinder or escape from it. However, $\phi_G (A,h)$ is not easy to deal with. Indeed, by the max-flow min-cut theorem, $\phi_G (A,h)$ is equal to the minimal capacity of a set of edges that cuts $B_1(A,h)$ from $B_2(A,h)$ in the cylinder. The dual of this set of edges is a surface of plaquettes whose boundary on the sides of $\cyl (A, h)$ is completely free. This implies that the union of cutsets between the top and the bottom of two adjacent cylinders is not a cutset itself between the top and the bottom of the union of the two cylinders. Thus the maximal flow $\phi_G (A,h)$ does not have a property of subadditivity, which is the key tool in the study of classical first passage percolation. This is the reason why we define another maximal flow through $\cyl (A,h)$, for which subadditivity is recovered. The set $\cyl (A,h) \smallsetminus A$ has two connected components, denoted by $C_1 (A,h)$ and $C_2 (A,h)$. For $i=1,2$, we denote by $C_i' (A,h)$ the discrete boundary of $C_i (A,h)$ defined by
\begin{equation}
\label{e:deftau1}
 C_i' (A,h) \,=\, \{ x \in \ZZ^d \cap C_i (A,h) \,:\, \exists y  \notin \cyl(A,h) \,,\, \langle x ,y \rangle \in \EE^d  \}\,. 
 \end{equation}
We denote by $\tau_G (A,h)$ the maximal flow from the upper half part of the boundary of the cylinder to its lower half part, {\em i.e.},
\begin{equation}
\label{e:deftau2}
\tau_G (A,h)\,=\, \phi_G (C'_1(A,h) \rightarrow C'_2(A,h) \textrm{ in } \cyl(A,h) ) \,.
\end{equation}
By the max-flow min-cut theorem, $\tau_G (A,h)$ is equal to the minimal capacity of a set of edges that cuts $C_1'(A,h)$ from $C'_2 (A,h)$ inside the cylinder. To such a cutset $E$ corresponds a dual set of plaquettes $E^*$ whose boundary has to be very close to $\partial A$, the boundary of the hyperrectangle $A$. We say that a cylinder is straight if $\v = \v_0 := (0, 0 , \dots , 1)$ and if there exists $k_i, l_i , c \in \ZZ$ such that $k_i<l_i$ for all $i$ and $A = A (\vec k, \vec l) = \prod_{i=1}^{d-1} [k_i, l_i] \times \{c\}$. In this case, for $c=0$ and $k_i \leq 0 < l_i$, the family of variables $(\tau_G ( A(\vec k , \vec l), h) )_{\vec k, \vec l}$ is subadditive, since the minimal cutsets in adjacent cylinders can be glued together along the common side of these cylinders.


\subsection[Background]{Background on maximal flows}

A straightforward application of ergodic subadditive theorems in the multiparameter case (see Krengel and Pyke \cite{KrengelPyke} and Smythe \cite{Smythe}) leads to the following result.
\begin{prop}
Let $G$ be a probability measure on $[0,+\infty[$ such that $\int_{\RR^+} x \, dG(x) <\infty$. Let $A = \prod_{i=1}^{d-1} [k_i, l_i] \times \{0\}$ with $k_i \leq 0 < l_i \in \ZZ$. Let $h: \NN \rightarrow \RR^+$ such that $\lim_{p\rightarrow \infty} h(p) = +\infty$. Then there exists a constant $\nu_G (\v_0)$, that does not depend on $A$ and $h$, such that 
$$ \lim_{p\rightarrow \infty} \frac{\tau_G (pA, h(p))}{  \H^{d-1} (pA) } \,=\, \nu_G (\v_0) \quad \textrm{a.s. and in }L^1\,. $$
\end{prop}
This result has been stated in a slightly different way by Kesten in \cite{Kesten:flows}. He considered there the more general case of flows through cylinders whose dimensions goes to infinity at different speeds in each direction, but in dimension $d=3$. The constant $\nu_G (\v_0)$ obtained here is the equivalent of the time constant $\mu_G(e_1)$ defined in the context of random distances (see Section \ref{s:T}), and by analogy we call it the flow constant.

As suggested by classical first passage percolation, a constant $\nu_G (\v)$ can be defined in any direction $\v \in \SS^{d-1}$, where $\SS^{d-1} =\{x\in \RR^d \,:\, \|x\|_2 = 1  \}$. This is not that trivial, since a lack of subadditivity appears when we look at tilted cylinders, due to the discretization of the boundary of the cylinders. Moreover, classical ergodic subadditive theorems cannot be used if the direction $\v$ is not rational, {\em i.e.}, if there does not exist an integer $M$ such that $M \v$ has integer coordinates. However, these obstacles can be overcome and the two authors proved in \cite{RossignolTheret08b} the following law of large numbers.
\begin{thm}
\label{t:oldCV}
Let $G$ be a probability measure on $[0, +\infty[$ such that $\int_{\RR^+} x \, dG(x) <\infty$. For any $\v \in \SS^{d-1}$, for any non-degenerate hyperrectangle $A$ normal to $\v$, for any function $h : \NN \mapsto \RR^+$ satisfying $\lim_{p\rightarrow +\infty} h(p) =+\infty$, there exists a constant $\nu_G (\v) \in [0,+\infty[$ (independent on $A$ and $h$) such that 
$$ \lim_{p\rightarrow \infty} \frac{\tau_G (pA, h(p))}{\H^{d-1} (pA)} \,=\, \nu_G (\v) \quad \textrm{in }L^1\,. $$
If moreover the origin of the graph belongs to $A$, or if $\int_{\RR^+} x^{1+1/(d-1)} \, dG(x) <\infty$, then 
$$ \lim_{p\rightarrow \infty} \frac{\tau_G (pA, h(p))}{\H^{d-1} (pA)} \,=\, \nu_G (\v) \quad \textrm{a.s.} $$
If the cylinder is {\em flat}, {\em i.e.}, if $\lim_{p\rightarrow \infty}  h(p) /p =0$, then the same convergences hold also for $\phi_G (pA, h(p))$.
\end{thm}
When the origin of the graph belongs to $A$, and for an increasing function $h$ for instance, the cylinder $\cyl (pA, h(p))$ is completely included in the cylinder $\cyl ((p+1)A, h(p+1))$. The mean of the capacities of the edges inside $\cyl (pA, h(p))$ converges a.s. when $p$ goes to infinity as soon as $\int_{\RR^+} x \, dG(x) <\infty$ by a simple application of the law of large numbers, and Theorem \ref{t:oldCV} states that $\tau_G (pA, h(p))/\H^{d-1} (pA)$ converges a.s. under the same hypothesis. On the other hand, when the origin of the graph does not belong to $A$, the cylinders $\cyl (pA, h(p))$ and $\cyl ((p+1)A, h(p+1))$ may be completely disjoint. The a.s. convergence of the mean of the capacities of the edges included in $\cyl (pA, h(p))$ when $p$ goes to infinity is thus stated by some result about complete convergence of arrays of random variables, see for instance \cite{Gut92,Gut85}. This kind of results requires a stronger moment condition on the law of the random variables we consider, namely we need that $\int_{\RR^+} x^{1+1/(d-1)} \, dG(x) <\infty$. Theorem \ref{t:oldCV} states that $\tau_G (pA, h(p))/\H^{d-1} (pA)$ converges a.s. in this case under the same hypothesis on the moments of $G$.

Let $p_c(d)$ be the critical parameter of Bernoulli bond percolation on $(\ZZ^d, \EE^d)$. Zhang investigated in \cite{Zhang} the positivity of $\nu_G$ and proved the following result.
\begin{thm}
\label{t:oldnul}
Let $G$ be a probability measure on $[0, +\infty[$ such that $\int_{\RR^+} x \, dG(x) <\infty$. Then 
$$ \nu_G (\v) \,>\,0 \quad \iff \quad G(\{0\}) \,<\, 1-p_c(d) \,. $$
\end{thm}

The asymptotic behavior of the maximal flows $\phi_G (pA, h(p))$ in non-flat cylinders ({\em i.e.}, when $h(p)$ is not negligible in comparaison with $p$) is more difficult to study since these flows are not subadditive. In the case of straight cylinders (and even in a non isotropic case, {\em i.e.}, when the dimensions of the cylinders go to infinity at different speed in every directions), Kesten \cite{Kesten:flows} and Zhang \cite{Zhang2017} proved that $\phi_G (pA, h(p)) / \H^{d-1} (pA)$ converges a.s. towards $\nu_G (\v_0)$ also, under some moment condition on $G$. The behavior of $\phi_G (pA, h(p))$ is different in tilted and non-flat cylinders, we do not go into details and refer to \cite{RossignolTheret09b} (for $d=2$) and to \cite{CerfTheret09geoc,CerfTheret09supc,CerfTheret09infc,CerfTheretStream} in a more general setting.

We stress the fact that for all the results mentioned above, a moment assumption is required on the probability measure $G$ on $[0,+\infty[$: $G$ must at least have a finite mean.


\subsection{Main results}

Our first goal is to extend the previous results to probability measures $G$ on $[0,+\infty[$ that are not integrable, and even to probability measures $G$ on $[0,+\infty]$ under the hypothesis that $G(\{+\infty\}) < p_c(d)$.

For any probability measure $G$ on $[0, +\infty]$, for all $K>0$, we define $G^K = \ind{[0,K[} G + G([K,+\infty[) \delta_K$, {\em i.e.}, $G^K$ is the law of $\min(t_G(e), K)$ for any edge $e$. Then we define
\begin{equation}
\label{defnu}
   \forall \v \in \SS^{d-1} \,, \quad  \nu_G (\v) \,:=\, \lim_{K \rightarrow \infty} \nu_{G^K} (\v) \,. 
\end{equation}
Throughout the paper, we shall say that a function $h : \NN \mapsto \RR^+$ is \emph{mild} if 
\begin{equation}
\lim_{p\rightarrow +\infty} h(p) / \log p =+\infty\text{ and }\lim_{p\rightarrow \infty}  h(p) /p =0\;.
\end{equation}
We prove the following law of large numbers for cylinders with mild height functions.
\begin{thm}
\label{t:CV}
For any probability measure $G$ on $[0,+\infty]$ such that $G(\{+\infty\}) < p_c(d)$, for any $\v \in \SS^{d-1}$, for any non-degenerate hyperrectangle $A$ normal to $\v$, for any mild function $h$, we have
$$  \lim_{p\rightarrow \infty} \frac{\phi_G (pA, h(p))}{\H^{d-1} (pA)} \,=\, \nu_G (\v) \quad\textrm{a.s.}  $$
Moreover, for every $\v \in \SS^{d-1}$, 
$$\nu_G (\v) < +\infty$$
and 
$$ \nu_G (\v) \,>\,0 \quad \iff \quad G(\{0\}) \,<\, 1-p_c(d) \,. $$
\end{thm}
\begin{rem}
If $G$ is integrable, the constant $\nu_G$ defined in \eqref{defnu} is thus coherent with the definition given by Theorem \ref{t:oldCV}.
\end{rem}

We also want to establish the continuity of the function $G \mapsto \nu_G (\v)$ when we equip the set of probability measures on $[0,+\infty]$ with the topology of weak convergence - in fact these two questions are linked, as we will see in Section \ref{s:T}. More precisely, let $(G_n)_{n\in \NN}$ and $G$ be probability measures on $[0,+\infty]$. We say that $G_n$ converges weakly towards $G$ when $n$ goes to infinity, and we write $G_n \overset{d}{\rightarrow} G$, if for any continuous bounded function $f: [0,+\infty] \mapsto \RR^+$ we have
$$ \lim_{n \rightarrow +\infty} \int_{[0,+\infty]} f \, dG_n \,=\,  \int_{[0,+\infty]} f \, dG \,. $$
Equivalently, $G_n \overset{d}{\rightarrow} G$ if and only if $\lim_{n \rightarrow \infty} G_n([t,+\infty]) = G([t,+\infty])$ for all $t\in \RR^+$ such that $t\mapsto G([t,+\infty])$ is continuous at $t$. 
\begin{thm}
\label{thmcont}
Suppose that $G$ and $(G_n)_{n\in \NN}$ are probability measures on $[0,+\infty]$ such that $G(\{+\infty\}) < p_c(d)$ and for all $n\in \NN$, $G_n (\{+\infty\}) < p_c(d)$. If $G_n \overset{d}{\rightarrow} G$, then
$$   \lim_{n\rightarrow \infty} \sup_{\v \in \SS^{d-1}} \left\vert \nu_{G_n} (\vec v) - \nu_G (\vec v) \right\vert \,=\, 0 \,.$$
\end{thm}


\subsection[Time constant]{About the existence and the continuity of the time constant}
\label{s:T}

First passage percolation was introduced by Hammersley and Welsh \cite{HammersleyWelsh} in 1965 with a different interpretation of the variables associated with the edges. We consider the graph $(\ZZ^d, \EE^d)$ and we associate with the edges of the graph a family of i.i.d. random variables $(t_G(e), e\in \EE^d)$ with common distribution $G$ as previously, but we interpret now the variable $t_G(e)$ as the time needed to cross the edge $e$ (we call it the passage time of $e$). If $\gamma$ is a path, we define the passage time of $\gamma$ as $T_G(\gamma) = \sum_{e\in \gamma} t_G(e)$. Then the passage time between two points $x$ and $y$ in $\ZZ^d$, {\em i.e.}, the minimum time needed to go from $x$ to $y$ for the passage times $(t_G(e), e\in \EE^d)$, is given by
$$ T_G(x,y) \,=\, \inf \{ T_G(\gamma) \,:\, \gamma \textrm{ is a path from $x$ to $y$} \} \,.$$
This defines a random pseudo-distance on $\ZZ^d$ (the only property that can be missing is the separation property). This random distance has been and is still intensively studied. A reference work is Kesten's lecture notes \cite{Kesten:StFlour}. Auffinger, Damron and Hanson wrote very recently the survey \cite{AuffingerDamronHanson} that provides an overview on results obtained in the 80's and 90's, describes the recent advances and gives a collection of old and new open questions.

Fix $e_1 = (1,0, \dots , 0)$. Thanks to a subadditive argument, Hammersley and Welsh \cite{HammersleyWelsh} and Kingman \cite{Kingman68} proved that if $d=2$ and $F$ has finite mean, then $\lim_{n\rightarrow \infty} T_F(0,ne_1) / n $ exists a.s. and in $L^1$, the limit is a constant denoted by $\mu_F(e_1) $ and called the time constant. The moment condition was improved some years later by several people independently, and the study was extended to any dimension $d\geq 2$ (see for instance Kesten's Saint-Flour lecture notes \cite{Kesten:StFlour}). The convergence to the time constant can be stated as follows.
\begin{thm}
If $\EE[\min (t_F(1),\dots , t_{F}(2d))] < \infty$ where $(t_F(i), i\in \{1,\dots , 2d\})$ are i.i.d. with distribution $F$ on $[0,+\infty[$, there exists a constant $\mu_F(e_1) \in \RR^+$ such that
$$  \lim_{n\rightarrow \infty} \frac{T_F(0,n e_1)}{ n} \,=\, \mu_F(e_1) \quad \textrm{a.s. and in } L^1\,.$$
Moreover, the condition $\EE[\min (t_F(1),\dots , t_{F}(2d))] < \infty$ is necessary for this convergence to hold a.s. or in $L^1$.
\end{thm}
This convergence can be generalized by the same arguments, and under the same hypothesis, to rational directions : there exists a homogeneous function $\mu_F : \QQ^d \rightarrow \RR^+$ such that for all $x \in \ZZ^d$, we have $\lim_{n\rightarrow \infty} T_F(0, nx) / n = \mu_F (x)$ a.s. and in $L^1$. The function $\mu_F$ can be extended to $\RR^d$ by continuity (see \cite{Kesten:StFlour}). 

These results can be extended by considering a law $F$ on $[0,+\infty[$ which does not satisfy any moment condition, at the price of obtaining weaker convergence. This work was performed successfully by Cox and Durrett \cite{CoxDurrett} in dimension $d=2$ and then by Kesten \cite{Kesten:StFlour} in any dimension $d\geq 2$. More precisely, they proved that there always exists a function $\hat \mu_F : \RR^d \rightarrow \RR^+$ such that for all $x \in \ZZ^d$, we have $\lim_{n\rightarrow \infty} T_F(0, nx) / n =\hat \mu_F (x)$ in probability. If $\EE[\min (t_F(1),\dots , t_{F}(2d))] < \infty$ then $\hat \mu_F = \mu_F$. The function $\hat \mu_F$ is built as the a.s. limit of a more regular sequence of times $\hat T_F (0,nx) /n$ that we now describe roughly. They consider an $M\in \RR^+$ large enough so that $ F( [0,M])$ is very close to $1$. Thus the percolation $(\ind{\{ t_F(e) \leq M \}}, e \in \EE^d)$ is highly supercritical, so if we denote by $\C_{F,M}$ its infinite cluster, each point $x\in \ZZ^d$ is a.s. surrounded by a small contour $S(x) \subset \C_{F,M}$. They define $\hat T_F (x,y) = T_F(S(x), S(y))$ for $x,y \in \ZZ^d$. The times $\hat T _F(0,x)$ have good moment properties, thus $\hat \mu_F(x)$ can be defined as the a.s. and $L^1$ limit of $\hat T_F (0,nx) /n$ for all $x\in \ZZ^d$ by a classical subadditive argument; then $\hat \mu_F$ can be extended to $\QQ^d$ by homogeneity, and finally to $\RR^d$ by continuity. The convergence of $T_F(0, nx)/n$ towards $\hat \mu_F (x)$ in probability is a consequence of the fact that $T_F$ and $\hat T_F$ are close enough.

It is even possible to consider a probability measure $F$ on $[0,+\infty]$ under the hypothesis that $F([0,+\infty[) > p_c(d)$. This was done first by Garet and Marchand in \cite{GaretMarchand04} and then by Cerf and the second author in \cite{CerfTheretForme}. We concentrate on \cite{CerfTheretForme}, where the setting is closer to the one we consider here. To prove the existence of a time constant for a probability measure $F$ on $[0,+\infty]$ such that $F([0,+\infty[) > p_c(d)$, Cerf and the second author exhibit a quite intuitive object that is still subadditive. For $x\in \ZZ^d$, $\tilde \mu_F (x)$ is defined by a subadditive argument as the limit of $T_F(f_M(0), f_M(nx))/n$ a.s. and in $L^1$, where $M$ is a real number large enough such that $F([0,M]) >p_c(d)$, and for $z\in \ZZ^d$, $f_M(z)$ is the points of $\C_{F,M}$ which is the closest to $z$. The convergence of $T_F(0, nx)/n$ towards $\tilde \mu_F (x)$ still holds, but in a very weak sense: $T_F(0, nx)/n$ converges in fact in distribution towards $\theta_F^2 \delta_{\mu_F (x)} + (1-\theta_F^2) \delta_{+\infty}$, where $\theta_F$ is the probability that the connected component of $0$ in the percolation $(\mathds{1}_{t_F(e) < \infty}, e\in \EE^d)$ is infinite. For short, all these constants ($\hat \mu_F, \tilde \mu_F$ and $\mu_F$) being equal when they are defined, we denote all of them by $\mu_F$.

Once the time constant is defined, a natural question is to wonder if it varies continuously with the distribution of the passage times of the edges. This question has been answered positively by Cox and Kesten \cite{Cox,CoxKesten,Kesten:StFlour} for probability measures on $[0,+\infty[$.
\begin{thm}
Let $F$, $F_n$ be probability measures on $[0,+\infty[$. If $F_n$ converges weakly towards $F$, then for every $x\in \RR^d$,
$$\lim_{n\rightarrow \infty} \mu_{F_n}  (x) \,=\, \mu_F (x) \,.$$
\end{thm}
Cox \cite{Cox} proved first this result in dimension $d=2$ with an additional hypothesis of uniform integrability: he supposed that all the probability measures $F_n$ were stochastically dominated by a probability measure $H$ with finite mean. To remove this hypothesis of uniform integrability in dimension $d=2$, Cox and Kesten \cite{CoxKesten} used the regularized passage times and the technology of the contours introduced by Cox and Durrett \cite{CoxDurrett}. Kesten \cite{Kesten:StFlour} extended these results to any dimension $d\geq 2$. The key step of their proofs is the following lemma.
\begin{lem}
\label{l:key}
Let $F$ be a probability measure on $\RR^+$, and let $F^K = \mathds{1}_{[0,K)} F + F([K,+\infty)) \delta_K$ be the distribution of the passage times $t_F(e)$ truncated at $K$. Then for every $x\in \RR^d$,
$$ \lim_{K \rightarrow \infty} \mu_{F^K} (x) \,=\, \mu_F (x) \,. $$
\end{lem}
To prove this lemma, they consider a geodesic $\gamma$ from $0$ to a fixed vertex $x $ for the truncated passage times $\inf (t_{F}(e), K)$. When looking at the original passage times $t_F (e)$, some edges along $\gamma$ may have an arbitrarily large passage time: to recover a path $\gamma'$ from $0$ to $x$ such that $T_F (\gamma')$ is not too large in comparison with $T_{F^K} (\gamma)$, they need to bypass these bad edges. They construct the bypass of a bad edge $e$ inside the contour $S(e) \subset \C_{F,M}$ of the edge $e$, thus they bound the passage time of this bypass by $M \carde (S(e))$ where $\carde (S(e))$ denotes the number of edges in $S(e)$. More recently, Garet, Marchand, Procaccia and the second author extended in \cite{GaretMarchandProcacciaTheret} these results to the case where the probability measures considered are defined on $[0,+\infty]$ as soon as the percolation of edges with finite passage times are supercritical. To this end, they needed to perform a rescaling argument, since for $M$ large enough the percolation of edges with passage times smaller than $M$ can be choosen supercritical but not highly supercritical as required to use the technology of the contours.

The study of the existence of the time constant without any moment condition and the study of the continuity of the time constant with regard to the distribution of the passage times of the edges are closely related. Indeed, in the given proofs of the continuity of the time constant, the following results are used:
\begin{itemize}
\item the time constant $\mu_F$ is the a.s. limit of a subadditive process,
\item this subadditive process is integrable (for any distribution $F$ of the passage times, even with infinite mean),
\item this subadditive process is monotonic with regard to the distribution of the passage times.
\end{itemize}
Moreover, the technology used to prove the key Lemma \ref{l:key} (using the contours) is directly inspired by the study of the existence of the time constant without any moment condition.

The proof of the continuity of the flow constant, Theorem \ref{thmcont}, we propose in this paper is heavily influenced by the proofs of the continuity of the time constant given in \cite{CoxKesten,Kesten:StFlour,GaretMarchandProcacciaTheret}. The real difficulty of our work is to extend the definition of the flow constant to probability measure with infinite mean - once this is done, it is harmless to admit probability measures $F$ on $[0,+\infty]$ such that $F(\{+\infty\}) <p_c(d)$, we do not even have to use a renormalization argument. We choose to define the flow constant $\nu_F$ via \eqref{defnu} so that the result equivalent to Lemma \ref{l:key} in our setting is given by the precise definition of $\nu_F$. However, two major issues remain :
\begin{itemize}
\item[$(i)$] prove that $\nu_F$ is indeed the limit of some quite natural sequence of maximal flows,
\item[$(ii)$] prove that $\nu_F$ can be recovered as the limit of a nice subadditive process.
\end{itemize}

The first point, $(i)$, is precisely the object of Theorem \ref{t:CV}, that we prove in Section \ref{s:CV}. With no surprise, the difficulties we do not meet to prove the result equivalent to Lemma \ref{l:key} for the flow constant are found in the proof of this convergence, see Proposition \ref{p:tronquer}. The maximal flows that converge towards $\nu_G$ are maybe the most natural ones, {\em i.e.}, maximal flows from the top to the bottom of flat cylinders, and the convergence holds a.s., {\em i.e.}, in a strong sense, which is quite satisfying. It is worth noticing that in fact, to prove the a.s. convergence in tilted cylinders when $\nu_F=0$  (see Proposition \ref{p:zeroter}), we use the continuity of the flow constant - without this property, we obtain only a convergence in probability. However, to obtain a convergence (at least in probability) of these maximal flows towards $\nu_F$, we do not have to exhibit a subadditive process converging towards $\nu_F$. The existence of such a nice subadditive process, {\em i.e.}, the point $(ii)$ above, is nevertheless needed to prove the continuity of the flow constant. In Section \ref{s:ssadd}, we define such a process and prove its convergence towards $\nu_F$ (see Theorem \ref{t:ssadd}). Finally in Section \ref{s:cont} we prove the continuity of the flow constant, Theorem \ref{thmcont}.

Before starting these proofs, we give in the next section some additional notations.


\subsection{More notations}

We need to introduce a few more notations that will be useful. 

Given a unit vector $\v \in \SS^{d-1}$ and a non-degenerate hyperrectangle $A$ normal to $\v$, $\hyp (A)$ denotes the hyperplane spanned by $A$ defined by
$$ \hyp (A) \,=\, \{ x+ \vec w \,:\, x\in A \,,\, \vec w \cdot \v =0 \} $$
where $\cdot$ denotes the usual scalar product on $\RR^d$. For a positive real $h$, we already defined $\cyl (A,h)$ as the cylinder of height $2h$ with base $A-h\v$ and top $A+h\v$, see Equation \eqref{e:defcyl}. It will sometimes be useful to consider the cylinder $\cyl^{\v} (A,h)$ with height $h$, base $A$ and top $A+h\v$, {\em i.e.},
$$ \cyl^{\v} (A,h) \,=\, \{ x+t \v \,:\, x\in A\,,\, t\in [0,h] \} \,,$$
and the maximal flow $\phi_G^{\v} (A,h)$ from the discrete version of its top
$$ B_1^{\v}(A,h)  \,=\, \{ x \in \ZZ^d \cap \cyl^{\v} (A,h) \,:\, \exists y  \notin \cyl^{\v}(A,h) \,,\, \langle x ,y \rangle \in \EE^d \textrm{ and $\langle x ,y \rangle$ intersects } A+h \v \} $$
to the discrete version of its bottom
$$ B^{\v}_2(A,h)  \,=\, \{ x \in \ZZ^d \cap \cyl^{\v} (A,h) \,:\, \exists y  \notin \cyl^{\v}(A,h) \,,\, \langle x,y \rangle \in \EE^d \textrm{ and $\langle x ,y \rangle$ intersects } A-h \v \}\,.$$

Some sets can be seen as sets of edges or vertices, thus when looking at their cardinality it is convenient to specify whether we count the number of edges or the number of vertices in the set. The notation $\carde (\cdot)$ denotes the number of edges in a set whereas $\cardv (\cdot)$ denotes the number of vertices.

Given a probability measure $G$ on $[0,+\infty]$, a constant $K\in ]0,+\infty[$ and a vertex $x\in \ZZ^d$ (respectively an edge $f\in \EE^d$), we denote by $C_{G,K} (x)$ (resp. $C_{G,K} (f)$) the connected component of $x$ (resp. the union of the connected components of the two endpoints of $f$) in the percolation $(\mathds{1}_{t_G(e) > K}, e\in \EE^d)$, which can be seen as an edge set and as a vertex set. For any vertex set $C \subset \ZZ^d$, we denote by $\diam (C)$ the diameter of $C$, $\diam (C) = \sup\{ \|x-y\|_2 \,:\, x,y\in C \cap \ZZ^d \}$, by $\pe C$ its exterior edge boundary defined by
$$ \pe C \,=\, \{ e = \langle x,y \rangle \in \EE^d \,:\, x\in C \,,\, y\notin C \textrm{ and there exists a path from $y$ to infinity in } \ZZ^d \smallsetminus C\}\,, $$
and by $\pv C$  its exterior vertex boundary defined by
$$ \pv C \,=\, \{ x\in \ZZ^d \,:\, x\notin C \,,\, \exists y\in C \textrm{ s.t. } \langle x,y \rangle \in \EE^d\}\,. $$
Given a set $E$ of edges, we can define also its diameter $\diam (E)$ as the diameter of the vertex set made of the endpoints of the edges of $E$. We also define its exterior $\ext (E)$ by
$$ \ext (E) \,=\, \{ x\in \ZZ^d \,:\, \textrm{there exists a path from $x$ to infinity in }\EE^d \smallsetminus E \} $$
and its interior 
$$ \intt (E) \,=\, \ZZ^d \smallsetminus \ext (E)\,. $$
Notice that by definition, $C \subset \intt (\pe C)$ and if $C$ is bounded and $x \in \intt (\pe C)$, then $ \pe C$ separates $x$ from infinity. For any vertices $x$ and $y$, for any probability measure $G$ on $[0,+\infty]$ and any $K\in ]0,+\infty]$, one of the three following situation occurs:
\begin{itemize}
\item[$(i)$] $\pe C_{G,K} (x)=\pe C_{G,K} (y)$;
\item[$(ii)$] $\intt (\pe C_{G,K} (x)) \cap \intt (\pe C_{G,K} (y)) = \emptyset $;
\item[$(iii)$] $ \intt  (\pe C_{G,K} (x)) \subset \intt (\pe C_{G,K} (y))$, or $ \intt  (\pe  C_{G,K} (y) )\subset \intt (\pe C_{G,K} (x))$.
\end{itemize}
Case $(i)$ corresponds to the case where $x$ and $y$ are connected in the percolation $(\mathds{1}_{t_G(e) > K}, e\in \EE^d)$, whereas cases $(ii)$ and $(iii)$ correspond to the case where $x$ an $y$ are not connected, thus their connected components for this percolation are disjoint. Case $(iii)$ corresponds to the case where $x\in \intt( \pe C_{G,K} (y))$ (thus $C_{G,K} (x)$ is nested in a hole of $C_{G,K} (y)$ inside $\intt (\pe C_{G,K} (y))$) or conversely, whereas case $(ii)$ corresponds to the case where $x\in \ext( \pe C_{G,K} (y))$ and $y\in \ext (\pe C_{G,K} (x))$.

For any subset $\C$ of $\RR^d$ and any $h\in \RR^+$, we denote by $\E_{G,K} (\C, h)$ the following event
\begin{equation}
\label{e:E}
\E_{G,K} (\C, h) \,=\, \bigcap_{x\in \C \cap \ZZ^d} \{ \diam (C_{G,K} (x)) < h \}\,,
\end{equation}
and by $\E'_{G,K} (\C, h)$ the corresponding event involving edges instead of vertices
\begin{equation}
\label{e:E'}
\E'_{G,K} (\C, h) \,=\, \bigcap_{e\in \C \cap \EE^d} \{ \diam (C_{G,K} (e)) < h \}\,.
\end{equation}

In what follows $c_d$ denotes a constant that depends only on the dimension $d$ and may change from one line to another. Notice that for any finite and connected set $C$ of vertices, $\carde (\pe C) \leq c_d \cardv (C)$.

For two probability measures $H$ and $G$ on $[0,+\infty]$, we define the following stochastic domination relation:
\begin{equation}
\label{e:defdomsto}
G\preceq H \quad \iff \quad \forall t\in [0,+\infty) \quad G([t,+\infty]) \leq H([t,+\infty])\,.
\end{equation}
In what follows, we always build the capacities of the edges for different distributions by coupling, using a family of i.i.d. random variables with uniform distribution on $]0,1[$ and the pseudo-inverse of the distribution function of these distributions. Thus the stochastic comparison between probability measures $H$ and $G$ on $[0,+\infty]$ implies a simple comparison between the corresponding capacities of the edges:
\begin{equation}
\label{e:couplage}
 G\preceq H \quad \Longrightarrow  \forall e\in \EE^d \quad t_G(e) \leq t_H (e) \,. 
 \end{equation}


\section{Convergence of the maximal flows}
\label{s:CV}

This section is devoted to the proof of Theorem \ref{t:CV}.

\subsection{Properties of $\nu_G$}

First we investigate the positivity $\nu_G$ as defined by \eqref{defnu}.
\begin{prop}
\label{p:trivial}
Let $G$ be a probability measure on $[0, +\infty]$ such that $G(\{+\infty\})<p_c(d)$. For every $\v \in \SS^{d-1}$, we have
$$ \nu_G (\v) = 0 \quad \iff \quad G(\{0\})\geq 1-p_c(d)\,. $$
\end{prop}
\begin{dem}
By the above coupling, see Equation \eqref{e:couplage}, for any such probability $G$, for any $0<K_1\leq K_2$, for any $\v\in \SS^{d-1}$, for any hyperrectangle $A$ and any $h\in \RR^+$, we have
$$ \phi_{G^{K_1}} (A,h) \,\leq \,  \phi_{G^{K_2}} (A,h) \,.$$
By definition of $\nu_{G^K} (\v)$ (see Theorem \ref{t:oldCV}), this proves that $K\mapsto \nu_{G^K} (\v)$ is non-decreasing. Thus $\nu_G (\v) = 0$ if and only if for every $K\in \RR^+$, $\nu_{G^K} (\v) = 0$. By Theorem \ref{t:oldnul}, we know that $\nu_{G^K} (\v) = 0$ if and only if $G^K (\{0\}) \geq 1-p_c(d)$. But $G^K (\{0\})=G (\{0\})$ for all $K$, thus Proposition \ref{p:trivial} is proved.
\end{dem}

We now state a stochastic domination result, in the spirit of Fontes and Newman \cite{FontesNewman}, which will be useful to prove that $\nu_G$ is finite, and will be used again in section~\ref{subsec:truncating}.
\begin{lem}
\label{lem:domsto}
Let $W=\{x_1,\ldots,x_n\}$ be a finite subset of $\ZZ^d$. Consider an i.i.d. Bernoulli bond percolation on $\ZZ^d$. For $i=1,\ldots,n$, define $Z_i = Z(x_i)$ to be $\cardv (C(x_i))$, where $C(x_i)$ is the connected component of $x_i$ for the underlying percolation. Let $Y_1=Z_1$ and define recursively $Y_i$ for $i=2,\ldots,n$ by
$$Y_{i}=\left\lbrace\begin{array}{ll}Z_{i}&\text{ if }x_i\not\in \cup_{j=1}^{i-1}C(x_j)\\0 & \text{ if }x_i\in \cup_{j=1}^{i-1}C(x_j) \,.\end{array}\right.$$
Let also $(X_i, i\in \{1,\ldots,n\})$ be a family of i.i.d. random variables distributed as $Z_1 = \cardv (C(x_1))$. Then, for all $a$, $a_1,\ldots,a_n$ in $\RR$,
$$\PP\left(\sum_{i=1}^nY_i\geq a\text{ and }\forall i=1,\ldots ,n,\;Y_i\geq a_i\right)\leq \PP\left(\sum_{i=1}^n X_i\geq a\text{ and }\forall i=1,\ldots ,n,\;X_i\geq a_i \right)\,.$$
\end{lem}
\begin{dem}
For any $i$, let $\mathcal{F}_i$ be the sigma-field generated by the successive exploration of $C(x_1)$, $C(x_2)$, \ldots, $C(x_i)$. The conditional distribution of $C(x_{i})$ knowing $\mathcal{F}_{i-1}$ is the same as its conditional distribution knowing $\bigcup_{j=1}^{i-1}C(x_j)$. Then, conditionally on the event $\{\bigcup_{j=1}^{i-1}C(x_j)=B\}$, $Y_i=0$ if $x_i\in B$, and $Y_i$ is distributed like the cardinal of the cluster of $x_i$ in $\ZZ^d\setminus  B $ if $x_i\not \in B$. Thus the distribution of $Y_i$ conditionally on $\mathcal{F}_{i-1}$ is stochastically dominated by that of $X_i$. A straightforward induction gives the result.
\end{dem}

We now state that the constant $\nu_G$ is finite.
\begin{prop}
\label{p:finitude}
For any probability measure $G$ on $[0, +\infty]$ such that $G(\{+\infty\})<p_c(d)$, for any $\v\in \SS^{d-1}$, $\nu_G (\v) < +\infty$.
\end{prop}

\begin{dem}
Let $G$ be a probability measure on $[0,+\infty]$ such that $G(\{+\infty\})<p_c(d)$. Let $\v\in \SS^{d-1}$ be a unit vector, let $A$ be a non-degenerate hyperrectangle normal to $\v$ containing the origin $0$ of the graph, and let $h: \NN \mapsto \RR^+$ be mild. 

Let $K_0 <\infty$ be large enough such that $G(]K_0,+\infty]) < p_c(d)$. We recall that for every $x\in \ZZ^d$, $C_{G,K_0} (x)$ is the connected component of $x$ in the percolation $(\mathds{1}_{t_G(e)> K_0}, e \in \EE^d)$. We recall that $\E_{G,K_0}(\cyl(pA,h(p)),h(p))$ denotes the event
$$  \E_{G,K_0}(\cyl(pA,h(p)),h(p)) \,=\, \bigcap_{ x \in \cyl(pA, h(p))\cap \ZZ^d} \{ \diam (C_{G,K_0} (x)) < h(p)  \} \,.$$
To every $x\in B_2(pA, h(p))$, the bottom of the cylinder $\cyl(pA, h(p))$, we associate $S(x) =  \pe C_{G,K_0} (x)$. Some of the sets $S(x)$ may be equal, thus we denote by $(S_i)_{i=1,\dots , r}$ the collection of disjoint edge sets we obtain (notice that by construction for every $i\neq j$, $S_i \cap S_j = \emptyset$). For every $i\in \{1,\dots , r\}$, let $z_i \in B_2 (pA, h(p))$ be such that $S_i = S(z_i)$. 
We consider the set of edges
$$ E(p) \,=\, \bigcup_{i=1}^{r} \left(  S_i \cap \cyl(pA, h(p)) \right) \,. $$
On the event $\E_{G,K_0}(\cyl(pA,h(p)),h(p))$, the set $E(p)$ is a cutset that separates the top $B_1(pA, h(p))$ from the bottom $B_2(pA, h(p))$ of $\cyl(pA, h(p))$. Indeed, let $\gamma = (x_0, e_1, x_1 , \dots ,e_n, x_n )$ be a path from the bottom to the top of $\cyl(pA, h(p))$. There exists $i\in \{1,\dots ,r\}$ such that $x_0 \in \intt (S_i) = \intt ( \pe C_{G,K_0} (z_i))$. 
Since $z_i \in B_2(pA, h(p))$ and $x_n \in B_1(pA, h(p))$ we get $\|z_i-x_n\| \geq 2h(p) - 2 \geq h(p)$ (at least for $p$ large enough), thus on $\E_{G,K_0}(\cyl(pA,h(p)),h(p))$ we know that $x_n \notin \intt (\pe C_{G,K_0} (z_i))$. Let 
$$ k_0 \,=\, \min \{ k \in \{0,\dots , n\} \,:\, x_k \notin  \intt (\pe C_{G,K_0} (z_i)) \} \,.$$
Then $k_0 \in \{1, \dots ,  n\}$, $x_{k_0} \notin \intt (\pe C_{G,K_0} (z_i))$ and $x_{k_0 -1} \in \intt (\pe C_{G,K_0} (z_i))$, thus $e_{k_0} \in  \partial_e C_{G,K_0} (z_i) = S_i $. Since $e_{k_0} \in \gamma \subset \cyl(pA, h(p))$, we conclude that $e_{k_0} \in E(p) \cap \gamma$, thus $E(p)$ cuts the top from the bottom of $\cyl (pA, h(p))$.

For any vertex $x$, by definition of $C_{G,K_0} (x)$ we know that if $e\in \partial_e C_{G,K_0} (x)$ then $t_G(e) \leq K_0$. By definition of $\phi_G (pA, h(p))$, we deduce that on the event $\E_{G,K_0} (\cyl(pA,h(p)),h(p))$ we have
\begin{equation*}
 \phi_{G} (pA, h(p)) \,\leq \, T_G (E(p)) \,\leq \, K_0 \carde (E(p))  \,.
 \end{equation*}
For every $\beta >0$, we obtain that
\begin{align*}
\PP [\phi_{G} & (pA, h(p)) \geq \beta \H^{d-1} (pA)] \\
& \,\leq \, \PP[ \E_{G,K_0} (\cyl(pA,h(p)),h(p)) ^c] + \PP \left[\carde (E(p)) \geq \frac{\beta \H^{d-1} (pA)}{K_0}  \right]\\
& \,\leq \, \cardv (\cyl(pA, h(p))\cap \ZZ^d ) \PP [\diam (C_{G,K_0} (0)) \geq h(p)] + \PP \left[ \sum_{i=1}^{r} \carde (S_i) \geq \frac{\beta \H^{d-1} (pA)}{K_0}  \right]\,.
\end{align*}
We now want to use the stochastic comparison given by Lemma \ref{lem:domsto}. Consider the set of vertices $W=B_2(pA,h(p))$, the percolation $(\mathds{1}_{t_G(e) >K_0} , e\in \EE^d)$, and associate to each vertex $x\in W$ the variable $Z(x) = \cardv (C_{G,K_0}(x))$. We put an order on $W$ and build the variables $(Y(x), x\in V)$ as in Lemma \ref{lem:domsto}. Then
$$ \sum_{i=1}^r \card_e(S_i) \,=\, \sum_{i=1}^r \carde (\partial_e C_{G,K_0} (z_i)) \,\leq\, c_d  \sum_{i=1}^r \cardv (C_{G,K_0} (z_i) ) \,=\, c_d  \sum_{i=1}^r  Z (z_i) \,\leq\, c_d  \sum_{x\in V} Y(x) $$
since the vertices $z_i$ have been chosen in $V$ such that the sets $C_{G,K_0} (z_i)$ are disjoint. By Lemma \ref{lem:domsto}, noticing that $\cardv (W) \leq c_d  \lfloor \H^{d-1} (pA) \rfloor$, we obtain
\begin{align*}
\PP [\phi_{G} (pA, h(p)) \geq \beta \H^{d-1} (pA)] & \,\leq \, c_d \H^{d-1} (pA) h(p) \, \PP [\diam (C_{G,K_0} (0)) \geq h(p)] \\
& \quad + \PP \left[ \sum_{i=1}^{c_d \lfloor \H^{d-1} (pA) \rfloor} X_i \geq \frac{\beta \H^{d-1} (pA)}{K_0c_d}  \right]
\end{align*}
where the variables $X_i$ are i.i.d. with the same distribution as $\cardv ( C_{G,K_0} (0))$. Since $G(]K_0, +\infty]) < p_c(d)$, the percolation $(\mathds{1}_{t_G(e) > K_0}, e\in \EE^d)$ is sub-critical thus 
$$ \PP [X_1 \geq k] \,\leq \,  \kappa_1 e^{-\kappa_2 k} $$
and
\begin{equation}
\label{e:*}
\PP [\diam (C_{G,K_0} (0)) \geq k] \,\leq\,  \kappa_1 e^{-\kappa_2 k} \,,
\end{equation}
where $\kappa_i$ are constants depending only on $d$ and $G(]K_0, +\infty])$, see for instance Theorems (6.1) and (6.75) in \cite{grimmettt:percolation}. Thus there exists $\lambda (G,d) >0$ such that $\EE[\exp (\lambda X_1)] < \infty$, and we get
\begin{align}
\label{e:hop3}
\PP [\phi_{G} (pA, h(p)) \geq \beta \H^{d-1} (pA)] & \,\leq \, c_d \H^{d-1} (pA) h(p)  \kappa_1 e^{-\kappa_2 h(p)} \nonumber \\
&\quad  \quad + \EE[\exp (\lambda X_1)]^{c_d \H^{d-1} (pA)} e^{-\lambda \beta \H^{d-1} (pA)/K_0}\,.
\end{align}
Since $\lim_{p\rightarrow \infty} h(p)/\log p = +\infty$, the first term of the right hand side of \eqref{e:hop3} vanishes when $p$ goes to infinity. We can choose $\beta (G,d)$ large enough such that the second term of the right hand side of \eqref{e:hop3} vanishes too when $p$ goes to infinity, and we get
$$ \lim_{p\rightarrow \infty}  \PP [\phi_{G} (pA, h(p)) \geq \beta \H^{d-1} (pA)]  \,=\, 0\,.$$
Since for every $K\in \RR^+$, $\phi_{G^K} (pA, h(p)) \leq \phi_G (pA, h(p))$ by coupling (see Equation \eqref{e:couplage}), we get for the same $\beta$ that
\begin{equation}
\label{e:hop1}
 \forall K\in \RR^+\,, \quad  \lim_{p\rightarrow \infty}  \PP [\phi_{G^K} (pA, h(p)) \geq \beta \H^{d-1} (pA)]  \,=\, 0\,.
 \end{equation}
By Theorem \ref{t:oldCV}, we know that for every $K\in \RR^+$, 
\begin{equation}
\label{e:hop2}
\nu_{G^K} (\v) \,=\, \lim_{p\rightarrow \infty} \frac{\phi_{G^K} (pA, h(p)) }{ \H^{d-1} (pA)} \qquad \textrm{a.s.}
\end{equation}
Combining \eqref{e:hop1} and \eqref{e:hop2} we conclude that $ \nu_{G^K} (\v) \leq \beta$ for all $K$, thus $\nu_G(\v) = \lim_{K\rightarrow \infty} \nu_{G^K} (\v) \leq \beta < \infty $. This ends the proof of Proposition \ref{p:finitude}.
\end{dem}

Finally we state that $\nu_G$ satisfies some weak triangular inequality.
\begin{prop}
\label{p:cvx}
Let $G$ be a probability measure on $[0,+\infty]$ such that $G(\{+\infty\})<p_c(d)$. Let $(ABC)$ be a non-degenerate triangle in $\RR^d$ and let $\v_A, \v_B$ and $\v_C$ be the exterior normal unit vectors to the sides $[BC], [AC], [AB]$ in the plane spanned by $A,B,C$. Then
$$ \H^1 ([AB]) \nu_G (\v_C) \,\leq\, \H^1 ([AC]) \nu_G (\v_B) + \H^1 ([BC]) \nu_G (\v_A) \,.$$
As a consequence, the homogeneous extension of $\nu_G$ to $\RR^d$, defined by
$$ \nu_G(0) \,=\, 0 \quad \textrm{and} \quad \forall \w\in \RR^d\smallsetminus\{0\} \,, \,\, \nu_G(\w) \,=\, \|\w\|_2 \nu_G \left( \frac{\w}{\|\w\|_2} \right) $$
is a convex function.
\end{prop}
This proposition is a direct consequence of the corresponding property already known for $G^K$ for all $K$, see Proposition 4.5 in \cite{RossignolTheret08b} (see also Proposition 11.6 and Corollary 11.7 in \cite{Cerf:StFlour}).


\subsection{Truncating capacities}
\label{subsec:truncating}

We first need a new definition. Given a probability measure $G$ on $[0, +\infty]$, a unit vector $\v \in \SS^{d-1}$, a non-degenerate hyperrectangle $A$ normal to $\v$ and a height function $h:\NN \rightarrow \RR^+$, we denote by $E_{G} (pA, h(p))$ the (random) cutset that separates the top from the bottom of the cylinder $\cyl (pA, h(p))$ with minimal capacity, {\em i.e.}, $\phi_G (pA, h(p)) = T_G (E_G (pA, h(p)))$, with minimal cardinality among them, with a deterministic rule to break ties.

Furthermore, in this section, if $E\subset \EE^d$ is a set of edges and $C=\cyl(A,h)$ a cylinder, we shall say that $E$ \emph{ cuts $C$ efficiently} if it cuts the top of $C$ from its bottom and no subset of $E$ does. Notice that $E_{G} (pA, h(p))$ cuts $\cyl (pA, h(p))$ efficiently.

\begin{prop}
\label{p:tronquer}
Let $G$ be a probability measure on $[0,+\infty]$ such that $G(\{+\infty\})<p_c(d)$. Then, for any $\eps>0$ and $\alpha >0$, there exist constants $K_1 $ and $C<1$ such that for every $K\geq K_1$, every unit vector $\v\in \SS^{d-1}$, every non-degenerate hyperrectangle $A$ normal to $\v$, every mild height function $h: \NN \mapsto \RR^+$, and for every $p\in \NN^+$ large enough, we have
$$\PP\left[ \phi_{G}(pA, h(p))\geq \phi_{G^K}(pA, h(p))+\eps p^{d-1}\text{ and }\carde (E_{G^K}(pA, h(p)) )\leq \alpha p^{d-1} \right] \leq  C^{h(p)}\,.$$
\end{prop}

Let us say a few words about the proof before starting it. Proposition \ref{p:tronquer} is the equivalent of Lemma \ref{l:key} in the study of the time constant. The proof of Proposition \ref{p:tronquer} is thus inspired by the proof of Lemma \ref{l:key}. The spirit of the proof is the following: we consider a cutset $E$ which is minimal for the truncated $G^K$-capacities. Our goal is to construct a new cutset $E'$ whose $G$-capacity is not much larger than the $G^K$-capacity of $E$. To obtain this cutset $E'$, we remove from $E$ the edges with huge $G$-capacities, and replace them by some local cutsets whose $G$-capacity is well behaved. In fact, the construction of these local modifications of $E$ is in a sense more natural when dealing with cutsets rather than geodesics.

Before embarking to the proof of Proposition~\ref{p:tronquer}, let us state a lemma related to renormalization of cuts. For a fixed $L\in  \NN^*$, we define $\Lambda_L = [-L/2,L/2]^d $, and we define the family of $L$-boxes by setting, for $\mathbf{i}\in \ZZ^d$,
\begin{equation}
\label{e:Lambda}
\Lambda_L (\mathbf{i}) \,=\, \{ x+L \mathbf{i} \,:\, x\in \Lambda_L \} \,.
\end{equation}
The box $\Lambda_L (\mathbf{i})$ is the translated of the box $\Lambda_L$ by a translation of vector $L \mathbf{i}\in \ZZ^d$. A lattice animal is a finite set which is $\ZZ^d$-connected. For $E\subset\EE^d$, let 
$$\Gamma(E) \,=\, \{ \mathbf{j} \in \ZZ^d \,:\, E \cap \Lambda_L(\mathbf{j}) \neq \emptyset  \} $$
be the set of all $L$-boxes that $E$ intersects. 
\begin{lem}
\label{lem:controleanimal}
Let $A$ be a non-degenerate hyperrectangle and $h$ a positive real number. Let $l(A,h)$ denote the minimum of the edge-lengths of $\cyl (A,h)$. Suppose that $E\subset\EE^d$ cuts $\cyl(A, h)$ efficiently. Then, $\Gamma(E)$ is a lattice animal. Furthermore, there exists a constant $c_d$ depending only on $d$ such that if $l(A,h)\geq c_d L$,
\begin{equation}
\label{e:controleanimal}
 \cardv (\Gamma (E)) \,\leq\, c_d \frac{\carde (E)}{L} 
 \end{equation}
\end{lem}
\begin{dem}
Let us prove that $\Gamma(E)$ is a lattice animal. Since $E$ cuts  $\cyl (A,h)$ efficiently, we know that $E$ is somehow connected. More precisely, let us associate with any edge $e\in \EE^d$ a small plaquette that is a hypersquare of size length $1$, that is normal to $e$, that cuts $e$ in its middle and whose sides are parallel to the coordinates hyperplanes. We associate with $E$ the set $E^*$ of all the plaquettes associated with the edges of $E$, and we can see $E^*$ as a subset of $\RR^d$. Then $E^*$ is connected in $\RR^d$ (see \cite{Kesten:flows} Lemma 3.17 in dimension $3$, but the proof can be adapted in any dimension). Thus $\Gamma(E)$ is $\ZZ^d$-connected. 

Now, let us prove~\eqref{e:controleanimal}. We shall denote by $\Lambda'_L$ the enlarged box 
\begin{equation}
\label{e:Lambda'}
\Lambda'_L (\mathbf{i}) \,=\, \{ x+L \mathbf{i} \,:\, x\in \Lambda_{3L} \} \,=\, \bigcup_{ \mathbf{j} \in \ZZ^d \,:\, \|\mathbf{i}-\mathbf{j} \|_{\infty} \leq 1} \Lambda_L (\mathbf{j})  \,.
\end{equation}
First of all we prove that for every $\mathbf{i}\in \ZZ^d$, if $ E \cap \Lambda_L(\mathbf{i}) \neq \emptyset$, then $\carde (E \cap \Lambda'_L(\mathbf{i})) \geq L/2$. Let $e$ be an edge in $E \cap \Lambda_L(\mathbf{i})$. Since $E\smallsetminus \{e\}$ is not a cutset, there exists a path $\gamma = (x_0, e_1, x_1, \dots , e_n, x_n)$ in $\cyl(A, h)$ from the top to the bottom of $\cyl(A, h)$ such that $\gamma $ does not intersect $E\smallsetminus \{e\}$. Since $E$ is a cutset, this implies that $e\in \gamma$. We will prove that locally, inside $\Lambda'_L(\mathbf{i})\smallsetminus \Lambda_L(\mathbf{i})$, the set $E$ must contain at least $L/2$ edges. To do so,  we shall remove $e$ from $\gamma$ and construct of order $L$ possible bypaths of $e$ for $\gamma$ inside $\Lambda'_L(\mathbf{i})\smallsetminus \Lambda_L(\mathbf{i})$, {\em i.e.}, $L/2$ disjoint paths $\gamma'$ such that $\gamma' \subset \Lambda'_L(\mathbf{i})\smallsetminus \Lambda_L(\mathbf{i})$ and the concatenation of the two parts of $\gamma\smallsetminus \{e\}$ and of $\gamma'$ creates a path in $\cyl(A, h)$ from the top to the bottom of $\cyl(A, h)$, see Figure \ref{f:trunc2}.

\begin{figure}[h!]
\centering
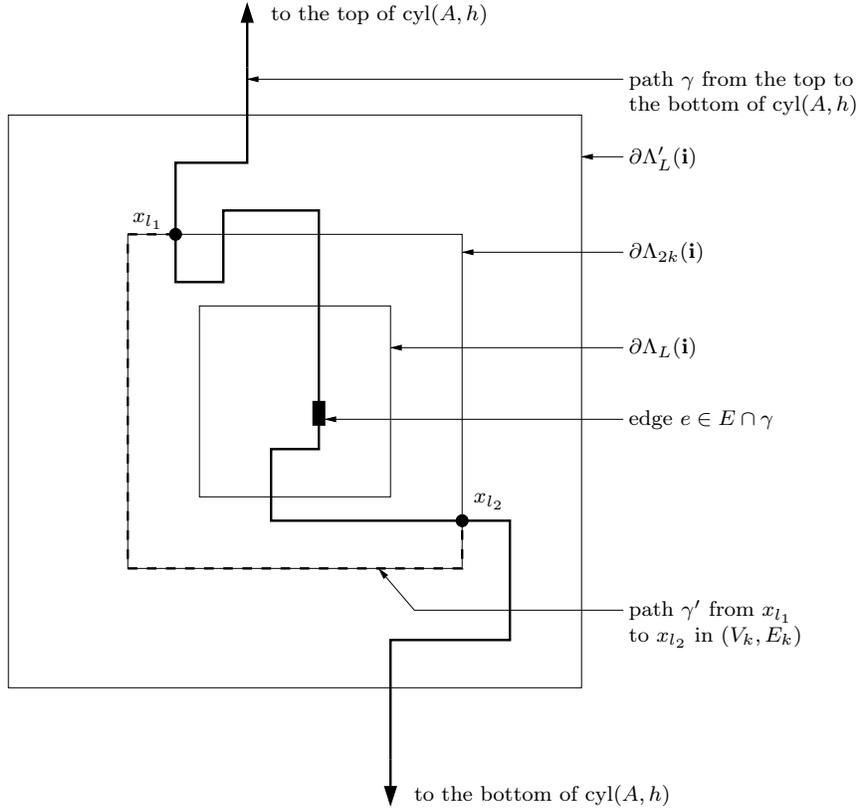
\caption{The path $\gamma$ and the sets $V_k$ ($d=2$).}
\label{f:trunc2}
\end{figure}

For all $k\in [L/2, 3L/2] \cap  \NN$, let $V_k$ be the set of vertices that lies on the faces of $\Lambda_{2k} (\mathbf{i})$, {\em i.e.},
$$ V_k \,=\, \{x+L\mathbf{i} \,:\, x\in \partial \Lambda_{2k} \} $$
and let $E_k$ be the set of edges between vertices in $V_k$,
$$ E_k \,=\, \{ \langle x, y \rangle \in \EE^d \,:\,  x,y \in V_k\} \,.$$
When looking at Figure \ref{f:trunc2}, one sees that the graph $(V_k,E_k)$ forms a kind of shell that surrounds the box $\Lambda_L(\mathbf{i})$. Then any two points $x,y\in V_k$ are connected by a path in $(V_k,E_k)$, and if $x,y$ also belong both to $\cyl (A,h)$ and $h$ is at least $c_dL$ for some constant $c_d$ depending only on the dimension, $x$ and $y$ are also connected by a path in $(V_k\cap \cyl (A,h),E_k\cap \cyl (A,h))$. Let $k\in [L/2, 3L/2] \cap \NN$. 

We claim that the set $(\gamma\smallsetminus \{e\} )\cup (E_k \cap \cyl(A, h) )$ contains a path from the top to the bottom of $\cyl(A, h)$. Let us assume this for the moment, and finish the proof of the lemma. Since the set $(\gamma\smallsetminus \{e\} )\cup (E_k \cap \cyl(A, h) )$ contains a path from the top to the bottom of $\cyl(A, h)$, we know that $E_k$ must intersect the cutset $E$. Since the sets $E_k$ are disjoint, we conclude that
$$ \carde (E \cap \Lambda'_L(\mathbf{i})) \,\geq\, \card ([L/2, 3L/2] \cap \NN) \,\geq L/2\,. $$
This implies that
\begin{align*}
\frac{L}{2} \cardv (\Gamma (E)) & \,\leq\, \sum_{\mathbf{i}\in \Gamma(E)}  \carde (E \cap \Lambda'_L(\mathbf{i}))  \,\leq\, \sum_{\mathbf{i} \in \ZZ^d}  \carde (E \cap \Lambda'_L(\mathbf{i}))\\
& \,\leq\, \sum_{\mathbf{i} \in \ZZ^d}  \sum_{\mathbf{j} \in \ZZ^d\,:\, \|\mathbf{i} - \mathbf{j} \|_{\infty} \leq 1} \carde (E \cap \Lambda_L(\mathbf{j}))\\
& \,\leq\, \sum_{\mathbf{j} \in \ZZ^d} \carde (E \cap \Lambda_L(\mathbf{j})) \card (\{ \mathbf{i} \in \ZZ^d \,:\, \|\mathbf{i} - \mathbf{j} \|_{\infty} \leq 1\})\\
& \,\leq\, 3^d \sum_{\mathbf{j} \in \ZZ^d} \carde (E \cap \Lambda_L(\mathbf{j})) \,\leq \, c_d \carde (E)\,.
\end{align*}
It remains to prove the claim we have left aside, {\it i.e.}, that the set $(\gamma\smallsetminus \{e\} )\cup (E_k \cap \cyl(A, h) )$ contains a path from the top to the bottom of $\cyl(A, h)$. Suppose first that $x_0$ and $x_n$ do not belong to $\Lambda_{2k} (\mathbf{i})$. Then let 
$$ l_1 \,=\, \min \{ l \,:\, x_l \in \Lambda_{2k} (\mathbf{i}) \}  \quad \textrm{and} \quad l_2 \,=\, \max \{ l\,:\, x_l \in \Lambda_{2k} (\mathbf{i}) \} \,,$$
see Figure \ref{f:trunc2}. There exists a path $\gamma'$ from $x_{l_1}$ to $x_{l_2}$ in $(V_k\cap \cyl (A,h),E_k\cap \cyl (A,h))$. We can now concatenate the paths $(x_0, e_1, \dots, x_{l_1})$, $\gamma'$ and $(x_{l_2}, \dots, x_n)$ to obtain a path from the top to the bottom of $\cyl(A, h)$. Suppose now that $x_0 \in \Lambda_{2k} (\mathbf{i})$. Thus, if $l(A,h)$ is at least $c_d L$ for some constant $c_d$ depending only on the dimension, $x_n \notin \Lambda_{2k} (\mathbf{i})$ and there exists a vertex $y \in V_k \cap B_1 (A,h)$ ($B_1 (A, h)$ is the top of the cylinder). We define as previously
$$l_2 \,=\, \max \{ l\,:\, x_l \in \Lambda_{2k}(\mathbf{i}) \} \,.$$
There exists a path $\gamma''$ from $y$ to $x_{l_2}$ in $(V_k\cap \cyl (A,h),E_k\cap \cyl (A,h))$, and we can concatenate $\gamma'$ with $(x_{l_2}, \dots, x_n)$ to obtain a path from the top to the bottom of $\cyl(A, h)$. We can perform the symmetric construction if $x_n \in \Lambda_{2k} (\mathbf{i})$. Thus for every $k\in [L/2, 3L/2] \cap \NN$ the set $(\gamma\smallsetminus \{e\} )\cup (E_k \cap \cyl(A, h) )$ contains a path from the top to the bottom of $\cyl(A, h)$.
\end{dem}

\begin{figure}[h!]
\centering
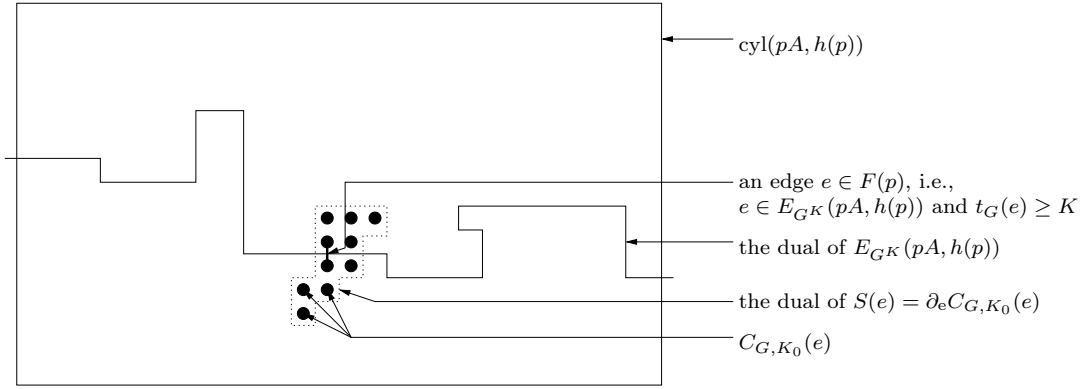
\caption{The cutset $E_{G^K}(pA, h(p))$ and the set $S(e)$ for $e\in F(p)$ ($d=2$).}
\label{f:trunc1}
\end{figure}

{\bf Proof of Proposition~\ref{p:tronquer}:} 
Let $G$ be a probability measure on $[0,+\infty]$ such that $G(\{+\infty\})<p_c(d)$. We use the natural coupling $t_{G^K} (e) = \min (t_G(e),K)$ for all $e\in \EE^d$. Let $K_0$ be such that $G(]K_0,+\infty])<p_c(d)$. We shall modify $E_{G^K}(pA, h(p))$ around the edges having too large $G$-capacities in order to obtain a cut whose capacity is close enough to $\phi_{G^K}(pA,h(p))$ (for $K$ large enough). We recall that $C_{G,K_0}(f)$ is the connected component of the edge $f$ in the percolation $(\mathds{1}_{t_{G}(e)> K_0}, e\in \EE^d)$. For short, we write $S(e) = \pe C_{G,K_0} (e)$, the edge-boundary of $C_{G,K_0}(e)$ separating $e$ from infinity, see Figure \ref{f:trunc1}.

Define also
$$F(p) \,=\, F_{G,K}(pA, h(p)) \,=\,\{e\in E_{G^K}(pA, h(p))\;:\;t_G(e)\geq K\}\;.$$
We collect all the sets $(S(e), e\in F(p))$. As in the proof of Proposition \ref{p:finitude}, from this collection we keep only one copy of each distinct edge set. 
We obtain a collection $(S_i)_{i=1, \dots , r}$ of disjoint edge sets. 
For every $i\in \{1,\dots, r\}$, let $f_i \in F(p)$ such that $S_i = S(f_i)$. Let us define
$$E'(p)\, =\, E_{G,K}'(pA,h(p))= \left( E_{G^K}(pA, h(p))\setminus F(p) \right) \cup \bigcup_{i=1}^{r}\left(S_i \cap\cyl(pA,h(p))\right)\;.$$
We consider the event
$$\E'_{G,K_0} (\cyl(pA,h(p)), h(p))\,=\,\bigcap_{e\in \cyl(pA,h(p))\cap \EE^d} \{\diam C_{G,K_0}(e) < h(p)\}\,.$$
First, we claim that on the event $\E'_{G,K_0} (\cyl(pA,h(p)), h(p))$, the set $E'(p)$ cuts the top from the bottom of $\cyl(pA,h(p))$. 
Indeed, suppose that $\gamma$ is a path in $\cyl(pA,h(p))$ joining its bottom to its top. Since $E_{G^K}(pA, h(p))$ is a cutset, there is an edge $e$ in $E_{G^K}(pA, h(p))\cap\gamma$. If $e$ does not belong to $F(p)$, then $e$ belongs to $E'(p)$ and thus $\gamma$ intersects $E'(p)$. If $e$ belongs to $F(p)$, denote by $x$ (resp. $y$) a vertex belonging to $\gamma$ and the top of $\cyl(pA,h(p))$ (resp. to $\gamma$ and the bottom of $\cyl(pA,h(p))$), and let $i\in \{1,\dots , r\}$ such that $e\in \intt (S_i) =  \intt (S(f_i)) $. On the event $\E'_{G,K_0} (\cyl(pA,h(p)), h(p))$, $x$ and $y$ cannot belong both to $\intt S(f_i)$, otherwise $\diam C_{G,K_0}(f_i)$ would be at least $2h(p)-2 \geq h(p)$ (at least for $p$ large enough). Thus, $\gamma$ contains at least one vertex in $\ext S(f_i)$ and one vertex (any endpoint of $e$) in $ \intt (S(f_i))$. Thus, at least one edge $e'$ of $\gamma$ must be in $S(f_i)$, and since $\gamma$ is included in $\cyl(pA,h(p))$, $e'$ must be in $\cyl(nA,h(p))\cap S(f_i)$. Thus $e'\in E'(p)$ and this proves that $E'(p)$ cuts the top from the bottom of $\cyl(pA, h(p))$.

Now, on the event $\E'_{G,K_0} (\cyl(pA,h(p)), h(p))$ we get
\begin{equation}
\label{e:raf1}
\phi_G(pA,h(p))\leq \phi_{G^K}(pA,h(p))+K_0 \, \sum_{i=1}^r \carde ( S_i ) \,.
\end{equation}
Moreover, still on the event $\E'_{G,K_0} (\cyl(pA,h(p)), h(p))$, we notice that if we replace a single edge $e$ of $F(p)$ by $\left(S(e) \cap\cyl(pA,h(p))\right)$ in $E_{G^K}(pA, h(p))$ we obtain a new set of edges that is still a cutset between the top and the bottom of $\cyl(pA, h(p))$ (this could be proved by a similar but simpler argument than the one presented to prove that $E'(p)$ is a cutset). By minimality of the capacity of $E_{G^K}(pA, h(p))$ among such cutsets, we deduce that
\begin{equation}
\label{e:raf2}
\forall e \in F(p) \,,\quad K_0 \carde (S(e)) \,\geq \, K \,.
\end{equation}
We recall that $\carde (S(e)) \leq c_d \cardv (C_{G,K_0} (e))$. Furthermore, notice that if $e_i=\langle z^1_i,z^2_i \rangle$, then $C_{K_0}(e) = C_{K_0}(z^1_i) \cup C_{K_0}(z^2_i) $ and  $\cardv (C_{K_0}(e))\leq \cardv (C_{K_0}(z^1_i))+\cardv (C_{K_0}(z^2_i))$. Consequently, 
$$\max_{j=1,2} \cardv (C_{K_0}(z^j_i)) \geq \cardv (C_{K_0}(e)) /2 $$
Let us denote by $B$ the event whose probability we want to bound from above:
$$B:=\{ \phi_G (pA,h(p))\geq \phi_{G^K}(pA,h(p))+\eps p^{d-1}\textrm{ and }\carde (E_{G^K}(pA, h(p)) )\leq \alpha \, p^{d-1}\}$$
and for positive  $\beta$, $\gamma$ and for $x_1,\ldots,x_k$ in $\ZZ^d$, let
$$B_{G,K_0}(x_1,\ldots,x_k;\beta,\gamma):=\left[\begin{array}{c}(  C_{G,K_0}(x_i) )_{i\geq 1}\text{ are pairwise disjoint,}\\
 \sum_{i=1}^k \cardv (C_{G,K_0}(x_i)) \geq \beta \\
 \text{ and }\forall i=1,\ldots,k,\; \cardv (C_{G,K_0}(x_i)) \geq  \gamma \end{array} \right]\;.$$
If $E\subset \EE^d$ and $x\in \ZZ^d$, we say that $x\in E$ if and only if $x$ is the endpoint of an edge $e$ that belongs to $E$. We obtain this way
\begin{align*}
\PP & \left[B\cap \E'_{G,K_0} (\cyl(pA,h(p)), h(p)) \right] \\
& \,\leq \,\PP\left[ \begin{array}{c} \exists E \subset \EE^d \;:\; \carde (E )\leq \alpha\,  p^{d-1}\,,\, E \textrm{ cuts }\cyl (pA,h(p))\textrm{ efficiently}\,,\\ \exists k\geq 0\,,\, \exists x_1,\ldots, x_k\in E\text{ such that } 
B_{G,K_0}\left(x_1,\ldots,x_k;\frac{ \eps p^{d-1}}{2 c_d K_0},\frac{K}{2 c_d K_0}\right)\textrm{ holds} \end{array}\right]\,.
\end{align*}
As in the proof of the continuity of the time constant given by Cox and Kesten in \cite{CoxKesten}, we need a renormalization argument to localize these vertices $x_1,\ldots, x_k$ in a region of the space whose size can be controlled. For a given $\mathbf{i}\in \ZZ^d$ and $k\in \NN^*$, we denote by $\mathcal{A}(\mathbf{i}, k)$ the set of all lattice animals of size $k\in \NN^*$ containing $\mathbf{i}$. If $k\in \mathbb{R}^+$, then we write $\mathcal{A}(\mathbf{i}, k)$ instead of $\mathcal{A}(\mathbf{i}, \lfloor k \rfloor)$ for short, where $\lfloor k \rfloor \in \mathbb{N}$ and satisfies $\lfloor k \rfloor \leq k < \lfloor k \rfloor+1$. 

Let $E \subset \EE^d$ such that $E$ cuts  $\cyl (pA,h(p))$ efficiently. Let us denote by $u\in \RR^d$ one of the corners of $A$. We can find a path from the top to the bottom of $\cyl (pA, h(p))$ that is located near any of the vertical sides of $\cyl (pA, h(p))$, more precisely there exists a constant $c_d$ depending only on $d$ such that the top and the bottom of $\cyl(pA,h(p))$ are connected in $V(u, h(p)) := \{ u+l\v + \vec w \,:\, l\in [-h(p), h(p)] \,,\, \|\vec w\|_2 \leq c_d \} \cap \cyl (pA, h(p))$. Thus any custset $E$ must contain at least one edge in $V(u,h(p))$. We denote by $\mathbf{I} (pA, h(p))$ the set of $L$-boxes that intersect $V(u,h(p))$:
$$ \mathbf{I} (pA, h(p)) \,=\, \{ \mathbf{i}\in \ZZ^d \,:\, \Lambda_L(\mathbf{i}) \cap V(u,h(p))\neq \emptyset \} \,.$$
Then $\Gamma(E)$ must intersect $ \mathbf{I} (pA, h(p))$, and $\cardv ( \mathbf{I} (pA, h(p))) \leq c_d h(p)/L$. 

Furthermore, Lemma~\ref{lem:controleanimal} ensures that for $p$ large enough,
$$
\cardv (\Gamma (E)) \,\leq\, c_d \frac{\carde (E)}{L} \,.
$$

From all these remarks, we conclude that if $E$ cuts $\cyl (pA,h(p))$ efficiently and if $\carde (E )\leq \alpha\,  p^{d-1}$, then 
$$\Gamma(E) \,\in\,  \bigcup_{\mathbf{i}\in \mathbf{I} (pA, h(p))} \bigcup_{k\leq c_d \alpha p^{d-1}/L} \mathcal{A}(\mathbf{i}, k) \,.$$
Notice that for any $\Gamma \in \mathcal{A}(\mathbf{i}, k)$ with $k\leq c_d \alpha p^{d-1}/L$, there exists $\Gamma' \in \mathcal{A}(\mathbf{i}, c_d \alpha p^{d-1}/L)$ such that $\Gamma\subset \Gamma'$. Thus if $\carde (E )\leq \alpha\,  p^{d-1}$, we obtain that there exists $\mathbf{i}\in \mathbf{I} (pA, h(p))$ and $\Gamma \in \mathcal{A}(\mathbf{i}, c_d \alpha p^{d-1}/L)$ such that $\Gamma (E) \subset \Gamma$. For any lattice animal $\Gamma$, we denote by $\Gamma_L$ the uion of the boxes associated to this vertex set, {\em i.e.},
$$ \Gamma_L \,=\,\bigcup_{\mathbf{j} \in \Gamma}  \Lambda_L (\mathbf{j}) \,.$$
We obtain
\begin{align*}
\PP & \left[ B \cap \E'_{G,K_0} (\cyl(pA,h(p)), h(p)) \right] \\
&\,\leq \,  \sum_{\mathbf{i}\in \mathbf{I} (pA, h(p))} \sum_{\Gamma \in \mathcal{A}(\mathbf{i}, c_d \alpha p^{d-1}/L)}
\PP \left[ \begin{array}{c}  \exists k\geq 0\,,\, \exists \textrm{ vertices }x_1,\ldots, x_k\in \Gamma_L\text{ such that }\\
B_{G,K_0}\left(x_1,\ldots,x_k;\frac{ \eps p^{d-1}}{2 c_d K_0},\frac{K}{2 c_d K_0}\right)\textrm{ holds} \end{array} \right]\,.\\
&\,\leq \,  \sum_{\mathbf{i}\in \mathbf{I} (pA, h(p))} \sum_{\Gamma \in \mathcal{A}(\mathbf{i}, c_d \alpha p^{d-1}/L)} \sum_{k\in \NN^*} \sum_{x_1,\ldots, x_k\in \Gamma_L}  
\PP \left[B_{G,K_0}\left(x_1,\ldots,x_k;\frac{ \eps p^{d-1}}{2 c_d K_0},\frac{K}{2 c_d K_0}\right) \right] \,,
 \end{align*}
We now use the stochastic comparison given by Lemma \ref{lem:domsto}. We consider the set of vertices $W=\{x_1,\dots, x_k\} $, the percolation $(\mathds{1}_{t_G(e)>K_0}, e\in \EE^d)$ and associate to each vertex $x_i$ the variable $Z_i = Z (x_i) = \cardv (C_{G,K_0}(x_i))$. We build the variables $(Y_i, 1\leq i \leq k)$ as in Lemma \ref{lem:domsto} and let $(X_i, 1\leq i \leq k)$ be i.i.d. random variables with distribution $\cardv (C_{G,K_0}(e))$. Then by Lemma \ref{lem:domsto} we obtain for $\lambda>0$,
\begin{align*}
\PP  \left[ B_{G,K_0}(x_1,\ldots,x_k;\beta,\gamma)  \right] 
& \,= \, \PP \left[    \sum_{i=1}^k  Y_i \geq \beta  \text{ and }\forall i=1,\ldots,k,\; Y_i \geq  \gamma  \right] \\
& \,\leq \, \PP \left[\sum_{i=1}^k X_i \geq \beta  \text{ and }\forall i=1,\ldots,k,\; X_i \geq  \gamma \right] \\
& \,\leq \, e^{-\lambda \beta} \EE\left[ e^{\lambda X_1 }\II_{X_1 \geq \gamma} \right]^k\\
& \,\leq \, e^{-\lambda \beta} \EE\left[ e^{2\lambda X_1 }\right]^{\frac{k}{2}}\PP\left[X_1 \geq \gamma \right]^{\frac{k}{2}}
\end{align*}
where we used the Cauchy-Schwartz inequality. Thus, we get
\begin{align}
\label{e:raf4}
\PP & \left[ B \cap \E'_{G,K_0} (\cyl(pA,h(p)), h(p)) \right] \nonumber\\
&\,\leq \,  \sum_{\mathbf{i}\in \mathbf{I} (pA, h(p))} \sum_{\Gamma \in \mathcal{A}(\mathbf{i}, c_d \alpha p^{d-1}/L)} \sum_{k\in \NN^*} \sum_{x_1,\ldots, x_k\in \Gamma_L}  
e^{-\lambda \frac{\eps p^{d-1}}{2 c_d K_0}} \EE \left[ e^{2 \lambda X_1 } \right]^{\frac{k}{2}} \PP \left[ X_1 \geq \frac{K}{2 c_d K_0} \right] ^{\frac{k}{2}}\nonumber  \\
& \,\leq \,e^{-\lambda \frac{\eps p^{d-1}}{2 c_dK_0}}   c_d \frac{h(p)}{L} c_d^{ \frac{\alpha p^{d-1}}{L}} \sum_{k\in \NN^*} \left( \begin{array}{c} c_d \alpha p^{d-1}L^{d-1} \\ k  \end{array} \right) \EE \left[ e^{2 \lambda X_1 } \right]^{\frac{k}{2}} \PP \left[X_1 \geq \frac{K}{2 c_d  K_0} \right] ^{\frac{k}{2}}\nonumber  \\
& \,\leq \, c_d \frac{h(p)}{L}  \left[ e^{-\lambda \frac{\eps }{2 c_d K_0}}  c_d^{ \frac{\alpha }{L}} \left( 1 +  \EE \left[ e^{2 \lambda X_1 } \right]^{\frac{1}{2}}  \PP \left[X_1 \geq \frac{K}{2 c_d  K_0} \right] ^{\frac{1}{2}} \right)^{c_d \alpha L^{d-1}} \right]^{p^{d-1}} \,. 
\end{align}
Now, since $G(]K_0,+\infty])<p_c(d)$ we can choose first $\lambda = \lambda (G, d)$ such that 
$$\EE\left[ e^{2\lambda X_1 } \right]<\infty \,.$$
Then, we choose $L = L(G, d, \alpha, \eps)$ large enough such that
$$e^{-\lambda \frac{\eps }{2 c_d K_0}}c_d^\frac{\alpha}{L}<1\,.$$
And finally we choose $K_1 = K_1 (G,d,\alpha, \eps)$ such that:
\begin{equation}
\label{e:raf5}
C(G,d,\alpha, \eps)\,:=\, e^{-\lambda \frac{\eps }{2 c_d K_0}}c_d^\frac{\alpha}{L}\left( 1 +  \EE[e^{2 \lambda X_1}]^{\frac{1}{2}}  \PP \left[c_d X_1\geq \frac{K_1}{2 K_0} \right] ^{\frac{1}{2}} \right)^{c_d \alpha L^{d-1}}\,<\,1 \,.
\end{equation}
Combining \eqref{e:raf4} and \eqref{e:raf5} we get
\begin{eqnarray*}
\PP \left[ B \cap \E'_{G,K_0} (\cyl(pA,h(p)), h(p)) \right] &\leq& c_d \frac{h(p)}{L} C(G,d,\alpha, \eps)^{p^{d-1}}\\
& \leq & C(G,d,\alpha, \eps)^{\frac{1}{2} p^{d-1}}
\end{eqnarray*}
for some $C(G,d,\alpha, \eps)<1$, for every $K\geq K_1 (G,d,\alpha, \eps)$ and for every $p$ large enough, since $\lim_{p\rightarrow \infty}  h(p)/p =0$. For every edge $e=\langle x,y \rangle$, since 
$$\diam (C_{G,K_0} (e)) \,\leq\, \diam (C_{G,K_0} (x)) + \diam (C_{G,K_0} (y)) +1 \,, $$ 
it is easy to show that
$$ \PP[\E'_{G,K_0} (pA, h(p))^c]\,\leq\, c_d \H^{d-1}(pA) h(p) \kappa_1 e^{-\kappa_2 h(p)} \,\leq \,e^{-\frac{\kappa_2}{2} h(p)}  $$
for some positive constants $\kappa_i $ (see \eqref{e:*}) and for $p$ large enough since $\lim_{p\rightarrow \infty} h(p) / \log p = +\infty$.
Since $\lim_{p\rightarrow \infty}h(p)/p =0$, this ends the proof of Proposition \ref{p:tronquer}.
{\flushright$\blacksquare$\\}

\begin{rem}
The result of Proposition \ref{p:tronquer} could apply, with the same constants depending on $G$, to any probability measure $H$ on $[0,+\infty]$ such that $H \preceq G$ (we recall that the stochastic comparison $H\preceq G$ is defined in \eqref{e:defdomsto}).
\end{rem}


\subsection[Proof of the convergence I]{Proof of the convergence I: case $G(\{0\})<1-p_c(d)$}

To prove Theorem \ref{t:CV}, we shall consider the two situations $G(\{0\})<1-p_c(d)$ and $G(\{0\})\geq 1-p_c(d)$. The purpose of this section is to prove the following proposition, that corresponds to the statement of Theorem \ref{t:CV} in the case where $G(\{0\})<1-p_c(d)$.
\begin{prop}
\label{p:positif}
For any probability measure $G$ on $[0, +\infty]$ such that $G(\{+\infty\})<p_c(d)$ and $G(\{0\}) < 1-p_c(d)$, for any $\v \in \SS^{d-1}$, for any non-degenerate hyperrectangle $A$ normal to $\v$, for any mild function $h : \NN \mapsto \RR^+$, we have
$$  \lim_{p\rightarrow \infty} \frac{\phi_G (pA, h(p))}{\H^{d-1} (pA)} \,=\, \nu_G(\v) \quad \textrm{a.s.}$$
\end{prop}

\begin{dem}
Let $\v\in \SS^{d-1}$, let $A$ be a non-degenerate hyperrectangle normal to $\v$, let $h:\NN^* \mapsto \RR^+$ be mild. Let $G$ be a probability measure on $[0,+\infty]$ such that $G(\{+\infty\})<p_c(d)$. Since $d,G,\v, A, h$ are fixed, we will omit in the notations a lot of dependences in these parameters.

In this section, we suppose that $G(\{ 0 \}) < 1-p_c(d)$. For any fixed $K\in \RR^+$, we know by Theorem~\ref{t:oldCV} that a.s.
$$ \liminf_{p\rightarrow \infty} \frac{\phi_G (pA, h(p))}{\H^{d-1} (pA)} \,\geq\, \lim_{p\rightarrow \infty} \frac{\phi_{G^K} (pA, h(p))}{\H^{d-1} (pA)} \,=\, \nu_{G^K} (\v) \,, $$
thus
$$ \liminf_{p\rightarrow \infty} \frac{\phi_G (pA, h(p))}{\H^{d-1} (pA)} \,\geq\, \sup_K \nu_{G^K} (\v) \,=\, \nu_G (\v) \,. $$
It remains to prove that a.s., 
$$ \limsup_{p\rightarrow \infty} \frac{\phi_G (pA, h(p))}{\H^{d-1} (pA)} \,\leq\, \nu_G(\v)\,. $$
We claim that it is sufficient to prove that 
\begin{equation}
\label{e:flop1}
  \forall \eps >0\,,\, \exists K(\eps)\quad \sum_{p\geq1} \PP [\phi_G(pA, h(p)) \geq \phi_{G^K} (pA, h(p)) + \eps \H^{d-1} (pA)] \,<\, +\infty \,.
\end{equation}
Indeed, if \eqref{e:flop1} is satisfied, by Borel-Cantelli and Theorem \ref{t:oldCV} it implies that
$$   \forall \eps >0\,,\, \exists K(\eps)\,,\, \textrm{a.s.}\quad \limsup_{p\rightarrow \infty}   \frac{\phi_G (pA, h(p))}{\H^{d-1} (pA)} \,\leq\, \nu_{G^{K}} (\v)+ \eps \,\leq\, \nu_G(\v) + \eps $$
and Proposition \ref{p:positif} is proved. Now for every $\eps>0$, for every $\alpha >0$ and every $\beta >0$, we have
\begin{align}
\label{e:flop2}
\PP [ & \phi_G(pA, h(p)) \geq \phi_{G^K} (pA, h(p)) + \eps \H^{d-1} (pA)] \nonumber \\
&\,\leq\, \PP [\{\phi_G(pA, h(p)) \geq \phi_{G^K} (pA, h(p)) + \eps \H^{d-1} (pA)\} \cap \{ \carde (E_{G^K} (pA, h(p)) \leq \alpha p^{d-1}  \}] \nonumber\\
&\qquad + \PP [\{\carde (E_{G^K} (pA, h(p)) > \alpha p^{d-1} \} \cap \{ \phi_{G^K} (pA,h(p)) \leq \beta \H^{d-1} (pA )\}] \nonumber\\
&\qquad + \PP [\phi_{G^K} (pA,h(p))> \beta \H^{d-1} (pA )  ]
\end{align}

By \eqref{e:hop3}, since $\phi_{G^K} (pA,h(p))\leq \phi_{G} (pA,h(p))$ by coupling (see Equation \eqref{e:couplage}) and $\lim_{p\rightarrow \infty} h(p)/\log p = +\infty$, we know that we can choose $\beta = \beta (G,d)$ such that for any $K\in \RR^+$, the last term of the right hand side of \eqref{e:flop2} is summable in $p$.

Given this $\beta (G,d)$, by Zhang's Theorem 2 in \cite{Zhang2017}, as adapted in Proposition 4.2 in \cite{RossignolTheret08b}, we know that since all the probability measures $G^K$ coincide on a neighborhood of $0$, we can choose a constant $\alpha (G,d)$ such that for any $K\in \RR^+$ the second term of the right hand side of \eqref{e:flop2} is summable in $p$.

Given this $\alpha (G,d)$, by Proposition \ref{p:tronquer}, we know that there exist some constants $C=C(G,d,\eps)<1$ and $K_1(G,d,\eps)$ such that for every $K\geq K_1(G,d,\eps)$ and for all $p$ large enough 
\begin{equation}
\label{e:flop3}
 \PP [\{\phi_G(pA, h(p)) \geq \phi_{G^K} (pA, h(p)) + \eps \H^{d-1} (pA)\} \cap \{ \carde (E_{G^K} (pA, h(p)) \leq \alpha p^{d-1}  \}] \,\leq\,  C^{h(p)}\,.
\end{equation}
The right hand side of \eqref{e:flop3} is summable in $p$ since $\lim_{p\rightarrow \infty} h(p)/\log p = +\infty$. This concludes the proof of \eqref{e:flop1}, thus the convergence in Theorem \ref{t:CV} is proved when $G(\{0\})<1-p_c(d)$.
\end{dem}


\subsection[Proof of the convergence II]{Proof of the convergence II : case $G(\{0\})\geq 1-p_c(d)$}

It remains to prove that the convergence in Theorem \ref{t:CV} holds when $G(\{0\})\geq 1-p_c(d)$, {\em i.e.}, when $\nu_G =0$. We first deal with straight cylinders. For $A = \prod_{i=1}^{d-1} [0,k_i] \times \{0\}$ (with $k_i>0$ for all $i$) and $h \in \NN$, we denote by $\phi^{\v_0}_G (A,h) $ the maximal flow from the top $B_1^{\v_0}(A,h) = (\prod_{i=1}^{d-1} [0,k_i] \times \{h\}) \cap \ZZ^d$ to the bottom $B_2^{\v_0}(A,h)=(\prod_{i=1}^{d-1} [0,k_i] \times \{0\})\cap \ZZ^d$ in the cylinder $\cyl^{\v_0}(A,h)=\prod_{i=1}^{d-1} [0,k_i] \times [0,h]$ for $\v_0 = (0,\dots, 0,1)$. We recall the definition of the event $\E_{G,K} (\C, h)$, for any subset $\C$ of $\RR^d$ and any $h\in \RR^+$, that was given in \eqref{e:E}:
$$ \E_{G,K} (\C, h) \,=\, \bigcap_{x\in \C \cap \ZZ^d} \{ \diam (C_{G,K} (x)) < h \}\,.$$

\begin{prop}
\label{p:zero}
Let $G$  be a probability measure on $[0, +\infty]$ such that $G(\{+\infty\})<p_c(d)$ and $G(\{0\}) \geq 1-p_c(d)$. Let $A = \prod_{i=1}^{d-1} [0,k_i] \times \{0\}$ (with $k_i>0$ for all $i$). For any function $h : \NN \mapsto \NN$ satisfying $\lim_{p\rightarrow +\infty} h(p) /\log p =+\infty$, we have
$$  \lim_{p\rightarrow \infty} \frac{\phi^{\v_0}_G (pA, h(p))}{\H^{d-1} (pA)} \,=\, 0 \quad \textrm{a.s.}$$
Moreover, if $G(]K_0,+\infty]) < p_c(d)$, then we also have  
$$ \lim_{p\rightarrow \infty} \frac{\phi^{\v_0}_G (pA, h(p)) \mathds{1}_{ \E_{G,K_0}( \cyl^{\v_0}(pA,h(p)) ,h(p))}}{ \H^{d-1} (pA)} \,=\, 0 \quad \textrm{in }L^1\,. $$
\end{prop}
This result is in fact a generalization of Zhang's Theorem \ref{t:oldnul}, and the strategy of the proof is indeed largely inspired by Zhang's proof. However, we need to work a little bit harder, because we do not have good integrability assumptions. We thus re-use here some ideas that appeared in the proof of Proposition \ref{p:finitude}. Notice that $\phi_G (pA, h(p))$ itself may not be integrable in general (it can even be infinite with positive probability).

\begin{dem}
We shall construct a particular cutset with an idea quite similar to the one we used in the proof of Proposition \ref{p:finitude}. Let $K_0$ be large enough to have $G(]K_0,+\infty]) < p_c(d)$. Let $h\in \NN^*$. Let $\HH$ be the half-space $\ZZ^{d-1}\times \NN$. For any $x\in \ZZ^d$, $D_1 \subset \ZZ^d$ and $D_2 \subset \ZZ^d$, let us denote by $\left\{x \overset{D_2}{\underset{G,0}{\longleftrightarrow}} D_1 \right\}$ the event
$$ \left\{x \overset{D_2}{\underset{G,0}{\longleftrightarrow}} D_1 \right\} \,=\, \{ x \textrm{ is connected to $D_1$ by a path $\gamma \subset D_2$ s.t. } \forall e\in \gamma \,,\, t_G(e)>0 \} $$ 
For any $x\in \ZZ^{d-1}\times \{0\}$, we define the event
$$ F_{x,\ell} \,=\, \left\{ x \overset{\HH}{\underset{G,0}{\longleftrightarrow}} \ZZ^{d-1} \times \{\ell\} \right\} \,.$$
Let $x\in A$. If $F_{x,\ell}^c$ occurs, we associate with $x$ the set $\pe C_{G,0} (x)$, that is by definition made of edges with null capacity. If $F_{x,\ell}$ occurs, we associate with $x$ the set $ \pe C_{G,K_0} (x)$, see Figure \ref{f:zero}. 
\begin{figure}[h!]
\centering
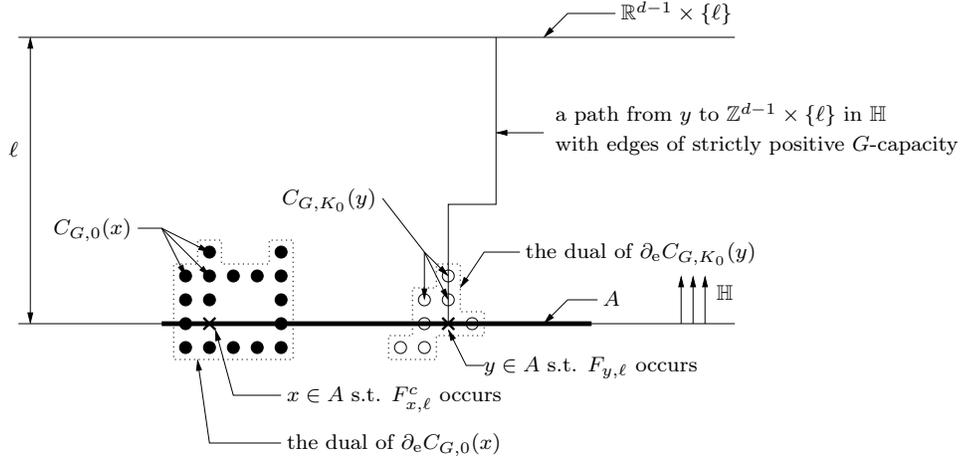
\caption{The construction of the cutset $E'(A,\ell)$ in $\cyl^{\v_0} (A,\ell)$ ($d=2$).}
\label{f:zero}
\end{figure}
We consider the set
$$ E'(A,\ell) \,=\, \left [\left(\bigcup_{x\in A\cap \ZZ^d \,,\, F_{x,\ell}^c} \pe C_{G,0} (x) \right) \cup \left( \bigcup_{x\in A\cap \ZZ^d \,,\, F_{x,\ell}} \pe C_{G,K_0} (x)  \right) \right] \cap \cyl^{\v_0} (A,\ell) \,. $$
We consider the good event
$$   \E_{G,K_0}\left( \cyl^{\v_0}(A, \ell) ,\ell \right)\,=\, \bigcap_{ x \in \cyl^{\v_0}(A, \ell)  \cap \ZZ^d} \{ \diam (C_{G,K_0} (x)) < \ell  \} \,.$$
We claim that on $   \E_{G,K_0}( \cyl^{\v_0}(A, \ell) ,\ell)$, the set $E'(A,\ell)$ cuts the top $B_1^{\v_0}(A,\ell)$ from the bottom $B_2^{\v_0}(A,\ell)$ in the cylinder $\cyl^{\v_0}(A,\ell)$. Let $\gamma = (x_0, e_1, x_1 , \dots ,e_n, x_n )$ be a path from the bottom to the top of the cylinder $\cyl^{\v_0}(A,\ell)$. If $F_{x_0,\ell}^c$ occurs, since $x_n \in \ZZ^{d-1} \times \{\ell\}$, then $x_n \in \ext (\pe C_{G,0} (x_0))$, thus $\gamma$ has to use an edge in $\pe C_{G,0} (x)\cap \cyl^{\v_0} (A,\ell)$. If $F_{x_0,\ell}$ occurs, on $ \E_{G,K_0}(\cyl^{\v_0} (A,\ell), \ell)$ we know that $x_n \in \ext (\pe C_{G,K_0} (x_0))$, thus $\gamma$ must contain an edge in $\pe C_{G,K_0} (x_0)  \cap \cyl^{\v_0}(A,\ell)$. We conclude that on $ \E_{G,K_0}(\cyl^{\v_0}(A,\ell),\ell)$, $E'(A,\ell)$ is indeed a cutset from the top to the bottom of $\cyl^{\v_0}(A,\ell)$. Thus on $  \E_{G,K_0}( \cyl^{\v_0}(A,\ell),\ell)$,
\begin{equation}
\label{e:pif1}
 \phi^{\v_0}_G (A,\ell) \,\leq\, K_0 \sum_{x\in A\cap \ZZ^d} \mathds{1}_{F_{x,\ell}} \carde \left(  \pe C_{G,K_0} (x) \right)  \,\leq\, c_d K_0 \sum_{x\in A\cap \ZZ^d} \mathds{1}_{F_{x,\ell}} \cardv \left(  C_{G,K_0} (x) \right) \,.
 \end{equation}
For every $x\in \ZZ^{d-1} \times \{0\}$, let us define
$$ R^\ell _x \,=\,  \mathds{1}_{F_{x,\ell}} \cardv \left(  C_{G,K_0} (x) \right) \,, $$
and for every $D \in \mathcal J = \{ \prod_{i=1}^{d-1} [l_i,l_i']\,:\, \forall i \,,\,0\leq  l_i\leq l_i'  \}$, let us define
$$  X^\ell_D \,=\, \sum_{x\in D\times \{0\} \cap \ZZ^d} R^\ell_x \,.$$
For every $\ell$, the process $(X^\ell_D, D \in \mathcal J)$ is a discrete additive process. By classical multiparameter ergodic Theorems (see for instance Theorem 2.4 in \cite{Akcoglu} and Theorem 1.1 in \cite{Smythe}), if $\EE[R^\ell_0] <\infty$, then there exists an integrable random variable $X^\ell$ such that for every $D = \prod_{i=1}^{d-1} [0,k_i]$ (with $k_i\in \NN^*$ for all $i$), 
\begin{equation}
\label{e:pif2}
\lim_{p\rightarrow \infty} \frac{1}{\cardv (pD \cap \ZZ^{d-1})} X^\ell_{pD} =X^\ell \quad \textrm{a.s. and in }L^1
\end{equation}
and $\EE[X^\ell] = \EE[R_0^\ell]$. Moreover, by ergodicity $X^\ell$ is constant a.s., so 
\begin{equation}
\label{e:pif3}
X^\ell \,=\, \EE[R_0^\ell] \qquad \textrm{a.s.} 
\end{equation}
We need to control the expectation of $R_0^\ell$ to apply these ergodic theorems. For all $r>0$ we have by the independence of the edge weights
\begin{align*}
\PP [  R_0^\ell = r] 
& \,\leq \, \PP[ F_{0,\ell} \cap \{   \cardv (  C_{G,K_0} (0)) = r \} ] \\
& \,\leq \, \sum_{C\,:\, \cardv (C)= r} \,\sum_{v\in  \pv C \cap \HH} \PP \left[ \{C_{G,K_0} (0) = C\} \cap \left\{ v \overset{\HH\smallsetminus C}{\underset{G,0}{\longleftrightarrow}} \ZZ^{d-1} \times \{\ell\} \right\} \right] \\
& \,\leq \, \sum_{C\,:\, \cardv (C) = r} \,\sum_{v\in  \pv C \cap \HH} \PP \left[ C_{G,K_0} (0) = C \right] \PP\left[ v \overset{\HH\smallsetminus C}{\underset{G,0}{\longleftrightarrow}} \ZZ^{d-1} \times \{\ell\} \right] \\
& \,\leq \, \sum_{C\,:\, \cardv (C)= r} \,\sum_{v\in  \pv C \cap \HH} \PP \left[ C_{G,K_0} (0) = C \right] \PP\left[ 0 \overset{\HH}{\underset{G,0}{\longleftrightarrow}} \ZZ^{d-1} \times \{\ell-(r+1)\}  \right] \\
& \,\leq \, c_d\, r\,  \PP\left[ \cardv (  C_{G,K_0} (0)) = r \right]\PP\left[ 0 \overset{\HH}{\underset{G,0}{\longleftrightarrow}} \ZZ^{d-1} \times \{\ell-(r+1)\}   \right]\,.
\end{align*} 
Using the fact that $k\mapsto \PP\left[ 0 \overset{\HH}{\underset{G,0}{\longleftrightarrow}} \ZZ^{d-1} \times \{k\}   \right] $ is non-increasing, we get
\begin{align*}
\EE[R_0^\ell] & \,=\, \sum_{r\in \NN} r \,\PP [  R_0^\ell = r]\\
& \,\leq \, \sum_{r\leq \ell/2} c_d \,r^2 \, \PP\left[ \cardv (  C_{G,K_0} (0)) = r \right]  \PP\left[ 0 \overset{\HH}{\underset{G,0}{\longleftrightarrow}} \ZZ^{d-1} \times \{\ell-(r+1)\}  \right]\\
& \qquad + \sum_{r> \ell/2} c_d\, r^2\,  \PP\left[ \cardv (  C_{G,K_0} (0)) = r \right]  \PP\left[ 0 \overset{\HH}{\underset{G,0}{\longleftrightarrow}} \ZZ^{d-1} \times \{\ell-(r+1)\}  \right]\\
& \,\leq\, c_d \,\PP\left[ 0 \overset{\HH}{\underset{G,0}{\longleftrightarrow}} \ZZ^{d-1} \times \{\ell/2 - 1\}  \right] \sum_{r\in \NN} r^2\, \PP\left[ \cardv (  C_{G,K_0} (0)) = r \right]\\
&\qquad  + c_d \sum_{r> \ell/2}  r^2  \,\PP\left[ \cardv (  C_{G,K_0} (0)) = r \right] \\
& \,\leq\, c_d \,\PP\left[ 0 \overset{\HH}{\underset{G,0}{\longleftrightarrow}} \ZZ^{d-1} \times \{\ell/2 - 1\}  \right] \EE \left[ \cardv (  C_{G,K_0} (0)) ^2 \right]  \\
& \qquad +  \EE \left[ \cardv (  C_{G,K_0} (0)) ^2  \mathds{1}_{ \cardv (  C_{G,K_0} (0)) > \ell/2 }\right] \,. 
\end{align*}
Since we choose $K_0$ such that $G(]K_0,+\infty]) < p_c(d)$, we know that $\EE \left[ \cardv (  C_{G,K_0} (0)) ^2 \right] < \infty$ (see for instance Theorem (6.75) in \cite{grimmettt:percolation}), thus $\EE[R_0^\ell] < \infty$ and the multiparameter ergodic theorems mentioned above apply to get \eqref{e:pif2} and \eqref{e:pif3}. Moreover, by a dominated convergence theorem, we obtain that 
$$ \lim_{\ell\rightarrow \infty} \EE \left[ \cardv (  C_{G,K_0} (0)) ^2  \mathds{1}_{ \cardv (  C_{G,K_0} (0)) > \ell/2 }\right]  \,=\, 0\,. $$
It is known that at criticality, there is no infinite cluster in the percolation in half space, see \cite{grimmettt:percolation}, Theorem (7.35). Thus $G(\{0\}) \geq 1-p_c(d)$ implies that 
$$ \lim_{\ell\rightarrow \infty} \PP\left[ 0 \overset{\HH}{\underset{G,0}{\longleftrightarrow}} \ZZ^{d-1} \times \{\ell/2 - 1\} \right]  \,=\, \PP \left[ 0 \textrm{ is connected to $\infty$ in }(\mathds{1}_{t_G(e)>0},e\in\HH) \right] \,=\, 0\,.$$
Thus for all $\eta >0$ we can choose $\ell ^\eta$ large enough so that for every $\ell\geq \ell^{\eta}$, $\EE[R_x^\ell]<\eta$. For every height function $h:\NN \mapsto \RR^+$ such that $\lim_{n\rightarrow \infty} h (n) =+\infty$, let $p_0$ be large enough such that for all $p\geq p_0$, $h(p) \geq \ell ^{\eta}$. The function $\ell\mapsto R_x^\ell$ is non-increasing, thus for every $D=\prod_{i=1}^{d-1} [0,k_i]$ ($k_i>0$) we have a.s. 
\begin{align*}
0 & \,\leq \, \limsup_{p\rightarrow \infty} \frac{1}{\cardv (pD \cap \ZZ^{d-1})} X_{pD}^{h(p)} \,\leq \, \limsup_{p\rightarrow \infty} \frac{1}{\cardv (pD \cap \ZZ^{d-1}) } X_{pD}^{\ell^{\eta}} \,=\, \EE [X^{\ell^\eta}] \,\leq\, \eta \,,
\end{align*}
thus 
\begin{equation}
\label{e:step2}
 \lim_{p\rightarrow \infty} \frac{X_{pD}^{h(p)} }{\H^{d-1}(pD\times \{0\})} \,=\, 0 \quad \textrm{a.s.} 
 \end{equation}
We turn back to the study of $\phi^{\v_0}_G (pA, h(p))$. We recall that we supposed $\lim_{p\rightarrow \infty} h(p) / \log p =+\infty$. As in the proof of Proposition \ref{p:finitude} (see \eqref{e:*}), since $G(]K_0,+\infty])<p_c(d)$, we have
\begin{align*}
 \sum_{p\in \NN} \PP [   \E_{G,K_0}( \cyl^{\v_0}(pA,h(p)),h(p))^c ] &\,\leq \, \sum_{p\in \NN}c_d \H^{d-1} (pA) h(p) \PP [\diam (C_{G,K_0} (0)) \geq h(p)]\\
 & \,\leq \, \sum_{p\in \NN}c_d \H^{d-1} (pA) h(p) \kappa_1 e^{-\kappa_2 h(p)} \\ 
& \,<\, +\infty \,,
 \end{align*}
thus by Borel-Cantelli we know that 
\begin{equation}
\label{e:step3}
\textrm{a.s., for all $p$ large enough,} \quad  \E_{G,K_0}( \cyl^{\v_0}(pA,h(p)),h(p)) \textrm{ occurs.}
\end{equation}
Proposition \ref{p:zero} is proved by combining  \eqref{e:pif1}, \eqref{e:step2} and \eqref{e:step3}.
\end{dem}

We now extend Proposition \ref{p:zero} to the study of any tilted cylinder. We will bound the maximal flow through a tilted cylinder by maximal flows through straight boxes at an intermediate level. Unfortunately, at this stage, we could not prove that the convergence holds almost surely. However, we prove that the convergence holds in a weaker sense, namely in probability. We will upgrade this convergence in Proposition \ref{p:zeroter}. 
\begin{prop}
\label{p:zerobis}
For any probability measure $G$ on $[0, +\infty]$ such that $G(\{+\infty\})<p_c(d)$ and $G(\{0\}) \geq 1-p_c(d)$, for any $\v \in \SS^{d-1}$, for any non-degenerate hyperrectangle $A$ normal to $\v$, for any function $h : \NN \mapsto \RR^+$ satisfying $\lim_{p\rightarrow +\infty} h(p) / \log p =+\infty$, we have
$$  \lim_{p\rightarrow \infty} \frac{\phi_G (pA, h(p))}{\H^{d-1} (pA)} \,=\, 0 \quad \textrm{in probability.}$$
\end{prop}

\begin{dem}
Fix $A$, $\v$ and $h$ and consider a $p$ large enough. Since $\phi(A,h)$ is non increasing in $h$, we can suppose that $h(p) \leq p$ for all $p$. We will bound $\phi_G (pA, h(p))$ by maximal flows through straight boxes at an intermediate level $L$, $1\leq L \leq p$ (in what follows $L$ will depend on $p$). For a fixed $L\in 2 \NN^*$, we chop $\ZZ^d$ into (almost) disjoint $L$-boxes as we already did in the proof of Proposition \ref{p:tronquer}. We recall the definitions of the $L$- and $3L$-boxes given in \eqref{e:Lambda} and \eqref{e:Lambda'}: let $\Lambda_L = [-L/2,L/2]^d $, for $\mathbf{i}\in \ZZ^d$ we have
$$ \Lambda_L (\mathbf{i}) \,=\, \{ x+L \mathbf{i} \,:\, x\in \Lambda_L \} \quad \textrm{and} \quad \Lambda'_L (\mathbf{i}) \,=\, \{ x+L \mathbf{i} \,:\, x\in \Lambda_{3L} \} \,.$$
For every $L\in \NN^*$, for every $\mathbf{i}\in \ZZ^d$, let $\phi_G (L,\mathbf{i})$ be the maximal flow from $\partial \Lambda_L (\mathbf{i}) \cap \ZZ^d$ to $\partial \Lambda'_L (\mathbf{i})  \cap \ZZ^d$ in $\Lambda'_L(\mathbf{i}) \smallsetminus \Lambda_L(\mathbf{i})$. By the max-flow min-cut Theorem, 
$$ \phi_G (L,\mathbf{i}) \,=\, \min \{ T_G (E) \,:\, E \subset \EE^d \,,\, E \textrm{ cuts $\partial \Lambda_L (\mathbf{i}) \cap \ZZ^d$ from $\partial \Lambda'_L (\mathbf{i})  \cap \ZZ^d$ in $\Lambda'_L(\mathbf{i}) \smallsetminus \Lambda_L(\mathbf{i})$} \} \,,$$
{\em i.e.}, roughly speaking, $ \phi_G (L,\mathbf{i})$ is the minimal capacity of a cutset in the annulus $\Lambda'_L (\mathbf{i}) \smallsetminus \Lambda_L (\mathbf{i})$. For every $\mathbf{i}\in \ZZ^d$, let $E_G(L,\mathbf{i})$ be a minimal cutset for $\phi_G (L,\mathbf{i})$. We choose $L =L(p)$ such that $h(p)$ is large in comparison with $L$, in the sense that no $3L$-box can intersect both the top and the bottom of $\cyl (pA, h(p)$). Thus we can choose $L(p) = 2 \lfloor h(p)/c_d \rfloor$ for some constant $c_d$ depending only on the dimension. Let $\mathbf{J}(pA, L)$ be the indices of all the $L$-boxes that are intersected by the hyperrectangle $pA$ (see Figure \ref{f:zerobis}):
$$ \mathbf{J} \,=\, \{ \mathbf{i}\in \ZZ^d \,:\, pA \cap \Lambda_L (\mathbf{i}) \neq \emptyset \}\,. $$
\begin{figure}[h!]
\centering
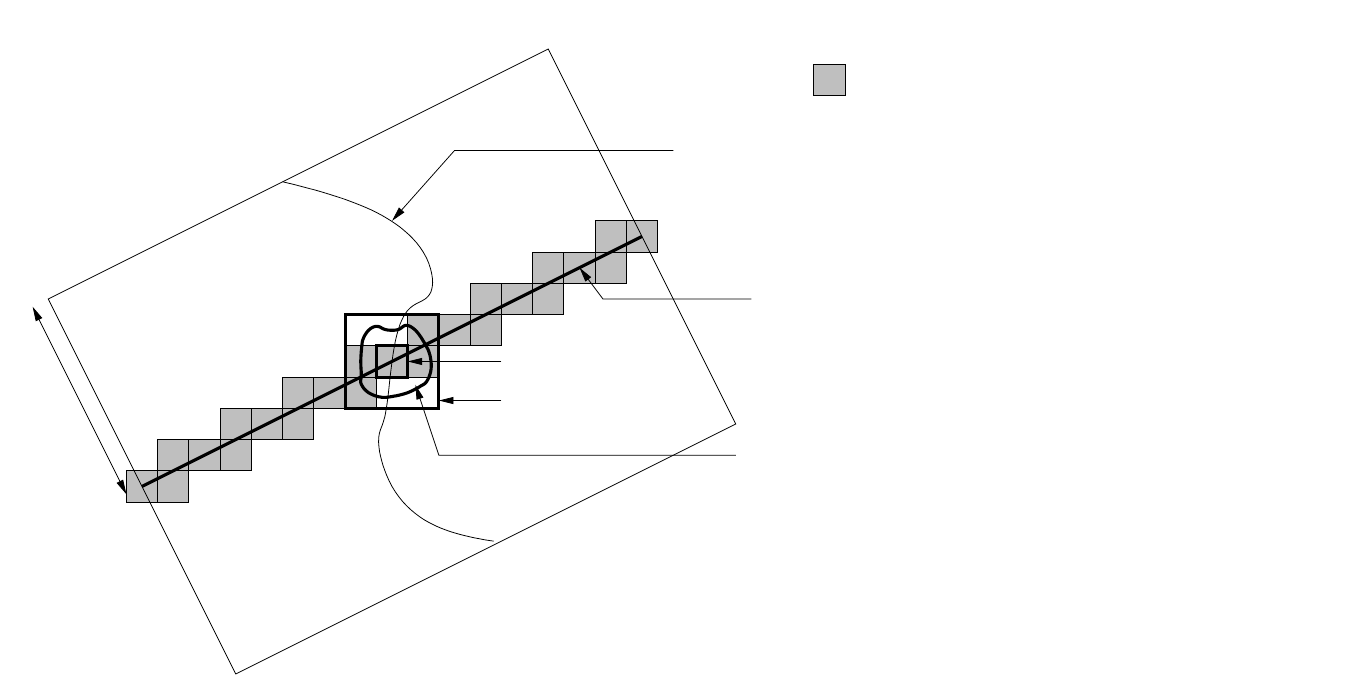
\caption{The cylinders $\cyl(pA,h(p)$ and $\Lambda_L (\mathbf{j}), \mathbf{j}\in \mathbf{J}$ ($d=2$).}
\label{f:zerobis}
\end{figure}

Let us prove that inequality \eqref{e:comp1} holds:
\begin{equation}
\label{e:comp1}
\phi_G (pA, h(p)) \,\leq \, \sum_{\mathbf{i}\in \mathbf{J}}  \phi_G (L,\mathbf{i})
\end{equation}
by proving that $\cup_{\mathbf{i}\in \mathbf{J}} E_G(L,\mathbf{i})  $ is a cutset for $\phi_G(pA, h(p))$. Let $\gamma = (x_0, e_1, x_1, \dots , e_n, x_n)$ be a path from the top to the bottom of $\phi_G(pA, h(p))$. Since $pA \subset \cup_{\mathbf{i}\in \mathbf{J}} \Lambda_L (\mathbf{i})  $ and $\gamma$ (seen as a continuous curve) must intersect $pA$, then there exists $\mathbf{i} \in J$ such that $\gamma \cap \Lambda_L (\mathbf{i}) \neq \emptyset $. Since $h(p)$ is large in comparison with $L$, $\gamma$ cannot be included in $ \Lambda'_L (\mathbf{i}) $, thus $\gamma$ contains a path from $\partial \Lambda_L (\mathbf{i}) \cap \ZZ^d$ to $\partial \Lambda'_L (\mathbf{i})  \cap \ZZ^d$ in $\Lambda'_L(\mathbf{i}) \smallsetminus \Lambda_L(\mathbf{i})$, thus by definition it must intersect $E_G(L,\mathbf{i})$ (see Figure \ref{f:zerobis}). This proves that $\cup_{\mathbf{i}\in \mathbf{J}} E_G(L,\mathbf{i})  $ is a cutset for $\phi_G(pA, h(p))$, thus
$$\phi_G (pA, h(p)) \,\leq \,  \sum_{\mathbf{i}\in \mathbf{J}} T_G(E_G(L,\mathbf{i}) )\,=\, \sum_{\mathbf{i}\in \mathbf{J}}  \phi_G (L,\mathbf{i})  \,.$$

It remains to compare $\phi_G (L,\mathbf{i})$ for any fixed $\mathbf{i}\in \ZZ^d$ with maximal flows through straight cylinders. Fix $\mathbf{i}\in \ZZ^d$. For every $k\in \{1,\dots , d\}$, define
$$ \widetilde{\Lambda}_L (\mathbf{i}, k , +) \,=\, \{ x+L\mathbf{i} \,:\, x\in [-3L/2, 3L/2]^{k-1} \times [L/2, 3L/2] \times [-3L/2, 3L/2]^{d-k} \}  $$
and
$$ \widetilde{\Lambda}_L (\mathbf{i}, k , -) \,=\, \{ x+L\mathbf{i} \,:\, x\in [-3L/2, 3L/2]^{k-1} \times [-3L/2, -L/2] \times [-3L/2, 3L/2]^{d-k} \}  \,,$$
see Figure \ref{f:zerobis2}.
\begin{figure}[h!]
\centering
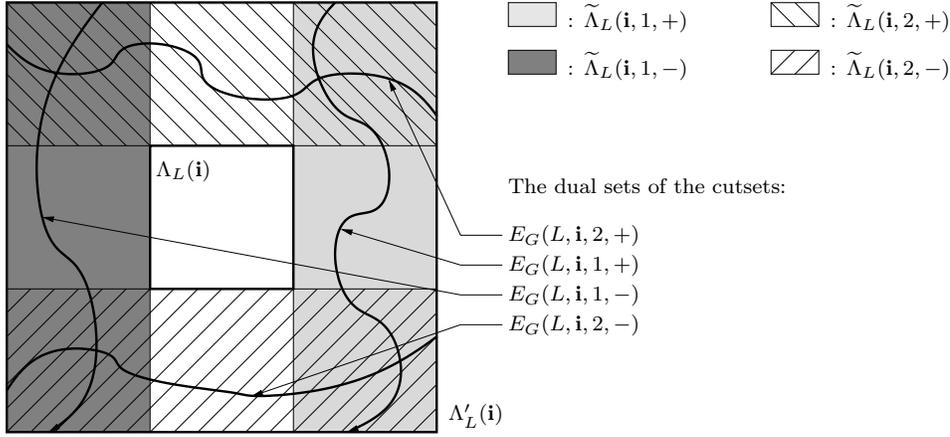
\caption{The cylinders $\Lambda_L (\mathbf{i}), \Lambda_L'(\mathbf{i})$ and $\widetilde{\Lambda}_L (\mathbf{i}, k , l)$ for $k\in \{1,2\}$ and $l\in\{+,-\}$ ($d=2$).}
\label{f:zerobis2}
\end{figure}
Let $\phi_G (L, \mathbf{i}, k , +)$ (rep. $\phi_G (L, \mathbf{i}, k , -)$) be the maximal flow in $ \widetilde{\Lambda}_L (\mathbf{i}, k , +)$ (resp. $ \widetilde{\Lambda}_L (\mathbf{i}, k , -)$) from its top $[-3L/2, 3L/2]^{k-1} \times\{ 3L/2 \} \times [-3L/2, 3L/2]^{d-k}\cap \ZZ^d$ (resp. $[-3L/2, 3L/2]^{k-1} \times\{ -3L/2 \} \times [-3L/2, 3L/2]^{d-k}\cap \ZZ^d$) to its bottom $[-3L/2, 3L/2]^{k-1} \times\{ L/2 \} \times [-3L/2, 3L/2]^{d-k}\cap \ZZ^d$ (resp. $[-3L/2, 3L/2]^{k-1} \times\{ -L/2 \} \times [-3L/2, 3L/2]^{d-k}\cap \ZZ^d$), and let $E_G (L, \mathbf{i}, k , +)$ (rep. $E_G (L, \mathbf{i}, k , -)$) be a corresponding minimal cutset.

We claim that for every $L\in \NN^*$, for every $\mathbf{i}\in \ZZ^d$, $\cup_{l=+,-} \cup_{k=1}^d E_G (L, \mathbf{i}, k , l)$ cuts $\partial (\Lambda_L (\mathbf{i})) \cap \ZZ^d$ from $\partial (\Lambda'_L (\mathbf{i}) ) \cap \ZZ^d$ in $\Lambda'_L(\mathbf{i}) \smallsetminus \Lambda_L(\mathbf{i})$, we have
\begin{equation}
\label{e:comp2}
\phi_G (L,\mathbf{i})  \,\leq\, \sum_{l=+,-} \sum_{k=1}^d \phi_G (L, \mathbf{i}, k , l) \,.
\end{equation}
We now prove this claim. Let $\gamma = (x_0, e_1, x_1, \dots , e_n, x_n)$ be a path from $\partial (\Lambda_L (\mathbf{i})) \cap \ZZ^d$ to $\partial (\Lambda'_L (\mathbf{i}) ) \cap \ZZ^d$ in $\Lambda'_L(\mathbf{i}) \smallsetminus \Lambda_L(\mathbf{i})$. Let 
$$ j \,=\, \inf \{ i\in \{0,\dots , n \} \,:\, x_i  \in \partial (\Lambda'_L(\mathbf{i})) \} \,. $$
Then there exists $k\in \{1,\dots, d\}$, $l\in\{+,-\}$ such that
$$ x_{j} \,=\,  [-3L/2, 3L/2]^{k-1} \times\{ l 3L/2 \} \times [-3L/2, 3L/2]^{d-k} \,,$$
thus $x_{j} \in  \widetilde{\Lambda}_L (\mathbf{i}, k , l)$. Let
$$ j' \,=\, \inf \{ i\leq j \,:\, \forall i' \in \{i,\dots , j \} \,,\, x_{i'} \in   \widetilde{\Lambda}_L (\mathbf{i}, k , l) \} \,.$$
Then by continuity of $\gamma$ we know that $x_{j'}\in \widetilde{\Lambda}_L (\mathbf{i}, k , l)$ but $x_{j'}$ has a neighbor outside $\widetilde{\Lambda}_L (\mathbf{i}, k , l)$. By definition of $j$, since $j'<j$, $x_{j'}$ can be only on one side of the boundary of $\widetilde{\Lambda}_L (\mathbf{i}, k , l)$, precisely $x_{j'}\in [-3L/2, 3L/2]^{k-1} \times\{ l L/2 \} \times [-3L/2, 3L/2]^{d-k}$. Thus the subset of $\gamma$ between $x_{j'}$ and $x_{j}$ is a path from the bottom to the top of $\widetilde{\Lambda}_L (\mathbf{i}, k , l)$, thus it must intersect $E_G (L, \mathbf{i}, k , l)$. This proves that $\cup_{l=+,-} \cup_{k=1}^d E_G (L, \mathbf{i}, k , l)$ cuts $\partial (\Lambda_L (\mathbf{i})) \cap \ZZ^d$ from $\partial (\Lambda'_L (\mathbf{i}) ) \cap \ZZ^d$ in $\Lambda'_L(\mathbf{i}) \smallsetminus \Lambda_L(\mathbf{i})$, thus \eqref{e:comp2} is proved.

Combining \eqref{e:comp1} and \eqref{e:comp2}, we obtain that for every $L$, for every $p$ large enough, 
\begin{equation}
\label{e:comp3}
\phi_G (pA, h(p)) \,\leq \, \sum_{\mathbf{i}\in \mathbf{J}} \sum_{l=+,-} \sum_{k=1}^d \phi_G (L, \mathbf{i}, k , l) \,.
\end{equation}
For short, we denote by $\mathcal E (L,\mathbf i, k, l)$ the event $\E_{G,K_0}(\widetilde{\Lambda}_L (\mathbf{i}, k, l), L  )$ defined by
$$ \mathcal E (L,\mathbf i, k, l) := \bigcap_{x\in \widetilde{\Lambda}_L (\mathbf{i}, k, l) \cap \ZZ^d } \{ \diam (C_{G,K_0} (x)) < L \} $$
for any $\mathbf{i}\in \mathbf{J}$, $k\in \{1,\dots, d\}$ and $l\in\{+,-\}$. On one hand, by symmetry and invariance of the model by translations of integer coordinates, we have, for any such $(\mathbf i, k, l)$,
$$ \EE [ \phi_G (L, \mathbf{i}, k , l) \mathds{1}_{\mathcal E (L,\mathbf i, k, l)}]  \,=\, \EE [ \phi_G (L, \mathbf{0}, d , +) \mathds{1}_{\mathcal E (L,\mathbf 0, d, +)}] \,=\, \EE [ \phi^{\v_0}_G (LD,L ) \mathds{1}_{\mathcal E_{G,K_0} (\cyl^{\v_0}(LD,L),L)}] $$
with $D = [0, 3]^{d-1} \times \{0\}$ and $\phi^{\v_0}$ defined as in Proposition \ref{p:zero}. By Proposition \ref{p:zero} we know that
\begin{equation}
\label{e:comp4}
\lim_{L \rightarrow \infty } \frac{\EE [ \phi^{\v_0}_G (LD,L ) \mathds{1}_{\mathcal E_{G,K_0} (\cyl^{\v_0}(LD,L),L)}] }{L^{d-1}} \,=\, 0 \,.
\end{equation}
On the other hand, let $A^1$ by a hyperrectangle a bit larger than $A$, namely
$$ A^1 \,=\, \{ x+\vec w \,:\, x\in A \,,\, \|\vec w \|_{\infty} \leq 1 \,,\, \vec w \cdot \v =0 \} \,. $$
We recall that the event $\E_{G,K_0 }(\cyl(pA^1,h(p)), L)$ is defined by
$$ \E_{G,K_0 }(\cyl(pA^1,h(p)), L) \,=\, \bigcap_{x\in \cyl  (pA^1,h(p)) \cap \ZZ^d } \{ \diam (C_{G,K_0} (x)) <  L \} \,.$$
Notice that for all $\mathbf{i}\in \mathbf{J}$, we have $\Lambda'_L(\mathbf{i}) \subset \cyl  (pA^1,h(p))$ at least for $p$ large enough since we choose $L=L(p) = 2 \lfloor h(p)/c_d \rfloor$ for some constant $c_d$ depending only on the dimension. Then 
$$ \E_{G,K_0 }(\cyl(pA^1,h(p)), L)\,\subset \,     \bigcap_{\mathbf{i}\in \mathbf{J}} \bigcap_{l=+,-} \bigcap_{k=1}^d  \E(L,\mathbf i, k, l) \,.$$
By \eqref{e:comp3} we obtain
\begin{align*}
\frac{\EE [\phi_G (pA, h(p)) \mathds{1}_{\E_{G,K_0 }(\cyl(pA^1,h(p)), L) }]} {p^{d-1}}
& \,\leq\, \frac{ \EE [\phi_G (pA, h(p)) \mathds{1}_{ \cap_{\mathbf{i}\in \mathbf{J}} \cap_{l=+,-} \cap_{k=1}^d  \E(L,\mathbf i, k, l)}]}{p^{d-1}} \\
& \,\leq\, \frac{1}{p^{d-1}}  \sum_{\mathbf{i}\in \mathbf{J}} \sum_{l=+,-} \sum_{k=1}^d \EE[ \phi_G (L, \mathbf{i}, k , l) \mathds{1}_{ \E(L,\mathbf i, k, l)} ]\\
& \,\leq\, \frac{2d L^{d-1} \card (\mathbf{J})}{p^{d-1}}  \frac{\EE [ \phi^{\v_0}_G (LD,L ) \mathds{1}_{\mathcal E_{G,K_0} (\cyl^{\v_0}(LD,L),L)}] }{L^{d-1}}\,,
\end{align*}
and it remains to notice that $L=L(p)$ goes to infinity and $\card (\mathbf{J}) L(p)^{d-1} / p^{d-1}$ remains bounded when $p$ goes to infinity to conclude by \eqref{e:comp4} that
\begin{equation}
\label{e:comp5}
\lim_{p\rightarrow \infty} \frac{\EE [\phi_G (pA, h(p)) \mathds{1}_{ \E_{G,K_0 }(\cyl (pA^1,h(p)), L) } ]}{p^{d-1}} \,=\, 0 \,.
\end{equation}
For every $\eta >0$, we obtain as in \eqref{e:*} that
\begin{align*}
\PP [\phi_G (pA, h(p))&  \geq \eta p^{d-1}] \\
&\,\leq\, \PP [\E_{G,K_0 }(\cyl(pA^1,h(p)), L) ^c] + \PP [\phi_G (pA, h(p))  \mathds{1}_{ \E_{G,K_0 }(\cyl(pA^1,h(p)), L) }  \geq \eta p^{d-1}] \\
&\,\leq\, c_d \H^{d-1} (pA^1) h(p) \kappa_1 e^{-\kappa_2 L(p)}+ \eta^{-1} \frac{\EE [\phi_G (pA, h(p))  \mathds{1}_{ \E_{G,K_0 }(\cyl(pA^1,h(p)), L) } ]}{p^{d-1}}
\end{align*}
that goes to zero when $p$ goes to infinity since $h(p) \leq p $, $L(p) = 2 \lfloor h(p)/c_d \rfloor$ and $\lim_{p\rightarrow \infty} h(p)/\log (p) =+\infty$.
\end{dem}


\subsection[Proof of the convergence III]{Proof of the convergence III: end of the proof of Theorem \ref{t:CV}}
\label{s:retourzero} 

At this stage, what remains to prove to finish the proof of Theorem \ref{t:CV} is to strengthen the mode of convergence in Proposition \ref{p:zerobis}. This can be done easily using the continuity of $G\mapsto \nu_G$, {\em i.e.}, Theorem \ref{thmcont}. 

\begin{prop}
\label{p:zeroter}
We suppose that $G\mapsto \nu_G(\v)$ is continuous, {\em i.e.}, if $G$ and $(G_n)_{n\in \NN}$ are probability measures on $[0,+\infty]$ such that $G(\{+\infty\}) < p_c(d)$ and for all $n\in \NN$, $G_n (\{+\infty\}) < p_c(d)$ and $G_n \overset{d}{\rightarrow} G$, then
$$   \lim_{n\rightarrow \infty} \sup_{\v \in \SS^{d-1}} \left\vert \nu_{G_n} (\vec v) - \nu_G (\vec v) \right\vert \,=\, 0 \,.$$
For any probability measure $G$ on $[0, +\infty]$ such that $G(\{+\infty\})<p_c(d)$ and $G(\{0\}) \geq 1-p_c(d)$, for any $\v \in \SS^{d-1}$, for any non-degenerate hyperrectangle $A$ normal to $\v$, for any function $h : \NN \mapsto \RR^+$ satisfying $\lim_{p \rightarrow +\infty} h(p) / \log p =+\infty$, we have
$$  \lim_{p\rightarrow \infty} \frac{\phi_G (pA, h(p))}{\H^{d-1} (pA)} \,=\, 0 \quad \textrm{a.s.}$$
\end{prop}

\begin{dem}
Let $\v\in \SS^{d-1}$. Let $G$ be a probability measure on $[0,+\infty]$ such that $G(\{+\infty\})<p_c(d)$ and $G(\{ 0 \}) \geq 1-p_c(d)$. By Proposition \ref{p:trivial}, this implies that $\nu_G(\v) =0$.  let $A$ be a non-degenerate hyperrectangle normal to $\v$, and let $h:\NN^* \mapsto \RR^+$ such that $\lim_{p\rightarrow \infty} h(p)/\log p =+\infty$. 
Suppose first that $h$ also satisfies $\lim_{p\rightarrow \infty} h(p)/ p =0$.

For any $\eps >0$, we denote by $^\eps G $ the distribution of the variables $t_G(e) + \eps$. Obviously $^\eps G (\{0\}) = 0$ thus Proposition \ref{p:positif} states that
$$ \lim_{p\rightarrow \infty} \frac{\phi_{^\eps G}(pA, h(p))}{\H^{d-1} (pA)}  \,=\, \nu_{^\eps G} (\v) \quad \textrm{a.s.}$$
Moreover, by coupling (see Equation \eqref{e:couplage}), $ \phi_{G}(pA, h(p)) \leq \phi_{^\eps G}(pA, h(p))$, thus
\begin{equation}
\label{e:casnul}
\forall \eps >0 \quad  \limsup_{p\rightarrow \infty} \frac{\phi_{G}(pA, h(p))}{\H^{d-1} (pA)}  \,\leq\, \nu_{^\eps G} (\v) \quad \textrm{a.s.}
\end{equation}
To conclude the proof, we will use the continuity of $G\mapsto \nu_G$: since $^\eps G \overset{d}{\longrightarrow} G$ when $\eps $ goes to $0$ we obtain that 
\begin{equation}
\label{e:casnul2}
\lim_{\eps \rightarrow 0} \nu_{^\eps G} (\v) \,=\, \nu_G(\v) \,=\, 0 \,.
\end{equation}
Combining \eqref{e:casnul} and \eqref{e:casnul2} we obtain that
$$ \lim_{p\rightarrow \infty} \frac{\phi_{G}(pA, h(p))}{\H^{d-1} (pA)} \,=\, 0 \qquad \textrm{a.s.}$$

If $h$ does not satisfy $\lim_{p\rightarrow \infty} h(p)/ p =0$, define $\tilde{h} (p) = \min (h(p), \sqrt{p})$ for all $p\in \NN^*$. Then $\tilde h$ is mild, thus we just proved that
$$ \lim_{p\rightarrow \infty} \frac{\phi_{G}(pA, \tilde{h}(p))}{\H^{d-1} (pA)} \,=\, 0 \qquad \textrm{a.s.}$$
Moreover, since $\tilde{h} (p) \leq h(p)$ for all $p\in \NN^*$, any cutset from the top to the bottom of $\cyl(pA, \tilde{h}(p))$ is also a cutset from the top to the bottom of $\cyl(pA, h(p))$, thus by the max-flow min-cut Theorem we obtain that $\phi_{G}(pA, \tilde{h}(p)) \geq \phi_{G}(pA, h(p))$ for all $p\in \NN^*$. This allows us to conclude that 
$$ \limsup_{p\rightarrow \infty} \frac{\phi_{G}(pA, h(p))}{\H^{d-1} (pA)} \,\leq\, \lim_{p\rightarrow \infty} \frac{\phi_{G}(pA, \tilde{h}(p))}{\H^{d-1} (pA)} \,=\, 0 \qquad \textrm{a.s.}$$
This ends the proof of Proposition \ref{p:zeroter}.
\end{dem}

\begin{rem}
It is worth noticing that this proof does not use Propositions \ref{p:zero} or Proposition \ref{p:zerobis} directly. However, we need these intermediate results to prove the continuity of $G\mapsto \nu_G$ that we use here.
\end{rem}


\section{Subadditivity}
\label{s:ssadd}

As mentioned in section~\ref{s:T}, expressing the flow constant as the limit of a subadditive and integrable object is crucial to prove its continuity. This is the purpose of the present section. The first idea is to take the capacity of a cut which in a sense separates a hyperrectangle $A$ from infinity in a half-space. This will ensure subadditivity. However, in order to have a chance to compare it to the flows that we used so far, one needs the cut to stay at a small enough distance from $A$ so that it will be flat in the limit. In addition, to ensure good integrability properties, one needs this distance to be large enough so that one may find enough edges with bounded capacity to form a cutset. These constraints lead to searching for a cutset in a \emph{slab} of random height, which height is defined in \eqref{e:defH}.

\begin{figure}[h!]
\centering
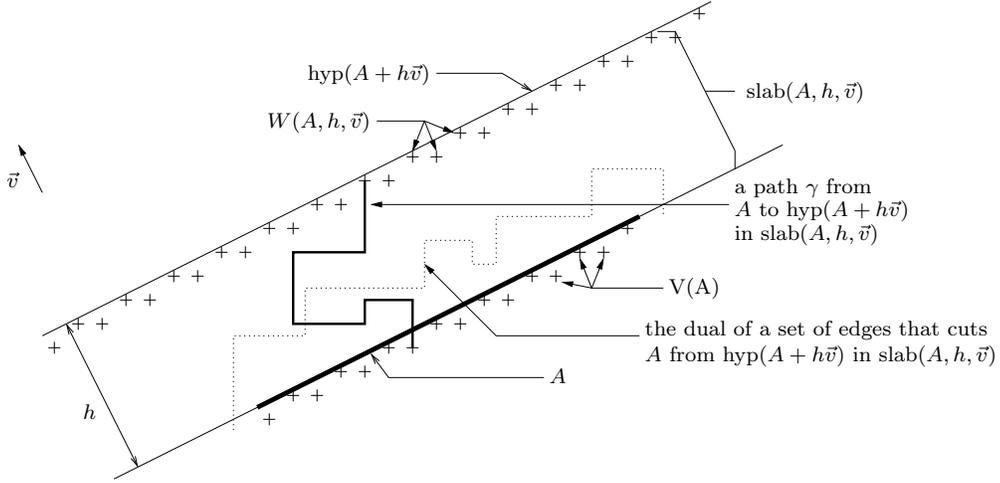
\caption{A path $\gamma $ from $A$ to $\hyp (A+h\v)$ in $\slab (A,h,\v)$ and a corresponding cutset ($d=2$).}
\label{f:sub1}
\end{figure}

Let $\v \in \SS^{d-1}$, and let $A$ be any non-degenerate hyperrectangle normal to $\v$. For any $h$, we denote by $\slab (A,h, \v)$ the cylinder whose base is the hyperplane spanned by $A$ and of height $h$ (possibly infinite), {\em i.e.}, the subset of $\RR^d$ defined by
$$ \slab (A,h,\v) \,=\, \{ x + r \v \,:\, x \in \hyp (A)  \,,\, r\in [0,h] \} \,.$$
Let $V(A)$ be the following set of vertices in $\ZZ^d$, which is a discretized version of $A$ :
$$ V(A) \,=\, \{ x\in \ZZ^d \cap \slab (A,\infty,\v)^c \,:\, \exists y \in \ZZ^d \cap \slab (A,\infty,\v) \,,\, \langle x,y \rangle \in \EE^d \textrm{ and $\langle x,y \rangle$ intersects $A$}\}\,. $$
Let $W(A,h, \v)$ be the following set of vertices in $\ZZ^d$, which is a discretized version of $\hyp (A + h\v)$ :
$$ W(A,h, \v) \,=\, \{ x\in \ZZ^d \cap \slab (A,h,\v) \,:\, \exists y \in \ZZ^d \cap ( \slab (A,\infty,\v) \smallsetminus \slab (A,h,\v)) \,,\, \langle x,y \rangle \in \EE^d \}\,. $$
We say that a path $\gamma = (x_0, e_1, x_1, \dots, e_n, x_n)$ goes from $A$ to $\hyp (A+h\v)$ in $\slab (A,h,\v)$ if $x_0 \in V(A)$, $x_n \in W(A,h, \v)$ and for all $ i\in\{1,\dots , n\}$, $x_i \in \slab (A,h,\v)$ (see Figure \ref{f:sub1}).
We say that a set of edges $E$ cuts $A$ from $\hyp (A+h\v)$ in $\slab (A,h,\v)$ if $E$ contains at least one edge of any path $\gamma $ that goes from $A$ to $\hyp (A+h\v)$ in $\slab (A,h,\v)$.

For any probability measure $F$ on $[0,+\infty]$ such that $F(\{+\infty\}) < p_c(d)$, for any $K_0\in \RR$ such that $F([K_0, +\infty])<p_c(d)$, we define the random height $H_{F,K_0}(A)$ as
\begin{equation}
\label{e:defH}
 H_{F,K_0}(A) \,=\, \inf \left\{ h\geq \H^{d-1} (A) ^{\frac{1}{2(d-1)}} \,:\, \begin{array}{c}\exists E \subset \EE^d \textrm{ s.t. } \forall e\in E \,,\, t_F(e) \leq K_0\\ \text{ and $E$ cuts $A$ from $\hyp(A+h\v)$}\\ \text{ in $\slab(A,h,\v)$} \end{array}\right\}\,.
 \end{equation}
We will say a few words about the definition of $H_{F,K_0} (A)$ in Remark \ref{r:H} after the proof of the first result of this section, namely Theorem \ref{t:ssadd}. We finally define the random alternative maximal flow $\wphi_{G, F,K_0} (A)$ by
\begin{equation}
\label{e:varphi}
\wphi_{G,F,K_0} (A) \,=\, \inf \left\{T_G(E) \,:\, \begin{array}{c} E\subset \EE^d \textrm{ and $E$ cuts $A$ from $\hyp(A+H_{F,K_0}(A) \v)$}\\ \textrm{in $\slab(A,H_{F,K_0} (A),\v)$ }  \end{array} \right\} \,.
\end{equation}
Notice that we do not know if the infimum in the definition \eqref{e:varphi} of $\wphi_{G,F,K_0} (A)$ is achieved. The purpose of using two different distributions $F$ and $G$ in the definition of $ \wphi_{G,F,K_0} (A)$ is to have monotonicity in $G$, which will be used later, in the proof of Proposition~\ref{propupper}.

Finally, we say that a direction $\v \in \SS^{d-1}$ is rational if there exists $M\in \RR^+$ such that $M\v$ has rational coordinates.

Now, we will prove that these flows $\wphi_{G,F,K_0} (A)$ properly rescaled converge for large hyperplanes towards $\nu_G(\v)$ (as defined in \eqref{defnu}), and thus obtain an alternative definition of $\nu_G(\v)$. This will be done in two steps. First we prove the convergence of $\wphi_{G,F,K_0} (pA)/\H^{d-1} (pA)$ towards some limit $\wnu_{G,F,K_0} (\v)$ by some subadditive argument in Theorem~\ref{t:ssadd}. Then, we compare $\wphi_{G,F,K_0} (pA)$ with $\phi_{G} (pA, h(p))$ to prove that $\wnu_{G,F,K_0}(\v) = \nu_G (\v)$, and this is done in Proposition~\ref{prop:ssadd}.
\begin{thm}
\label{t:ssadd}
Let $G$ be a probability measure on $[0,+\infty]$ such that $G(\{+\infty\}) < p_c(d)$. For any probability measure $F$ on $[0,+\infty]$ such that $F(\{+\infty\}) < p_c(d)$ and $G\preceq F$, for any $K_0\in \RR$ such that $F(]K_0, +\infty])<p_c(d)$, for any rational $\v \in \SS^{d-1}$, there exists a non-degenerate hyperrectangle $A$ (depending on $\v$ but neither on $G, F$ nor on $K_0$) which is normal to $\v$ and contains the origin of the graph $\mathbb{Z}^d$ such that
$$\wnu_{G,F,K_0}(\v):= \inf_{p\in\NN^*} \frac{\EE[\wphi_{G, F,K_0} (pA)]}{\H^{d-1} (pA)}<\infty$$
and
$$  \lim_{p\rightarrow \infty} \frac{\wphi_{G,F,K_0} (pA)}{\H^{d-1} (pA)} \,=\,\wnu_{G,F,K_0}(\v)\quad \text{a.s. and in }L^1\;.$$
\end{thm}
\begin{dem}
  Let $G$ be a probability measure on $[0,+\infty]$ such that $G(\{+\infty\}) < p_c(d)$. Let $F$ be a probability measure on $[0,+\infty]$ such that $F(\{+\infty\}) < p_c(d)$ and $G\preceq F$. Let $K_0\in \RR$ such that $F(]K_0, +\infty])<p_c(d)$. We consider a fixed rational $\v \in \SS^{d-1}$ and $H$, the hyperplane normal to $\v$ containing $0$. Since $\v$ is rational, there exists a orthogonal basis of $H$ of vectors with integer coordinates, let us call it $(\vec{f}_1,\ldots,\vec{f}_{d-1})$. Then, we take $A$ to be the hyperrectangle built on the origin and this basis: $A=\{\sum_{i=1}^{d-1}\lambda_i\vec{f}_i\;:\;\forall i, \;\lambda_i\in [0,1]\}$. Notice that the model is invariant under translations by $\vec{f}_i$ for any $i$, in the sense that the flow $\wphi_{G,F,K_0} (A+\vec{f}_i)$ with capacities $t$ is equal to the flow $\wphi_{G,F,K_0} (A)$ with capacities $t'$ defined by  $t'(\langle x,y\rangle)=t(\langle x+\vec{f}_i,y+\vec{f}_i\rangle)$. Moreover, if $A_1,\dots, A_k$ are hyperrectangles included in $\hyp (A)$ with disjoint interiors and such that $B=\cup_{i=1}^k A_i$ is also an hyperrectangle, we claim that
\begin{equation}
\label{e:ssadd}
 \wphi_{G,F,K_0} (B) \,\leq\, \sum_{i=1}^k \wphi_{G,F,K_0} (A_i) \,.
 \end{equation}
Indeed, first notice that if $B_1, B_2$ are hyperrectangles normal to $\v$ such that $B_1 \subset B_2$, then by definition any set of edges $E$ that cuts $B_2$ from $\hyp (B_2 + h\v)$ in $\slab (B_2, h \v)$ also cuts $B_1$ from $\hyp (B_1 + h\v) = \hyp (B_2 + h\v)$ in $\slab (B_1, h \v)  = \slab (B_2, h \v)$, thus $H_{F,K_0}(B_2) \geq H_{F,K_0} (B_1)$. Thus if $B=\cup_{i=1}^k A_i$ then $H_{F,K_0}(B) \geq \max_{1\leq i \leq k} H_{F,K_0} (A_i)$. For all $i\in \{1,\dots , k\}$,  let $E_i$ be a set of edges that cuts $A_i$ from $\hyp (A_i + H_{F,K_0} (A_i) \v)$ in $\slab (A_i, H_{F,K_0} (A_j) , \v)$. Let us prove that $\cup_{i=1}^k E_i$ cuts $B$ from $\hyp (B + H_{F,K_0} (B) \v)$ in $\slab (B, H_{F,K_0} (B) , \v)$. Let $\gamma =(x_0, e_1, x_1, \dots, e_n, x_n) $  be a path from $B$ to $\hyp (B + H_{F,K_0} (B) \v)$ in $\slab (B,H_{F,K_0}(B),\v)$. Since $B=\cup_{i=1}^k A_i$, we have $V(B) = \cup_{i=1}^k V(A_i)$ thus $x_0 \in V(A_j)$ for some $j\in \{1,\dots , k\}$. If $x_n\in W(A_j,H_{F,K_0} (A_j), \v) $ thus $\gamma$ is a path from $A_j$ to $\hyp (A_j + H_{F,K_0} (A_j) \v)$ in $\slab (A_j, H_{F,K_0} (A_j) , \v)$, thus $\gamma$ contains an edge of $E_j$. Otherwise since $H_{F,K_0} (A_j) \leq H_{F,K_0} (B)$, we know that $m:=\inf\{ p\in\{1,\dots ,n\} \,:\, x_p \notin \slab (A_j, H_{F,K_0} (A_j), \v)\} \leq n$, thus $\gamma$ contains a path $(x_0, e_1, x_1, \dots, x_{m-1})$ from $A_j$ to $\hyp (A_j + H_{F,K_0} (A_j) \v)$ in $\slab (A_j, H_{F,K_0} (A_j) , \v)$ for some $j\in \{1,\dots , k\}$, and we conclude also that $\gamma$ contains an edge of $E_j$. Inequality \eqref{e:ssadd} follows by optimizing on $T_G(E_i)$ for all $i$.

We now prove that $\wphi_{G,F,K_0} (A)$ has good integrability properties. For any $x\in V(A)$ we consider the connected component of $x$ in $\slab(A, \infty, \v)$ for the percolation $(\mathds{1}_{t_F(e)>K_0})$, {\em i.e.},
$$ C^{\v}_{F,K_0} (x) \,=\, \left\{ y\in \slab(A, \infty, \v) \,:\,\begin{array}{c} \textrm{$y$ is connected to $x$ by a path $\gamma = (x_0, e_0,\dots , x_n)$}\\ \textrm{s.t. $\forall i\in \{1,\dots, n\}$, $x_i \in \slab(A, \infty, \v) $ and $t_F(e_i) > K_0$}\end{array}  \right\}\,.$$
By definition of $H_{F,K_0} (A)$, we know that any path $\gamma$  in $ \slab(A,H_{F,K_0} (A),\v)$ from $A$ to $\hyp (A + H_{F,K_0} (A)\v)$ must contain at least one edge $e$ such that $t_F(e)\leq K_0$. Thus $\gamma$ cannot be included in $\cup_{x\in V(A)} C^{\v}_{F,K_0} (x)$, and this implies that $\gamma$ must contain at least one edge $e$ that belongs to $\pe (\cup_{x\in V(A)} C^{\v}_{F,K_0} (x)  )$. This edge $e$ satisfies (by the coupling relation \eqref{e:couplage}) $t_G(e)\leq t_F(e) \leq K_0$. Thus comparing clusters in the slab with clusters in the full space, we obtain
\begin{align*}
\EE[ \wphi_{G,F,K_0} (A)] & \,\leq \, K_0 \sum_{x\in V(A)} \EE[ \card_e (\pe  C^{\v}_{F,K_0} (x) )] \,\leq \, c_d K_0 \sum_{x\in V(A)}\EE[ \card_v ( C^{\v}_{F,K_0} (x) )]\\
 &  \,\leq \, c_d K_0 \sum_{x\in V(A)} \EE[\card_v ( C_{F,K_0} (x) )]  \,\leq \, c_d K_0 \card_v(V(A)) \EE[\card_v ( C_{F,K_0} (0) )] \,,
 \end{align*}
and
$$ \frac{\EE[ \wphi_{G,F,K_0} (A)] }{\H^{d-1} (A)} \,\leq \,  c_d K_0 \frac{\card_v(V(A))}{\H^{d-1} (A)} \EE[\card_v ( C_{F,K_0} (0) )]  \,\leq \,c'_d K_0 \EE[\card_v ( C_{F,K_0} (0) )] \,<\, +\infty $$
uniformly in $A$.

We can thus apply a multi-parameter ergodic theorem (see for instance Theorem~2.4 in \cite{Akcoglu} and Theorem~1.1 in \cite{Smythe}) to deduce that there exists a constant $\wnu_{G,F,K_0}(\v)$ (that depends on $\v$ but not on $A$ itself) such that
$$ \wnu_{G,F,K_0}(\v) \,=\, \inf_{p\in\NN^*} \frac{\EE[\wphi_{G,F,K_0} (pA)]}{\H^{d-1} (pA)}\,=\,  \lim_{p\rightarrow \infty} \frac{\wphi_{G,F,K_0} (pA)}{\H^{d-1} (pA)} \quad\textrm{a.s. and in }L^1 \,. $$
\end{dem}

We now state that the limit $\wnu_{G,F,K_0}(\v)$ appearing in Theorem \ref{t:ssadd} is in fact equal to $\nu_G(\v)$. We want to clarify the fact that in the proof of Proposition~\ref{prop:ssadd} below, we will use the convergence in probability of rescaled flows in tilted cylinders towards the flow-constant, stated above in Propositions~\ref{p:positif} and~\ref{p:zerobis}.

\begin{rem}
\label{r:H}
We will see in the proof of Proposition~\ref{prop:ssadd} below that with large probability, $H_{F,K_0}(pA)$ equals $\H^{d-1} (pA) ^{\frac{1}{2(d-1)}} $ for $p$ large. Since $\wphi_{G,F,K_0} (pA)$ depends on $F$ only through $\H^{d-1} (pA) ^{\frac{1}{2(d-1)}} $, it is thus natural that the limit in Theorem~\ref{t:ssadd} does not depend on $F$. Moreover, notice that the function $h:p\mapsto \H^{d-1} (pA) ^{\frac{1}{2(d-1)}}$ is mild: this is why we chose to make appear the lower bound $\H^{d-1} (pA) ^{\frac{1}{2(d-1)}}$ in the definition \eqref{e:defH} of $H_{F,K_0}(pA)$.
\end{rem}

\begin{prop}
\label{prop:ssadd}
For any fixed rational $\v \in \SS^{d-1}$, any probability measure $F$ on $[0,+\infty]$ such that $F(\{+\infty\}) < p_c(d)$ and $G\preceq F$, and any $K_0\in \RR$ such that $F(]K_0, +\infty])<p_c(d)$,
$$\wnu_{G,F,K_0}(\v) = \nu_G(\v)\;.$$
\end{prop}
\begin{dem}
We first prove that $\nu_G(\v) \leq \wnu_{G,F,K_0} (\v)$. We associate with a fixed rational $\v \in \SS^{d-1}$ the same hyperrectangle $A$ as in the proof of Theorem \ref{t:ssadd}. 
We consider the function $h(p) =  \H^{d-1} (pA) ^{\frac{1}{2(d-1)}} $. Then $h$ is mild. Thus we can apply Propositions \ref{p:positif} or \ref{p:zerobis} to state that
\begin{equation}
\label{e:sens11}
 \lim_{p\rightarrow \infty} \frac{\phi_G (pA, h(p))}{ \H^{d-1} (pA)} \,=\, \nu_G(\v) \quad \textrm{in probability.} 
 \end{equation}
Moreover, let $\gamma = (x_0, e_0, \dots, e_n, x_n)$ be a path from the bottom to the top of $\cyl(pA, h(p))$ inside $\cyl(pA,h(p))$. Let $k=\max \{ j\geq 0 \,:\, x_j \notin \slab (pA,\infty,\v) \}$. Then $x_j \in V(pA)$ and the truncated path $\gamma' = (x_k, e_{k+1}, \dots , x_n)$ is a path from $pA$ to $\hyp (pA+h(p) \v)$ in $\slab (pA,h(p),\v)$. On the event $\{H_{F,K_0} (pA) = h(p)\}$, we conclude that any set of edges $E$ that cuts $pA$ from $\hyp(pA+H_{F,K_0}(pA) \v)$ in $\slab(pA,H_{F,K_0} (pA),\v)$ also cuts any path from the bottom to the top of $\cyl (pA, h(p))$, thus on the event $\{H_{F,K_0} (pA) = h(p)\} $ we have 
\begin{equation}
\label{e:sens12}
 \phi_G (pA, h(p)) \leq \wphi_{G,F,K_0}(pA) \,. 
\end{equation}
Combining \eqref{e:sens12} and \eqref{e:*} we obtain
\begin{align}
\label{e:sens13}
 \PP[ \phi_G (pA, h(p)) > \wphi_{G,F,K_0}(pA) ]& \,\leq \, \PP [H_{F,K_0} (pA) > h(p) ] \nonumber \\
 & \,\leq\, \PP[ \exists x\in V(pA) \,:\, \diam (C_{F,K_0} (x) ) \geq h(p)/2 ] \nonumber \\
 & \,\leq \, \card_v (V(pA)) \PP[\diam (C_{F,K_0} (0) ) \geq h(p)/2] \nonumber \\
 & \,\leq \, c_d \H^{d-1} (pA) \kappa_1 e^{-\kappa_2 h(p)}  
 \end{align}
that goes to zero when $p$ goes to infinity since $\lim_{p\rightarrow \infty} h(p)/\log p =+\infty$. Combining  Theorem~\ref{t:ssadd}, \eqref{e:sens11} and \eqref{e:sens13} we conclude that $\nu_G(\v) \leq \wnu_{G,F,K_0} (\v)$.

We now prove that $\nu_G(\v) \geq \wnu_{G,F,K_0} (\v)$ for a fixed rational $\v \in \SS^{d-1}$. We associate with a fixed rational $\v \in \SS^{d-1}$ the same hyperrectangle $A$ as in the proof of Theorem~\ref{t:ssadd}, $A= \{ \sum_{i=1}^{d-1}\lambda_i\vec{f}_i\;:\;\forall i, \;\lambda_i\in [0,1] \}$. Let $(\w_1, \dots , \w_{d-1}) = (\vec{f}_1 / \|\vec{f}_1\|_2 , \dots , \vec{f}_{d-1} / \|\vec{f}_{d-1}\|_2)$: it is an orthonormal basis of the orthogonal complement of $\v$ made of rational vectors. 
We want to construct a set of edges that cuts $pA$ from $\hyp(pA+H_{F,K_0}(pA) \v)$ in $\slab(pA,H_{F,K_0} (pA),\v)$ by gluing together cutsets from the top to the bottom of different cylinders. For any fixed $\eta >0$, we slightly enlarge the hyperrectangle $A$ by considering 
$$ A^\eta \,=\, \{ x+\vec w \,:\, x\in A \,,\, \|\vec w \|_{\infty} \leq \eta \,,\, \vec w \cdot \v =0 \} \,.$$
Let $h(p) = \H^{d-1} (pA) ^{\frac{1}{2(d-1)}}  $ as previously. We consider the cylinder $\cyl^{\v}(pA^\eta , h(p))$, and a minimal cutset $E_0 (p,A,\eta )$ between the top $B_1^{\v}(pA^\eta , h(p))$ and the bottom $B_2^{\v} (pA^\eta , h(p))$ of this cylinder. To obtain a set of edges that cuts $pA$ from $\hyp (pA + H_{F,K_0} (pA) \v)$ in $\slab (pA, H_{F,K_0} (pA), \v)$, we need to add to $E_0 (p,A,\eta)$ some edges that prevent some flow to escape from $\cyl^{\v} (pA^\eta , h(p))$ by its vertical sides, see Figure \ref{f:sub2}.
\begin{figure}[h!]
\centering
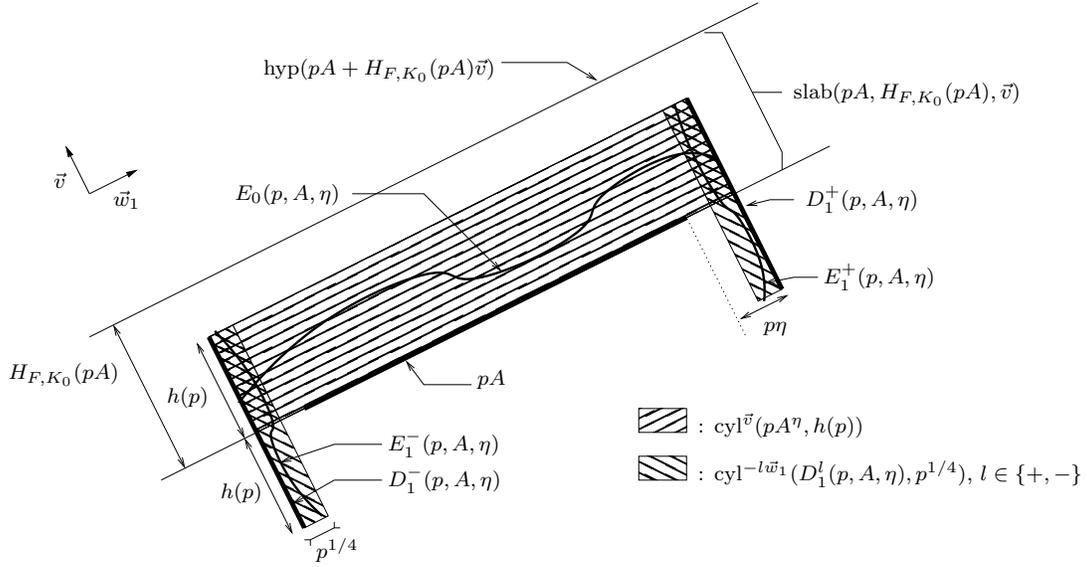
\caption{The construction of a cutset that separates $pA$ from $\hyp (pA + H_{F,K_0} (pA) \v)$ in $\slab (pA, H_{F,K_0} (pA), \v)$ (here $d=2$ and the cutsets $E_0 (p,A,\eta ), E^+_1 (p,A,\eta)$ and $E^+_1 (p,A,\eta)$ are represented {\em via} their dual as surfaces).}
\label{f:sub2}
\end{figure}
For $i\in \{1,\dots , d\}$, let $D^+_i (p,A,\eta)$ and $D^-_i (p,A,\eta)$ be the two $d-1$ dimensional sides of $\partial (\cyl(p A^\eta, h(p)))$ that are normal to $\w_i$, and such that $D^+_i (p,A,\eta)$ is the translated of $D^-_i (p,A,\eta)$ by a translation of vector $f (i,p,A,\eta) \w_i$ for some $f(i,p,A,\eta) >0$. We consider the cylinder $\cyl^{-\w_i} (D^+_i (p,A,\eta), p^{1/4} )$ (resp. $\cyl^{+\w_i} (D^-_i (p,A,\eta), p^{1/4} )$) and a minimal cutset $E^+_i (p,A,\eta)$ (resp. $E^-_i (p,A,\eta)$) from the top to the bottom of $\cyl^{-\w_i} (D^+_i (p,A,\eta), p^{1/4} )$ (resp. from the top to the bottom of $\cyl^{+\w_i} (D^-_i (p,A,\eta), p^{1/4} )$) in the direction $\w_i$. We emphasize the fact that the lengths of the sides of $\cyl^{-\w_i} (D^+_i (p,A,\eta), p^{1/4} )$ and $\cyl^{+\w_i} (D^-_i (p,A,\eta), p^{1/4} )$ do no grow to infinity at the same rate in $p$. We shall prove the three following properties:
\begin{itemize} 
\item[$(i)$] For every $\eta >0$, at least for $p$ large enough, the set of edges $F(p,A,\eta)$ defined by 
$$ F(p,A,\eta) \,:=\,  E_0 (p,A,\eta ) \cup \left(\bigcup_{i=1}^{d-1} E^+_i (p,A,\eta) \right)  \cup \left(\bigcup_{i=1}^{d-1} E^-_i (p,A,\eta) \right)$$
cuts $pA$ from $\hyp (pA+ H_{F,K_0} (pA) \v)$ in $\slab (pA, H_{F,K_0} (pA), \v)$.
\item[$(ii)$] For every $\eta >0$, 
$$ \lim_{p\rightarrow \infty} \frac{\phi^{\v} (p A^\eta , h(p) )}{\H^{d-1} ( p A^\eta)} \,=\,  \lim_{p\rightarrow \infty} \frac{T_G (E_0 (p,A,\eta ))}{\H^{d-1} ( p A^\eta)}  \,=\, \nu_G(\v)  \quad \textrm{in probability}\,.$$
\item[$(iii)$] For all $i\in \{1,\dots , d-1\}$, for all $l\in \{+,-\}$, for all $\eta>0$, we have 
$$ \lim_{p\rightarrow \infty} \frac{\phi^{-l \w_i} (D^l_i (p,A,\eta), p^{1/4} )}{p^{d-1}} \,=\,  \lim_{p\rightarrow \infty} \frac{ T_G (E^l_i (p,A,\eta)) }{p^{d-1}} \,=\, 0  \quad \textrm{in probability}\,.$$
\end{itemize}
Before proving these three properties, we show how they help us to conclude the proof. By property $(i)$ we know that for every $\eta>0$, for $p$ large enough,
\begin{equation}
\label{e:sens21}
\frac{\wphi_{G,F,K_0} (pA)}{\H^{d-1} (pA)} \,\leq \,  \frac{\H^{d-1} (A^\eta)}{ \H^{d-1} (A)} \frac{T_G (E_0 (p,A,\eta ))}{\H^{d-1} ( p A^{\eta})} +\frac{1}{\H^{d-1} (A)} \sum_{i=1}^{d-1} \sum_{l=+,-}  \frac{ T_G (E^l_i (p,A,\eta)) }{p^{d-1} } \,. 
\end{equation}
By Theorem~\ref{t:ssadd}, we know that the left hand side of \eqref{e:sens21} converges a.s. to $\wnu_{G,F,K_0}(\v)$ when $p$ goes to infinity. By properties $(ii)$ and $(iii)$, we know that the right hand side of \eqref{e:sens21} converges in probability to $\nu_G(\v) \H^{d-1} (A^\eta)/ \H^{d-1} (A)$, that is arbitrarily close to $\nu_G(\v)$ when $\eta$ goes to $0$. We conclude that $\wnu_{G,F,K_0}(\v) \leq \nu_G(\v)$.

To conclude the proof of Proposition~\ref{prop:ssadd} it remains to prove the properties $(i)$, $(ii)$ and $(iii)$. The proof of property $(i)$ is very similar to the proof of inequality \eqref{e:comp2} so we just recall the underlying idea of this proof without giving every details again. Indeed, if $\gamma$ is a path from $pA$ to $\hyp (pA+ H_{F,K_0} (pA) \v)$ in $\slab (pA, H_{F,K_0} (pA), \v)$, then $\gamma$ starts at a vertex of $V(pA)$, its next vertex is inside the cylinder $\cyl^{\v} (p A^{\eta/2} , h(p)/2) $ (for a fixed $\eta>0$, at least for $p$ large enough) and then after a finite number of steps it has to leave the cylinder $\cyl^{\v} (p A^{\eta} , h(p)) $ by its top or by one of its vertical faces. If it leaves by the top of this cylinder it must contain an edge of $E_0 (p,A, \eta)$, and if it leaves by one of its vertical faces it must contain an edge of one of the $E^l_i (p,A,\eta)$.

Property $(ii)$ is a straightforward application of Proposition \ref{p:positif} or \ref{p:zerobis}.

Property $(iii)$ is a bit more delicate to prove since we cannot apply Proposition \ref{p:positif} or \ref{p:zerobis} here. Indeed, for given $i\in \{1,\dots, d-1\}$ and $l\in\{+,-\}$, the base of the cylinder $\cyl^{-l \w_i} (D^l_i (p,A,\eta), p^{1/4} )$, namely $D^l_i (p,A,\eta)$, grows at speed $p$ in $(d-2)$ directions and at speed $h(p)$ (of order $ p^{\frac{1}{2(d-1)}}$) in one direction. We did not take into account this kind of anisotropic growth in our study (contrary to Kesten in \cite{Kesten:StFlour} and Zhang in \cite{Zhang,Zhang2017}). We can conjecture that $ \phi^{-l \w_i} (D^l_i (p,A,\eta), p^{1/4} )$ grows linearly with $p^{d-2} h(p)$, with a multiplicative constant given precisely by $\nu_G(\w_i)$, but this cannot be deduced easily from what has already been proved. However we do not need such a precise result. We recall that the definition of the event $\E_{G,K_0}(\cdot,\cdot)$ was given in \eqref{e:E}. Mimicking the proof of inequality \eqref{e:hop3} in the proof of Proposition \ref{p:finitude}, we obtain that for a constant $K(A,d,\eta)$, for variables $(X_i)$ that are i.i.d. with the same distribution as $\carde (C_{G,K_0} (0))$, we have
\begin{align}
\label{e:pfff}
\PP [ & \phi^{-l \w_i}  (D^l_i (p,A,\eta), p^{1/4} ) \geq \beta  p^{d-2} h(p) ]\nonumber\\
 & \,\leq\, \PP \left[\E_{G,K_0} \left( \cyl^{-l \w_i} (D^l_i (p,A,\eta), p^{1/4} ), \frac{p^{1/4}}{2} \right) ^c\right] + \PP \left[\sum_{i=1}^{K(A,d,\eta) \lfloor p^{d-2} h(p) \rfloor } X_i \geq c_dK_0^{-1} \beta p^{d-2} h(p) \right] \nonumber\\
& \,\leq\, K(A,d,\eta) p^{d-2} h(p) p^{1/4} \kappa_1 e^{-\kappa_2 p^{1/4} /2} + \EE[\exp(\lambda X_1)]^{K(A,d,\eta) p^{d-2} h(p)} e^{-\lambda \beta c_dK_0^{-1} p^{d-2} h(p)}
\end{align}
where $\lambda (G,d)>0$ satisfies $\EE[\exp(\lambda X_1)]<\infty$. Since $h(p)$ is of order $p^{\frac{1}{2(d-1)}}$, the first term of the right hand side of \eqref{e:pfff} vanishes when $p$ goes to infinity. We can choose $\beta (G,d,A)$ large enough such that the second term of the right hand side of \eqref{e:pfff} vanishes too when $p$ goes to infinity. This is enough to conclude that property $(iii)$ holds.
\end{dem}


\section{Continuity of $G \mapsto \nu_G $}
\label{s:cont}

This section is devoted to the proof of Theorem \ref{thmcont}. To prove this theorem we mimick the proof of the corresponding property for the time constant, see \cite{Cox}, \cite{CoxKesten}, \cite{Kesten:StFlour} and \cite{GaretMarchandProcacciaTheret}. We stress the fact that the proof relies heavily on these facts:
\begin{itemize}
\item[$(i)$] $\nu_G (\v) $ can be seen without any moment condition as the limit of a subadditive process that has good properties of monotonicity,
\item[$(ii)$] $ \lim_{K\rightarrow \infty} \nu_{G^K} (\v) = \nu_G(\v)$.
\end{itemize}
We stated Theorem \ref{t:ssadd} to get $(i)$. Property $(ii)$ is a direct consequence of the definition \eqref{defnu} of $\nu_G$ itself, but we had consequently to work to prove that the constant $\nu_G$ defined this way is indeed the limit of some rescaled flows. As a consequence, the following proof of Theorem \ref{thmcont} is quite classical and easy, since we already have in hand all the appropriate results to perform it efficiently.


\subsection{Preliminary lemmas}

Let $G$ (resp. $G_n,n\in \NN$) be a probability measure on $[0,+\infty]$ such that $G(\{+\infty\})<p_c(d)$ (resp. $G_n(\{+\infty\})<p_c(d)$). We define the function $\G : t \in [0,+\infty[ \mapsto G([t,+\infty])$ (respectively $\G_n(t) = G_n ([t,+\infty]) $) that characterizes $G$ (resp. $G_n$). Notice that the stochastic domination between probability measures can be easily characterized with these functions :
$$ G_1 \succeq G_2 \quad \Leftrightarrow \quad \forall t\in[0,+\infty[ \,, \,\, \G_1 (t) \geq \G_2(t) \,. $$
We recall that we always build the capacities of the edges for different laws by coupling, using a family of i.i.d. random variables with uniform law on $]0,1[$ and the pseudo-inverse of the distribution function of these laws. Thanks to this coupling, we get this classical result of convergence (see for instance Lemma 2.10 in \cite{GaretMarchandProcacciaTheret}).

\begin{lem}
\label{lemcouplage}
Let $G$, $(G_n)_{n\in \NN}$ be probability measures on $[0,+\infty]$. We define the capacities $t_G(e)$ and $t_{G_n} (e)$ of each edge $e\in \EE^d$ by coupling. If $G_n \overset{d}{\rightarrow} G$ then
$$ a.s.\,, \,\, \forall e \in \EE^d \,, \quad \lim_{n\rightarrow \infty} t_{G_n} (e) \,=\, t_G (e) \,. $$
\end{lem}

By the coupling relation \eqref{e:couplage}, Theorem \ref{t:oldCV} and the definition \eqref{defnu} of $\nu_G$, we also get trivially this monotonicity result.
\begin{lem}
\label{lemmonotonie}
Let $G$, $F$ be probability measures on $[0, +\infty]$ such that $G(\{+\infty\})<p_c(d)$ and $F(\{+\infty\})<p_c(d)$. If $G \preceq F$, then for all $\v \in \SS^{d-1}$ we have $\nu_G(\v) \leq \nu_{F} (\v)$.
\end{lem}

In what follows, it will be useful to be able to exhibit a probability measure that dominates stochastically (or is stochastically dominated by) any probability $G_n$ of a convergent sequence of probability measures, thus we recall this known result (see for instance Lemma 5.3 in \cite{GaretMarchandProcacciaTheret}) .
\begin{lem}
\label{lemsuiteprobas}
Suppose that $G$ and $(G_n)_{n\in \NN}$ are probability measures on $[0,+\infty]$ such that $G (\{+\infty\}) <p_c(d)$, $G_n (\{+\infty\}) <p_c(d)$ for every $n$ and $G_n \overset{d}{\rightarrow} G$. There exists a probability measure $F^+$ on $[0,+\infty]$ such that $F(\{+\infty\})<p_c(d)$, $G_n \preceq F^+$ for all $n\in \NN$ and $G \preceq F^+$.
\end{lem}


\subsection{Upper bound}

This is the easy part of the proof. It relies on the expression of $\nu_G (\v)$ as the infimum of a sequence of expectations.

\begin{prop}
\label{propupper}
Suppose that $G$ and $(G_n)_{n\in \NN}$ are probability measures on $[0,+\infty]$ such that $G(\{+\infty\}) < p_c(d)$ and for all $n\in \NN$, $G_n (\{+\infty\}) < p_c(d)$. If $G_n \succeq G$ for all $n\in \NN$ and $G_n \overset{d}{\rightarrow} G$, then for any rational direction $\v \in \SS^{d-1}$,
$$\limsup_{n\rightarrow \infty} \nu_{G_n} (\vec v) \,\leq \, \nu_G (\vec v)\,.$$
\end{prop}

\begin{dem}
Let $\v \in \SS^{d-1}$ be a rational vector, and let $A$ be the non-degenerate hyperrectangle normal to $\v$ given by Theorem~\ref{t:ssadd}. By Lemma \ref{lemsuiteprobas} we know that there exists a probability measure $F^+$ on $[0,+\infty]$ such that $F^+(\{+\infty\})<p_c(d)$, $G_n \preceq F^+$ for all $n\in \NN$ and $G \preceq F^+$. Let $K_0$ be large enough such that $F^+ (]K_0,+\infty])< p_c(d)$. Let $k\in \NN^*$. We recall that the definition of $H_{F^+,K_0} (kA)$ is given in \eqref{e:defH}. Let $E_k$ be a set of edges that cuts $kA$ from $\hyp (kA + H_{F^+,K_0} (kA) \v )$ in $\slab (kA,H_{F^+,K_0} (kA), \v )$. By coupling (see Equation \eqref{e:couplage}) we know that $\wphi_{G, F^+,K_0} (k A) \leq \wphi_{G_n, F^+,K_0} (k A)$, and by Lemma \ref{lemcouplage} we have a.s.
\begin{align*}
T_{G} (E_k) & \,=\, \sum_{e\in E_k} t_{G} (e) \,=\, \lim_{n\rightarrow \infty} \sum_{e\in E_k} t_{G_n} (e)\\
& \,\geq \, \limsup_{n\rightarrow \infty} \wphi_{G_n, F^+,K_0} (k A) \,\geq \,   \liminf_{n\rightarrow \infty} \wphi_{G_n, F^+,K_0} (k A) \, \geq\, \wphi_{G, F^+,K_0} (k A) \,,
\end{align*}
thus by optimizing on $E_k$ we obtain that
$$ \forall k\in \NN^* \,,\quad a.s. \,, \qquad  \lim_{n\rightarrow \infty} \wphi_{G_n, F^+,K_0} (k A) \,=\, \wphi_{G, F^+,K_0} (k A) \,.$$
Moreover, for all $k\in \NN^*$, we have also by coupling (see Equation \eqref{e:couplage}) that $\wphi_{G_n, F^+,K_0} (k A) \leq \wphi_{F^+, F^+,K_0} (k A)$ which is integrable. The dominated convergence theorem implies that
\begin{equation}
\label{eqTCD}
   \forall k\in \NN^*  \,, \qquad  \lim_{n\rightarrow \infty} \EE\left[  \wphi_{G_n, F^+,K_0} (k A)\right] \,=\, \EE \left[  \wphi_{G, F^+,K_0} (k A)\right] \,.
\end{equation}
By Theorem \ref{t:ssadd} we know that $\nu_G (\v) = \lim_{k\rightarrow \infty}  \EE \left[ \wphi_{G, F^+,K_0} (k A) \right] / \H^{d-1} (kA)$, thus for all $\varepsilon >0$ there exists $k_0$ such that
$$  \nu_G (\v) \,\geq\,   \frac{\EE \left[ \wphi_{G, F^+,K_0} (k_0 A) \right] }{\H^{d-1} (k_0 A)} - \varepsilon \,.$$
Using \eqref{eqTCD} we know that there exists $n_0$ such that for all $n\geq n_0$ we have
$$ \frac{\EE \left[  \wphi_{G, F^+,K_0} (k_0 A) \right] }{\H^{d-1} (k_0 A)} \,\geq \, \frac{\EE \left[ \wphi_{G_n, F^+,K_0} (k_0 A)  \right] }{\H^{d-1} (k_0 A)} - \varepsilon \,. $$
By Theorem \ref{t:ssadd} we also know that $\nu_{G_n} (\v) = \inf_{k}  \EE \left[\wphi_{G_n, F^+,K_0} (k A)  \right] / \H^{d-1} (k A)$, thus we obtain that for any $\varepsilon >0$, for all $n$ large enough,
$$ \nu_{G} (\v) \,\geq \, \nu_{G_n} (\v) - 2 \varepsilon \,.$$
This concludes the proof of Proposition \ref{propupper}.
\end{dem}

\subsection{The compact case}

This section is devoted to the proof of the continuity of the flow constant in the particular case where all the probability measures we consider have the same compact support $[0,R]$ for some finite $R$.

\begin{prop}
\label{propcompact}
Suppose that $G$ and $(G_n)_{n\in \NN}$ are probability measures on $[0,R]$ for some fixed $R\in [0,+\infty[$. If $G_n \overset{d}{\rightarrow} G$, then for any rational $\v \in \SS^{d-1}$,
$$ \lim_{n\rightarrow \infty} \nu_{G_n} (\vec v) \,= \, \nu_G (\vec v)\,.$$
\end{prop}

Notice that since the probability measure $G$ (resp. $G_n$) we consider has compact support, the flow constant $\nu_G (\vec{v})$ (resp. $\nu_{G_n} (\vec v)$) is defined via Theorem \ref{t:oldCV} as the limit when $k$ goes to infinity of the rescaled maximal flows $\tau_G (kA,k)$ (resp. $\tau_{G_n} (kA,k)$) defined in \eqref{e:deftau2} for a hyperrectangle $A$ normal to $\vec v$.
To prove Proposition \ref{propcompact}, we will follow the proof sketched by Kesten in the proof of the continuity of the time constant in \cite{Kesten:StFlour} to avoid the use of Cox's previous work \cite{Cox}. 
Let us describe briefly the proof of the main difficulty, namely that $\liminf_{n\rightarrow \infty} \nu_{G_n} (\v)  \geq \nu_G (\v)$. First, one reduces easily to the case where $G_n$ is stochastically dominated by $G$ for any $n$. Then, $0\leq \tau_{G_n} (kA,k) \leq \tau_{G} (kA,k)$ and 
$$\tau_{G} (kA,k)-\tau_{G_n} (kA,k)\leq \sum_{e\in E_n(k)}t_G(e)-t_{G_n}(e)\;,$$
where $E_n(k)$ is any minimal cutset for $\tau_{G_n} (kA,k)$. If one is able to control the size of $E_n(k)$, showing that the probability that it exceeds $\beta k^{d-1}$ decreases exponentially fast in $k^{d-1}$ for some $k$, then a standard union bound and a large deviation argument, with the fact that $\EE[\exp(t_G(e)-t_{G_n}(e))]$ is close to $1$ for $n$ large enough, will show that $\sum_{e\in E_n(k)}t_G(e)-t_{G_n}(e)$ is less than $\eps k^{d-1}$ uniformly in $n$ large enough. Thus we need a control on the size of a minimal cutset for $\tau_{ G_n} (kA,k)$ in the spirit of Theorem~1 in \cite{Zhang2017},  but uniformly in $n$ large enough. This was done in Proposition~4.2 in \cite{RossignolTheret08b}, but this proposition requires the sequence of distribution functions of $(G_n)$ to coincide on a neighborhood of $0$, at least for $n$ large enough, and this does not follow from our assumptions. Fortunately, inspecting the proof of Theorem~1 in \cite{Zhang2017} one may see that the conclusion of  Proposition~4.2 in \cite{RossignolTheret08b} holds uniformly in $n$ under weaker assumtions, stated in the following lemma. 
\begin{lem}
\label{l:zhang2017bis}
Suppose that $G_n$ is a sequence of distribution functions on $[0,+\infty[$ such that
\begin{equation}
  \label{eq:hypZhang1}
  \limsup_{n\rightarrow +\infty} G_n(]0,\eps[)\xrightarrow[\eps\rightarrow 0]{}0
\end{equation}
and
\begin{equation}
  \label{eq:hypZhang2}\limsup_{n\rightarrow +\infty} G_n(\{0\}){}<1-p_c(d)\;.
\end{equation}
For any $k$, let us denote by $E_n (k)$ a minimal cutset for $\tau_{\underline G_n} (kA,k)$ -- if there is more than one such cutset we choose (with a deterministic rule) one of those cutsets with minimal cardinality. Then, there exists positive constants $C$, $D_1$, $D_2$ and an integer $n_0$  such that
\begin{align}
\forall n\geq n_0\,, \forall \beta >0, \forall k\in\NN,\; \PP\left[\begin{array}{c}\carde (E_n(k)) \geq \beta k^{d-1} \\ \text{ and } \tau_{G_n} (kA,k) \leq \beta Ck^{d-1} \end{array}\right] \,\leq\, D_1 e^{-D_2 k^{d-1}}\,.\nonumber
\end{align}
\end{lem}

{\textbf{Proof of Proposition~\ref{propcompact}}~:\\}
Let $\v \in \SS^{d-1}$ be a rational direction. Let $G$, $(G_n)_{n\in \NN}$ be probability measures on $[0,R]$ for some $R \in [0,+\infty[$. We define $\underline \G _n = \min (\G, \G_n)$ (resp. $\overline \G _n = \max (\G, \G_n)$), and we denote by $\underline G_n$ (resp. $\overline G_n$) the corresponding probability measure on $[0,R]$. Then $\underline G_n \preceq G \preceq \overline G_n $ and $\underline G_n \preceq G_n \preceq \overline G_n $ for all $n\in \NN$, $\underline G_n \overset{d}{\rightarrow} G$ and $\overline G_n \overset{d}{\rightarrow} G$. To conclude that $\lim_{n\rightarrow \infty} \nu_{G_n} (\v) = \nu_G (\v)$, it is thus sufficient to prove that
\begin{itemize}
\item[$(i)$] $\limsup_{n\rightarrow \infty} \nu_{\overline G_n} (\v)  \leq \nu_G (\v)$, and
\item[$(ii)$] $\liminf_{n\rightarrow \infty} \nu_{\underline G_n} (\v)  \geq \nu_G (\v)$.
\end{itemize}
Inequality $(i)$ is a straightforward consequence of Proposition \ref{propupper}. If $\nu_G (\v) =0$, then inequality $(ii)$ is trivial and we can conclude the proof. From now on we suppose that $\nu_G (\v) >0$. By \cite{Zhang} (see Theorem \ref{t:oldnul} above), we know that $\nu_G (\v) >0 \iff G(\{0\}) < 1-p_c(d)$. Thanks to the coupling (see equation \eqref{e:couplage}), we know that for every edge $e\in \EE^d$, we have $ t_{\underline G_n} (e) \leq t_G (e)$.

Let $A$ be a non-degenerate hyperrectangle normal to $\v$ that contains the origin of the graph and such that $\H^{d-1} (A) = 1$. We recall that $\tau_G (kA, k)$ is defined in Equation \eqref{e:deftau2}. It denotes the maximal flow for the capacities $(t_G(e))$ from the upper half part $C_1'(kA, k)$ to the lower half part $C_2'(kA, k)$ of the boundary of $\cyl (kA, k)$ as defined in Equation \eqref{e:deftau1}, and it is equal to the minimal $G$-capacity of a set of edges that cuts the upper half part from the lower half part of the boundary of $\cyl (kA, k)$ in this cylinder. Since we work with integrable probability measures $G$ and $G_n$, we know by Theorem \ref{t:oldCV} that a.s. 
\begin{equation}
\label{e:hop}
\nu_G (\v) \,=\, \lim_{k\rightarrow \infty} \frac{\tau_G (kA, k)}{k^{d-1}} \quad \textrm{and} \quad \nu_{\underline G_n} (\v)\,=\,  \lim_{k\rightarrow \infty} \frac{\tau_{\underline G_n} (kA, k)}{k^{d-1}}  \,.
\end{equation}
%

Now, let us denote by $E_n (k)$ a minimal cutset for $\tau_{\underline G_n} (kA,k)$ as in Lemma~\ref{l:zhang2017bis}. According to Kesten's Lemma~3.17 in \cite{Kesten:StFlour}, any such minimal cutset with minimal cardinality is associated with a set of plaquettes which is a connected subset of $\RR^d$ - we will say that $E_n(k)$ is $\circ$-connected. Let $x\in \partial A$. There exists a constant $\hat c_d$, depending only on the dimension, such that for every $k\in \NN^*$ there exists a path from the upper half part to the lower half part of the boundary of $\cyl (kA, k)$ that lies in the Euclidean ball of center $kx$ and radius $\hat c_d$. We denote by $F(k)$ the set of the edges of $\EE^d$ whose both endpoints belong to this ball. Then $E_n (k)$ must contain at least one edge of $F(k)$, and $\carde (F(k)) \leq \hat c'_d$ for some constant $\hat c'_d$. Moreover, given a fixed edge $e_0$, the number of $\circ$-connected sets of $m$ edges containing $e_0$ is bounded by $\tilde c_d^m$ for some finite constant $\tilde c_d$ (see the proof of Lemma~2.1 in \cite{Kesten:StFlour}, that uses (5.22) in \cite{kesten:perco}). Thus for every $k \in \NN^*$, for every $\beta,C, \varepsilon >0$ we have
\begin{align}
\label{eqcompact}
\PP [& \tau_{\underline G_n} (kA,k)   \leq \tau_{G} (kA,k) - \varepsilon k^{d-1} ] \nonumber \\
& \,\leq \, \PP [\tau_{\underline G_n} (kA,k) > \beta C k^{d-1}]
+ \PP[ \carde (E_n (k) ) \geq \beta k^{d-1} \textrm{ and } \tau_{\underline G_n} (kA,k) \leq \beta Ck^{d-1} ]\nonumber \\
& \quad + \sum_{{\tiny \begin{array}{c}E  \textrm{ $\circ$-connected set of edges containing}\\ \textrm{some edge in } F(k) \textrm{ s.t. } \carde (E) \leq \beta k^{d-1} \end{array}}} \PP \left[ \sum_{e\in E} t_G(e) - t_{\underline G_n} (e) \geq \varepsilon k^{d-1}   \right]\nonumber \\
& \,\leq \, \PP [\tau_{G} (kA,k) > \beta C k^{d-1}]
+ \PP[ \carde (E_n (k)) \geq \beta k^{d-1} \textrm{ and } \tau_{\underline G_n} (kA,k) \leq \beta Ck^{d-1} ]\nonumber \\
& \quad +  c_d^{\beta k^{d-1}} \PP \left[ \sum_{i=1}^{\lfloor \beta  k^{d-1} \rfloor} t_G(e_i) - t_{\underline G_n} (e_i) \geq \varepsilon k^{d-1}   \right]\,,
\end{align}
where $(e_i)_{i\ge 1}$ is a collection of distinct edges and $c_d$ is a constant depending only on $d$. 

Let us prove that the sequence $(\underline G_n)$ satisfies conditions \eqref{eq:hypZhang1} and \eqref{eq:hypZhang2} of Lemma~\ref{l:zhang2017bis}. On one hand, since $\underline G_n \overset{d}{\rightarrow} G$ we have $\limsup_{n\rightarrow \infty} \underline G_n (\{0\}) \leq  G(\{0\})  < 1-p_c(d) $, thus condition  \eqref{eq:hypZhang2} is satisfied. On the other hand, we know that $\underline G_n  \preceq G$, {\em i.e.}, $\underline \G_n \leq \G$, thus for all $n\in \NN$ we have
$$\underline G_n (\{0\}) \,=\, 1-\lim_{p\rightarrow \infty} \underline \G_n (1/p) \,\geq\, 1-\lim_{p\rightarrow \infty}  \G (1/p) \,=\, G (\{0\}) \,,$$
and we conclude that
\begin{equation*}
\lim_{n\rightarrow \infty}\underline G_n (\{0\}) \,=\, G(\{0\}) \,.
\end{equation*}
Let $\eps\in ]0,+\infty[$ such that $G(\{\eps\}) =0$. Then $\underline G_n \overset{d}{\rightarrow} G$ implies that $\lim_{n\rightarrow \infty} \underline \G_n (\eps) =  \G(\eps) $, thus
$$ \underline G_n (]0,\eps[) \,=\,  1- \underline \G_n (\eps) - \underline G_n (\{0\})  \, \xrightarrow[n\rightarrow \infty]{}  1-  \G (\eps) -  G (\{0\})  \,=\, G (]0,\eps[)  \,,$$
and condition \eqref{eq:hypZhang1} follows from the fact that $\lim_{\eps \rightarrow 0}  G(]0,\eps[) = 0$. 
Now, we can use Lemma~\ref{l:zhang2017bis} to obtain the following uniform control:
\begin{align}
\label{eqzhang}
\exists  C, D_1, D_2,n_0  & \textrm{ such that } \forall n\geq n_0\,, \forall \beta , \forall k\nonumber \\
&  \PP[\carde (E_n (k)) \geq \beta k^{d-1} \textrm{ and } \tau_{\underline G_n} (kA,k) \leq \beta Ck^{d-1} ] \,\leq\, D_1 e^{-D_2 k^{d-1}}\,.
\end{align}
\medskip

We can easily bound $\tau_{G} (kA,k)$ by $\sum_{e\in E_k} t_G (e) \leq R \carde (E_k)$, where $E_k$ is a deterministic cutset of cardinality smaller than $c_d k^{d-1}$ - choose for instance $E_k$ as the set of all edges in $\cyl(kA, k)$ that are at Euclidean distance smaller than $2$ from $kA$. For any fixed $C>0$, since for every edge $e$ we have $t_{G} (e)\leq R$ there exists a constant $\beta $ such that
\begin{equation}
\label{eqborne}
\forall k\in \NN^* \,,\quad \PP [\tau_{G} (kA,k) > \beta C k^{d-1}] \,=\, 0 \,.
\end{equation}
By Markov's inequality, for any $\alpha >0 $ we have
$$ c_d^{\beta k^{d-1}} \PP \left[ \sum_{i=1}^{\lfloor \beta  k^{d-1} \rfloor} t_G(e_i) - t_{\underline G_n} (e_i) \geq \varepsilon k^{d-1}   \right] 
\,\leq \, \left( c_d \exp \left( \frac{-\alpha  \varepsilon}{\beta } \right)   \EE \left[ \exp \left( \alpha (t_{G} (e) - t_{\underline G_n} (e) ) \right) \right]\right)^{\beta  k^{d-1}}\,. $$
For any fixed $\varepsilon>0$ and $\beta <\infty$, we can choose $\alpha=\alpha (\varepsilon)$ large enough so that $c_d \exp (-\alpha  \varepsilon /\beta  ) \leq 1/4$. Then, by Lemma \ref{lemcouplage} we know that $\lim_{n\rightarrow \infty} t_{ \underline G_n}(e)  = t_G (e)$ a.s., and since $t_{G} (e) \leq R$ we can use the dominated convergence theorem to state that for $n$ large enough,
$$ \EE \left[ \exp \left( 2\alpha (t_{G} (e) - t_{ \underline G_n} (e) ) \right) \right] \,\leq\, 2\,.$$
We get
\begin{equation}
\label{eqsum1}
\sum_{k>0} c_d^{\beta k^{d-1}} \PP \left[ \sum_{i=1}^{\lfloor \beta  k^{d-1} \rfloor} t_G(e_i) - t_{\underline G_n} (e_i) \geq \varepsilon k^{d-1}   \right]  \,<\, +\infty\,.
\end{equation}
Combining \eqref{eqcompact}, \eqref{eqzhang}, \eqref{eqborne} and \eqref{eqsum1}, we obtain that for every $\varepsilon >0$, for all $n$ large enough,
$$ \sum_{k>0} \PP [  \tau_{\underline G_n} (kA,k)   \leq \tau_{G} (kA,k) - \varepsilon k^{d-1} ]\,<\,+\infty \,.$$
By Borel-Cantelli, we obtain that for every $\varepsilon >0$, for all $n$ large enough, a.s., for all $k$ large enough,
$$ \tau_{\underline G_n}(kA,k)  > \tau_{G}(kA,k) - \varepsilon k^{d-1}\,,$$
thus by \eqref{e:hop} for every $\varepsilon >0$, for all $n$ large enough, we have
$$ \nu_{\underline{G}_n} (\v) \,\geq \,\nu_G (\v) - \varepsilon \,. $$
This proves inequality $(ii)$ and ends the proof of Proposition \ref{propcompact}. 
{\flushright$\blacksquare$\\}


\subsection{Proof of Theorem \ref{thmcont}}

\begin{dem}
We first prove that convergence happens in a fixed rational direction $\v \in \SS^{d-1}$. We follow the structure of the proof of Proposition \ref{propcompact}. Let $G$, $(G_n)_{n\in \NN}$ be probability measures on $[0,+\infty[$. We want to prove that
\begin{equation}
\label{eqdirectionnelle}
\lim_{n\rightarrow \infty} \nu_{G_n} (\v) = \nu_G (\v)\,.
\end{equation}
We define $\underline G_n$ and $\overline G_n$ as in the proof of Proposition \ref{propcompact}, and we must show
\begin{itemize}
\item[$(i)$] $\limsup_{n\rightarrow \infty} \nu_{\overline G_n} (\v)  \leq \nu_G (\v)$, and
\item[$(ii)$] $\liminf_{n\rightarrow \infty} \nu_{\underline G_n} (\v)  \geq \nu_G (\v)$.
\end{itemize}
Inequality $(i)$ is still a straightforward consequence of Proposition \ref{propupper}. For every $K>0$, we define as previously $G^K = \ind{[0,K[} G + G([K,+\infty[) \delta_K$ and $\underline{G}_n^K = \ind{[0,K[} \underline{G}_n + \underline{G}_n([K,+\infty[) \delta_K$. Since $\underline{G}_n^K \preceq \underline{G}_n$, we know by Lemma \ref{lemmonotonie} that $\nu_{\underline{G}_n^K} \leq \nu_{\underline{G}_n}$. For every $K>0$, since $\underline G_n^K \overset{d}{\rightarrow} G^K$, using Proposition \ref{propcompact} we obtain that
\begin{equation}
\label{eqtroncfinal}
\liminf_{n\rightarrow \infty} \nu_{\underline G_n} (\v) \,\geq \, \lim_{n\rightarrow \infty} \nu_{\underline G_n^K} (\v) \,=\,  \nu_{\underline G^K} (\v)\,.
\end{equation}
By the definition \eqref{defnu} of $\nu_G (\v)$ we know that $\lim_{K\rightarrow \infty} \nu_{ G^K} (\v) = \nu_{ G} (\v)$. This concludes the proof of $(ii)$, thus \eqref{eqdirectionnelle} is proved.

We consider the homogeneous extension of $\nu_G$ to $\RR^d$ defined in Proposition~\ref{p:cvx}. By Proposition~\ref{p:cvx}, for all $x, y\in \RR^d$, we have $\nu_G ( x) \leq \nu_G (x-y) + \nu_G(y)$ and $\nu_G ( y) \leq \nu_G (x-y) + \nu_G(x)$ thus
\begin{equation}
\label{eqpropnu1}
|\nu_G (x) - \nu_G (y)| \,\leq \, \nu_G (x-y)\,.
\end{equation}
Moreover for all $x=(x_1,\dots ,x_d)$, we have
\begin{align}
\label{eqpropnu2}
\nu_G (x) & \,\leq \, \nu_G ((x_1, 0, \dots , 0)) + |\nu_G ((x_1, x_2, 0 , \dots , 0)) -  \nu_G ((x_1, 0, \dots , 0)) | \nonumber \\
& \qquad + \dots + |\nu_G ((x_1, \dots , x_d)) - \nu_G ((x_1, \dots , x_{d-1}, 0))|\nonumber \\
& \,\leq \, \nu_G (( x_1, 0, \dots , 0 )) + \nu_G ((0, x_2, 0 , \dots , 0 )) +\dots + \nu_G ((0 , \dots , 0, x_d ))\nonumber \\
& \,\leq \, \|x\|_1 \nu_G((1, 0, \dots , 0)) \,.
\end{align}
Combining \eqref{eqpropnu1} and \eqref{eqpropnu2}, we obtain that for all $x,y\in \RR^d$,
\begin{equation}
\label{e:ajout}
|\nu_G (x) - \nu_G (y)| \,\leq \, \|x-y\|_1\nu_G ((1, 0, \dots , 0))\,.
\end{equation}
The same holds for $\nu_{G_n}$. Since $\lim_{n\rightarrow \infty}\nu_{G_n}((1, 0, \dots , 0)) = \nu_G ((1, 0, \dots , 0))$, there exists $n_0$ such that for all $n\geq n_0$, we have $\nu_{G_n}((1, 0, \dots , 0)) \leq 2 \nu_G ((1, 0, \dots , 0))$. For every $n\geq n_0$, we have
\begin{equation}
\label{eqpropnu3}
 \forall \vec u, \v \in \SS^{d-1}\,, \quad  |\nu_{G_n} (\vec u) - \nu_{G_n} (\v)| \,\leq \, 2 \|\vec u - \v \|_1\nu_G ((1, 0, \dots , 0))\,. 
\end{equation}
Fix $\varepsilon >0$. Inequalities \eqref{e:ajout} and \eqref{eqpropnu3} imply that there exists $\eta >0$ such that
$$ \sup \{ |\nu_F (\vec u) - \nu_F (\v)| \,:\, \vec u , \, \v \in \SS^{d-1} ,\,\, \|\vec u - \v \|_1 \leq \eta,\,\, F \in \{G,G_n, n\geq n_0\}\} \,\leq\, \varepsilon \,. $$
There exists a finite set $(\v_1, \dots , \v_m)$ of rational directions in $\SS^{d-1}$ such that
$$ \SS^{d-1} \,\subset \, \bigcup_{i=1}^m  \{\vec u \in \SS^{d-1} \,:\, \|\vec u - \v_i \|_1  \leq \eta \} \,.$$
Thus 
$$ \limsup_{n\rightarrow \infty} \sup_{\vec u \in \SS^{d-1}} | \nu_{G_n} (\vec u)  - \nu_G (\vec u) | \,\leq \, 2 \varepsilon + \lim_{n\rightarrow \infty} \max_{i \in \{1,\dots , m\} }  | \nu_{G_n} (\v_i)  - \nu_G (\v_i) | \,,  $$
and thanks to \eqref{eqdirectionnelle} this ends the proof of Theorem \ref{thmcont}.
\end{dem}

\paragraph{Acknowledgements.} The second author would like to thank Rapha\"el Cerf for stimulating discussions on this topic. The authors would like to thank an anonymous referee for many valuable comments that helped to improve the quality of the paper.


\bibliographystyle{plain}
\bibliography{biblio}

\end{document}